\documentclass[14pt,letterpaper,reqno]{article}
\usepackage{epsf}
\usepackage{amsmath}
\usepackage{amsfonts}
\usepackage{amssymb}
\usepackage{amstext}
\usepackage{amsbsy}
\usepackage{amsthm}
\usepackage{amscd}
\usepackage[all]{xy}

\newcommand{\m}{\mathcal}
\newcommand{\w}{\widetilde}

\newcommand{\Zee}{\mathbb{Z}}
\newcommand{\Cee}{\mathbb{C}}

\newtheorem{teo}{Theorem}
\newtheorem{prop}[teo]{Proposition}

\newtheorem{rem}[teo]{Remark}
\newtheorem{dfn}[teo]{Definition}

\newtheorem{lem}[teo]{Lemma}

\title{Mathematics Underlying the F-Theory/Heterotic String Duality in 
Eight Dimensions} 
\author{Adrian Clingher
\thanks{The first author was supported by NSF grants DMS-97-29992 and 
PHY-00-70928.} 
\and 
John W. Morgan\thanks{
The second author was partially supported by NSF grant DMS-01-03877.}
}

\begin{document}
\maketitle

\vline
%\input{intro}
%
%
%     INTRODUCTION
%
\section{Introduction} 
\label{intro}
One of the dualities in string theory, the F-theory/heterotic string 
duality in eight dimensions \cite{vafa}, predicts an
interesting correspondence between two seemingly disparate geometrical objects. 
On one side of the duality there are elliptically fibered $K3$ surfaces with section. 
On the other side, one finds elliptic curves endowed with certain flat connections
and complexified Kahler classes. 
\par The F-theory \cite{sen} \cite{vafa} is a $12$-dimensional string theory which 
generally exists on 
elliptically fibered ambient manifolds with section. Heterotic string theory, on the other hand, 
exists on a $10$-dimensional space-time. In order to obtain effective $8$-dimensional 
models, one compactifies the two theories along elliptic $K3$ surfaces and 
elliptic curves, respectively. The duality mentioned above predicts then that the two 
theories are equivalent at the quantum level. In particular, their moduli spaces 
of quantum vacua should be isomorphic. As it is generally believed, in certain ranges 
of parameters the quantum corrections should be small and the quantum vacua should be 
well approximated by classical vacua. This leads one to expect that, the moduli 
spaces of classical vacua of the two theories should
resemble each other, at least on regions corresponding to insignificant quantum 
effects. 
%\par This paper is part of a project aimed at establishing the rigorous mathematical 
%results describing the geometry underlying the classical aspects of this $8$-dimensional duality. 
%The task involves understanding and comparing the classical moduli spaces 
%on the two sides.  
\par The classical vacua for the heterotic string theory compactified along a 
two-torus $E$ (there are two distinct such theories, one with structure Lie group 
$G_1= \left ( E_8 \times E_8 \right )\rtimes \Zee_2 $ and the other with 
$ G_2={\rm Spin}(32) / \Zee_2 $) consists of a flat $G_i$-connection on $E$, 
a flat metric and an extra one-dimensional field, the B-field. In the original physics
formulation \cite{narain1} \cite{narain2}, the B-field appears as a globally defined 
two-form $B$. The metric and $B$ fit together to form the imaginary and, respectively, 
the real part of the so-called complexified Kahler class. Each triplet $(A,g,B) $ 
determines a lattice of momenta $ \m{L}_{(A,g,B)} $ (after K. Narain \cite{narain1}) 
governing the 
associated physical theory. The lattices $ \m{L}_{(A,g,B)} $, turn out to be even, 
unimodular and of rank $20$. They are well-defined up to $ O(2) \times O(18) $ rotations 
and vary, according to the triplet parameter, in a fixed ambient real space 
$ \mathbb{R}^{2,18} $. The real group $ O(2,18) $ acts transitively on the set of all 
$ \m{L}_{(A,g,B)} $ and, in 
this light, one can regard the physical momenta as parameterized by the $36$-dimensional
real homogeneous space:
\begin{equation}
\label{coset} 
O(2,18) / O(2) \times O(18). 
\end{equation}         
One identifies then the configurations in $ (\ref{coset}) $ determining equivalent quantum 
theories. This amounts to factoring out the left-action of the group $ \Gamma $ of integral 
isometries of the lattice. However, not all identifications
so created are accounted for by classical geometry. Part of the $\Gamma $-action models 
the so-called quantum corrections \cite{morrison1} and results in
identifying momenta for pairs of triplets $(A,g,B) $ which are not isomorphic from the geometric
point of view. The quantum (Narain) moduli space of distinct heterotic string 
theories compactified on the two-torus appears as:
\begin{equation}
\label{narainmod}  
\m{M}^{\rm quantum}_{{\rm het}} \ = \ \Gamma \backslash O(2,18) / O(2) \times O(18).
\end{equation} 
\par The above physics-inspired Narain construction has a major flaw, though. 
It does not provide a holomorphic
description. Technically, one can endow the homogeneous quotient $ (\ref{narainmod}) $ 
with a natural complex structure, but holomorphic families of elliptic curves and flat
connections do not embed as holomorphic sub-varieties in  $ \m{M}^{\rm quantum}_{{\rm het}} $.        
\par In the recent years, it has been noted by a number of
authors (see for example \cite{witten1} or \cite{freed}) 
that, in order to fulfill various anomaly cancellation 
conditions required by heterotic
string theory, the B-field has to be understood within a gerbe-like formalism. 
In \cite{freed}, D. Freed introduces $B$ as a cochain in differential cohomology. Taking this point of view, 
one can ask then for a description of the space of triplets $ (A,g,B) $ up to natural 
geometric isomorphism.
This is the moduli space $ \m{M}^{G_i}_{{\rm het}} $ of classical vacua for 
$G_i$-heterotic string theory compactified over the two-torus. Freed's approach can be used to 
describe $ \m{M}^{G_i}_{{\rm het}} $ in an explicit holomorphic framework. It was shown 
in \cite{clingher} that: 
\begin{teo} \ 
\label{het1}
\begin{enumerate}
\item The classical $G_i$-heterotic moduli space $ \m{M}^{G_i}_{{\rm het}} $ can be given
the structure of a 18-dimensional complex variety with orbifold singularities. 
\item $ \m{M}^{G_i}_{{\rm het}} $ represents the total space of a holomorphic Seifert 
$ \mathbb{C}^* $-fibration
\begin{equation}
\label{fib}
\m{M}^{G_i}_{{\rm het}} \rightarrow  \m{M}_{E, G_i},
\end{equation} 
where $ \m{M}_{E, G_i} $ is the moduli space of isomorphism classes of 
pairs of elliptic curves and flat $ G_i$-bundles. 
\end{enumerate}
\end{teo}    
\noindent The holomorphic orbifold structure of $ \m{M}_{E, G_i} $ is described in \cite{mf1}. 
If one denotes by $ \mathbb{H} $ the upper half-plane, and by $ \Lambda $ 
the co-root lattice of $G_i$, then $ \m{M}_{E, G_i} $ is represented by a quotient of 
$\mathbb{H} \times \Lambda_{\Cee} $ through the action of a discrete group. Under 
this description, the fibration $ (\ref{fib}) $ appears as a well-known fibration
with complex lines over $ \m{M}_{E, G_i} $. This is, roughly
speaking, the line fibration supporting the holomorphic theta function:                     
\begin{equation}
\label{thh}
B_{G} \colon \mathbb{H} \times \Lambda_{\Cee}  \rightarrow \mathbb{C}, \ \ \ \ 
B_{G}(\tau, \ z) = \ \frac{1}{\eta(\tau)^{16}} \ \left ( \sum_{\gamma \in \Lambda} \ 
e^{ \pi i ( 2 (z,\gamma) + \tau (\gamma, \gamma))} \right )   
\end{equation}
where $ \eta(\tau) $ 
is Dedekind's eta-function. In this setting, one can prove:
\begin{teo} (\cite{clingher})
\label{het2}  
The $ \Cee^* $-fibration $ (\ref{fib}) $ is holomorphically identified with 
the complement of the zero-section in the complex line fibration induced by $(\ref{thh}) $.  
\end{teo}
\noindent Hence, the heterotic classical moduli space $ \m{M}^{G_i}_{{\rm het}} $ can be 
holomorphically identified with the total space of the theta fibration with 
the zero-section divisor removed. 
\par We turn now to the other side of the duality. The classical vacua for $8$-dimensional 
F-theory are simply elliptically fibered $K3$ surfaces with section. 
Using the period map and global Torelli theorem \cite{bpv} \cite{frthesis}, one can regard 
the moduli space $ \m{M}_{K3} $ of such structures as a moduli space of Hodge structures of weight two, 
i.e. as a quotient of an open $18$-dimensional 
hermitian symmetric domain $ \Omega $ by an arithmetic group of integral automorphisms.
In order to identify all equivalent classical vacua, there is one more factorization to 
be taken into account, identifying the complex conjugate structures. The classical 
$8$-dimensional F-theory moduli space obtained, denoted $ \m{M}_{F} $, is then double-covered 
by $  \m{M}_{K3} $ and can be seen to be isomorphic to an arithmetic
quotient of a symmetric domain:
\begin{equation}
\m{M}_{F} \ \simeq \  \Gamma \backslash O(2,18) / O(2) \times O(18).
\end{equation} 
The identification between the above description and $ (\ref{narainmod} ) $ is the usual
physics literature formulation of the F-theory/heterotic string duality in 
eight dimensions.  
\par The goal of this paper is to establish a rigorous geometric comparison between
the classical moduli spaces $ \m{M}_F $ and $ \m{M}^{G_i}_{{\rm het}} $. Our construction provides 
a natural holomorphic identifications between these classical moduli spaces, exactly 
on the regions where physics predicts the quantum effects are insignificant.     
\par The paper is structured as follows. In section $ \ref{first} $ we review various facts 
pertaining to the construction
of the moduli space $\m{M}_{K3} $ of elliptic $K3$ surfaces with section. This space is not compact. 
However, using a special case of Mumford toroidal compactification 
\cite{ash} \cite{friedman}, one can 
perform an arithmetic partial compactification: 
$$ \m{M}_{K3} \ \subset \ \overline{\m{M}}_{K3} $$
by adding two divisors at infinity $ \m{D}_1 $ and $ \m{D}_2 $,  
related to the two possible kinds of Type II maximal 
parabolic subgroups of $ O(2,18) $. The arithmetic machinery 
producing the partial compactification is reviewed in section $ \ref{arith} $. 
In section $ \ref{stablek3} $
we discuss the geometrical interpretation of the compactification. The points of 
$ \m{D}_1 $ and $ \m{D}_2 $ correspond to semi-stable degenerations of $ K3$ 
surfaces given by either a union of two rational elliptic surfaces glued together along a smooth
fiber or a union of rational surfaces glued along an elliptic curve with elliptic
fibration degeneration into two rational curves meeting at two points. Each of the
two configurations exhibits an elliptic curve $E$ (the double curve of the degeneration) 
and endows this elliptic curve with a flat
$ G$-connection. For $ \m{D}_1 $ the resulting Lie group $G$ turns out to be 
$ G_1 = \left ( E_8 \times E_8 \right ) \rtimes \Zee_2 $ whereas 
for $ \m{D}_2 $ one obtains $ G_2 = {\rm Spin}(32)/ \Zee_2 $. In the second case the flat 
connections obtained carry ``vector structure'', in the sense
that they can be lifted to flat $ {\rm Spin}(32) $-connections. Under this geometrically 
defined correspondence, one obtains a holomorphic isomorphism:
\begin{equation}
\label{compizo} 
\m{D}_i \ \simeq \ \m{M}_{E,G_i}. 
\end{equation} 
\par Next, each of the two types of parabolic groups determining the boundary components 
$ \m{D}_i $ produces an infinite
sheeted non-normal parabolic cover $ p \colon \m{P}_i \rightarrow \m{M}_{K3} $. The 
total space $ \m{P}_i $ fibers holomorphically $ \pi \colon \m{P}_i \rightarrow  \m{D}_i $ 
over the corresponding divisor at infinity, all fibers being 
copies of $ \Cee^* $. Under identification $ (\ref{compizo}) $, one obtains
therefore a pattern:
\begin{equation}
\begin{array}{ccccc}
\m{P}_i & \stackrel{\pi}{\rightarrow} & \m{D}_i & \simeq & \m{M}_{E,G_i} \\
\downarrow & & & & \\
\m{M}_{K3}. & & & & \\
\end{array} 
\end{equation}
It turns out, a neighborhood of infinity near the cusp $ \m{D}_i $ in $ \m{M}_{K3} $ 
is identified with a component of its pre-image in the parabolic cover $ \m{P}_i $. 
Moreover, this pre-image component is a neighborhood of the zero-section in 
\begin{equation} 
\label{ppp} 
\pi \colon \m{P}_i \rightarrow  \m{M}_{E,G_i}.  
\end{equation} 
Thus, a neighborhood of the boundary component $ \m{D}_i $ in $ \m{M}_{K3} $ can then be 
identified with a neighborhood of the zero-section of 
the parabolic fibration $ (\ref{ppp}) $. 
\par In section $ \ref{main} $ we give an explicit description of $ (\ref{ppp}) $. 
Based on this description we conclude:
\begin{teo}
Fibration $ (\ref{ppp}) $ is holomorphically isomorphic with the theta $\Cee$-fibration
induced by $ (\ref{thh}) $ with the zero-section removed. 
\end{teo} 
\noindent In the light of theorems $ \ref{het1} $  and $ \ref{het2} $, there is then 
a holomorphic isomorphism of $ \Cee^* $-fibrations, unique up to twisting with a unitary 
complex number: 
\begin{equation}
\begin{array}{ccc}
\m{P}_i & \simeq & \m{M}_{{\rm het}}^{G_i} \\
\downarrow & & \downarrow \\
\m{M}_{E,G_i} &= & \m{M}_{E,G_i} \\ 
\end{array} 
\end{equation}  
and that gives a natural explicit mathematical identification between 
the region in 
$ \m{M}_{{\rm het}}^{G_i} $ corresponding to large volumes with a 
region of $ \m{M}_F $ in the vicinity of the boundary component $D_i $. 
These are
exactly the regions that the physics duality predicts should be isomorphic. 
\par This paper belongs to a long project begun by the second author
jointly with R. Friedman and E. Witten \cite{rmw} in 1996 and continued
jointly with R. Friedman afterwards. 
The initial aim of the project was to give precise mathematical 
descriptions of various moduli spaces of principal $G$-bundles over 
elliptic curves in order to verify conjectures arising out of 
the F-theory/heterotic string theory duality in physics. Building on this 
earlier work, the present paper and \cite{clingher} establish 
the mathematical results allowing one to describe the 
duality completely when the two theories in 
question are compactified to eight dimensions.
\par The authors would like to thank Robert Friedman for many helpful 
conversations during the development of this work. The first author would 
also like to thank Charles Doran for many discussions regarding this work and 
the Institute for Advanced Study for its hospitality
and financial support during the course of the academic year 2002-2003.   
%
%
%
%\input{s1}
%
%
%  SECTION 1
%
\section{Review of the Compactification Procedure}
\label{first}
A coarse moduli space $ \m{M}_{K3} $ for isomorphism classes of elliptically fibered
$K3$ surfaces with section can be described using the period map. In this
section we review the Type II partial compactification of $ \m{M}_{K3} $. 
\subsection{Period Space}
\label{per1}
It is well-known that any two $K3$ surfaces are diffeomorphic. 
The second cohomology group over integers is torsion-free of rank $22$ and, 
when endowed with the symmetric bilinear form given by cup product, is an 
even unimodular 
lattice of signature $(3,19) $. Up to isometry, there exists a unique 
lattice with these properties. We pick a lattice of this type and denote
it by $ \mathbb{L} $. It happens then that for any $K3$ surface $X$ there
always exists an isometry:
\begin{equation}
\label{mark}
\varphi \colon H^2(X, \Zee) \rightarrow  \mathbb{L}.  
\end{equation}
Such a map is called a marking. 
\par An elliptic structure with section on $X$ induces naturally two particular
line bundles $ \m{F},  \m{S} \in {\rm Pic}(X) $ corresponding to the elliptic
fiber and section. %In particular, it can be shown that $ \m{F}$ and $ \m{S} $
%determine a quasi-ample polarization of degree two. Indeed, let 
%$ \m{L} = \m{S} \otimes \m{F}^{\otimes 2} \in {\rm Pic}(X) $. One checks that
%$ \m{L}^2=2 $ and $ \m{L} \cdot C \geq 0 $ for all irreducible curves
%$C$ on $X$. In such conditions, a theorem of Mayer \cite{mayer} asserts that
%the linear system $ \vert k \m{L} \vert $  has no base locus or fixed components 
%on $X$ for $ k \geq 2 $ and, moreover, for $ k \geq 2 $, the induced map:
%$$ \varphi_{\vert k\m{L}\vert } \colon X \rightarrow \mathbb{P}^N $$
%has as image a surface $ \overline{X} $ obtained by contracting all curves 
%$C$ on $X$ satisfying $ C^2=-2$ and $ C \cdot \m{L} = 0 $. In particular,
%$\overline{X} $ may have, as singularities, at most rational double points.
%\par However, in our discussion we shall mostly ignore this particular
%degree-two polarization and we shall rather consider the lattice-polarization 
%determined by both $ \m{F}$ and $ \m{S} $. 
Let $ f, s  \in H^2(X, \Zee) $ 
be the cohomology classes corresponding $ \m{F}$ and $ \m{S} $.         
These special classes intersect as $ f^2 =0, \ f.s=1, \ s^2=-2 $ and therefore   
span a hyperbolic type sub-lattice $ Q $ inside $ H^2(X, \Zee) $. The notion of 
marking can be adapted for this framework. Let 
$ H $ be a choice of hyperbolic sub-lattice in $ \mathbb{L} $.
All such choices are equivalent under the action of the group of 
isometries of $ \mathbb{L} $. Choose a basis $ \{F,S\} $ for
$ H $ with $ (F,F) = 0 $ , $ (F,S) = 1 $ and $ (S, S) = -2 $.  
A marking $ \varphi $ as in $ (\ref{mark}) $ is said to be compatible
with the elliptic structure if $ \varphi (f)=F $ and $ \varphi (s) = S $. 
In particular, a compatible marking transports the hyperbolic 
sub-lattice $Q \subset H^2(X, \Zee) $ isomorphically to $ H $.
Two marked pairs $ (X, \varphi) $ and $ (X', \varphi') $ are called 
isomorphic if there exists an isomorphism of surfaces 
$ g \colon X \rightarrow X' $ such that $ \varphi' = \varphi \circ g^* $.  
\par Let $\mathbb{L}_o$ be the sub-lattice of $\mathbb{L}$ orthogonal to $H$.    
The lattice $ \mathbb{L}_o $ is even, unimodular, and of signature $(2,18) $. 
By standard arguments, a marked pair $ (X, \varphi) $ determines 
a polarized Hodge structure of weight two on $ \mathbb{L}_o \otimes \Cee$ which 
is esentially determined by the period (2,0)-line  
$[\omega] \subset \mathbb{L}_o \otimes \Cee $. The periods satisfy 
the Hodge-Riemann bilinear relations $ (\omega , \omega) = 0, \  (\omega,\bar{\omega}) > 0 $. 
The classifying space of polarized Hodge structures of weight two on $ \mathbb{L}_o \otimes \Cee$
is then given by the period domain 
\begin{equation}
\Omega \ = \ \{ \ \omega \in \mathbb{P} \left (\mathbb{L}_o \otimes _{\Zee} \Cee 
  \right )    
\ \vert \ (\omega,\omega)=0, \ (\omega,\bar{\omega}) > 0 \  \}. 
\end{equation}
This is an open $18$-dimensional 
complex analytic variety embedded inside the compact complex quadric: 
$$ \Omega^{\vee} = \{ \omega  \ \vert \ (\omega,\omega) =0 \} \subset 
\mathbb{P} \left ( \mathbb{L}_o \otimes _{\Zee} \Cee  \right ). $$ 
\noindent One can equivalently regard the periods $ \omega \in \Omega $ as 
space-like, oriented two-planes in $ \mathbb{L}_o \otimes \mathbb{R} $. 
The real Lie group $ O(2,18) $ of real isometries of 
$ \mathbb{L}_o \otimes \mathbb{R} $ acts then transitively on $ \Omega $
leading to a description of the period domain in the form of a symmetric
bounded domain:    
\begin{equation}
\label{s}
\Omega \ \simeq \   O(2,18) / SO(2) \times O(18) .  
\end{equation} 
\par Following arguments of \cite{shapiro} \cite{looijenga}, one can prove the existence 
of a fine moduli space of marked elliptically fibered 
$K3$ surfaces with section, which is a $18$-dimensional complex manifold $ \m{M}^{\rm mark}_{K3} $. 
It follows then that a marked elliptic $ K3$ surface with section 
is uniquely determined by its period. The period map:
\begin{equation}
{\rm per} \colon \m{M}^{\rm mark}_{K3} \rightarrow \Omega 
\end{equation}  
is a  holomorphic isomorphism. However, in this setting, the period $ \omega \in \Omega $ 
clearly depends on the choice of marking. One removes the markings from the picture 
by dividing out the period domain by the action of the isometry 
group of the lattice. 
\par Let $ \Gamma $ be the group of isometries of $ \mathbb{L}_o  $. 
Two periods correspond to isomorphic marked surfaces if and only if they can be 
transformed one into the other through an isometry in $ \Gamma $. The arguments of 
global Torelli theorem allow one to conclude that:
\begin{equation}
\label{mmmm} 
\m{M}_{K3} \ = \ \Gamma \backslash \Omega 
\end{equation} 
is a coarse moduli space for elliptic $K3$ surfaces with section, without regard
to marking.
%
%  REMOVED
%
% Moreover, the points in $ \Gamma \backslash \Omega $ which are not 
%in  $ \Gamma \backslash \m{M}^m $ can be seen \cite{morrison} to correspond 
%to singular $K3$ surfaces with rational double points\footnote{Such a surface is 
%a compact complex surface $X$ with rational double points, whose minimal resolution
%$ \w{X} \rightarrow X $ is a $K3$ surface.}. Therefore, allowing rational double
%points as singularities, one has:
%\begin{equation}
%\label{mmmm} 
%\m{M} \ = \ \Gamma \backslash \Omega 
%\end{equation}    
%as moduli space for elliptically fibered $K3$ surfaces with section.
%
\par Let us briefly analyze the quotient $(\ref{mmmm}) $. First of all, $ \m{M}_{K3} $ is
connected. The period domain $ \Omega $ consists of two connected components, 
corresponding to the choice of 
%a continuously varying 
orientation in the set of positive two-planes in 
$ \mathbb{L}_o \otimes \mathbb{R} $. The two components are mapped into each 
other by complex conjugation. We choose either one and denote it by 
$ D $. Thus $ \Omega = D \sqcup \overline{D} $. However, there are isometries in 
$ \Gamma $ which exchange $D$ and $\overline{D} $ and therefore  
$(\ref{mmmm}) $ is connected. Secondly, the space $\m{M}_{K3} $ can be given a description
as a quotient of a bounded symmetric domain by a discrete, arithmetically defined
modular group. Indeed, the isomorphism $ (\ref{s}) $ is $\Gamma $-equivariant and
therefore:
$$ 
\label{s1}
\Gamma \backslash \Omega \ \simeq \   \Gamma \backslash O(2,18) / SO(2) \times O(18) .  
$$
%  
%
%  ADD THIS TO INTRODUCTION 
%
%The quotient variety $ \Gamma  \backslash  \Omega $ is far from being compact. However, 
%one can perform a partial compactification using essentially a 
%Mumford toroidal compactification procedure \cite{ash} \cite{friedman}. 
%This construction has a dual interpretation. Arithmetically, it amounts to  
%adding to the quotient $ \Gamma  \ \backslash  \Omega $ certain boundary 
%components $ \m{B}(P) $ associated with type II maximal rational parabolic subgroups 
%$P \subset O^{++}(2,18) $. From the Hodge theory point of view, the arithmetically 
%defined boundary components correspond to mixed Hodge structures extending the pure 
%Hodge structures associated to the periods in $ \Omega $. In this direction, the 
%construction incorporates deeper geometrical information. On the moduli side, 
%it amounts to enlarging the space $ \m{M} $ by allowing certain degenerations, 
%the stable $K3$ surfaces of Type II. In what follows, we summarize the 
%procedure and the main ideas involved, as they appear in \cite{ash} \cite{friedman} 
%\cite{frthesis}. 

\subsection{The Classical F-Theory Moduli Space} 
One obtains the moduli space $ \m{M}_F $ of classical vacua associated to F-theory 
compactified on a $K3$ surface by identifying conjugated complex structures in 
$ \m{M}_{K3} $. In the light of the previous discussion, one can assume then that:
\begin{equation}
\label{ss333}
\m{M}_{F} \ = \ \hat{\Gamma} \backslash \Omega
\end{equation}  
where $ \hat{\Gamma} $ is the semi-direct product 
$ \Gamma \rtimes \Zee_2  \subset  {\rm Aut} \left ( \mathbb{L}_o \otimes _{\Zee} \Cee \right ) $ 
with the $ \Zee_2 $ factor generated by complex conjugation.
%\begin{equation}
%\hat{\Gamma} \ = \ \Gamma \rtimes \Zee_2  \ \subset \ {\rm Aut} \left ( 
%\mathbb{L}_o \otimes _{\Zee} \Cee 
%\right )  
%\end{equation}
%with the $ \Zee_2 $ factor generated by complex conjugation.   
\par The moduli space $ (\ref{ss333}) $ can also be given a description as 
arithmetic quotient of a symmetric domain. The two connected components of 
the period domain, $D$ and $ \overline{D} $ are 
mapped one into each other by conjugation. This operation corresponds, on the right side of 
the isomorphism $ (\ref{s})$,  
%\begin{equation}
%\label{s2}
%\Omega \ \simeq \   O(2,18) / SO(2) \times O(18)   
%\end{equation}      
to changing the orientation of the positive two-plane. One obtains, therefore, an isomorphism:
\begin{equation}
\label{ma1} 
D \ \simeq \ O(2,18) / O(2) \times O(18).   
\end{equation}
Each isometry in $ \Gamma $, either 
preserves or exchanges the two connected components of $ \Omega $. One can precisely 
find the stabilizer $ \Gamma^+ = {\rm Stab}(D) $ as follows.  
%
%Let $ {\rm Diff}(S) $ be the group of diffeomorphisms of $S$ preserving 
%the elliptic structure and the section. 
%This group acts naturally on the family of compatible marked pairs 
%$ (S, \varphi) $, the quotient obtained after factoring out the action being 
%exactly the set of isomorphism classes of elliptic $K3$ surfaces with section, 
%$ \m{M} $.  Denote by $ \Gamma $ the group of isometries of the lattice $ \mathbb{L}_o $. 
%There is then an induced map: 
%\begin{equation}
%\label{map}
%{\rm Diff}(S) \rightarrow \Gamma .
%\end{equation}
%It was shown in [Bor][Matu][Don] that the image of $ (\ref{map}) $ is a subgroup
%$ \Gamma $ of index $2$ in $ O \left ( \mathbb{L}_o \right ) $. 
%
%We need to find the image of $ (\ref{map}) $ in $ \Gamma $. 
%
The orthogonal group $ O(2,18) $ of a real bilinear symmetric indefinite form 
of signature $(2,18) $ is a Lie group which has four connected components:
\begin{equation}
O(2,18) \ = \ O^{++}(2,18) \ \cup \ O^{+-}(2,18) \ \cup \ O^{-+}(2,18) \ \cup \ 
O^{--}(2,18). 
\end{equation} 
The upper signs refer to orientation behavior with respect to positive $2$-planes
and negative $18$-planes. The group of integral isometries can then 
be written as a disjoint union:
\begin{equation}
\Gamma \ = \ \Gamma ^{++} \ \cup \ \Gamma  ^{+-} \ \cup \  
\Gamma ^{-+} \ \cup \ \Gamma ^{--} 
\end{equation} 
by taking intersections of $ \Gamma $ with the components 
of the real orthogonal group. It follows then that the isometries preserving $D$ 
are exactly the ones preserving orientation on positive 2-planes:
%following arguments in \cite{borcea} \cite{matumoto} 
%cite{donald}, one concludes that the image of 
% $ (\ref{map}) $  is then given by 
\begin{equation}
\Gamma^+ \ = \ \Gamma ^{++} \ \cup \ \Gamma ^{+-} .  
\end{equation} 
The subgroup $ \Gamma^+ $ has index two in $ \Gamma $. We obtain then a model 
for the moduli space of elliptic $K3$ surfaces with section: 
\begin{equation}
\label{mmmmm} 
\m{M}_{K3} \ = \ \Gamma \backslash \Omega \ \simeq \ \Gamma^+ \backslash D   
\end{equation}
while the classical F-theory moduli spaces appears as:
\begin{equation}
\m{M}_F \ = \ \hat{\Gamma} \backslash \Omega \ \simeq \ \hat{\Gamma}^+ \backslash D   
\end{equation} 
where $ \hat{\Gamma}^+ =  \Gamma^+ \rtimes \Zee_2 $. However, it turns out that
$ \hat{\Gamma}^+ \simeq  \Gamma $ and, under this isomorphism, map
$ (\ref{ma1}) $ becomes an equivariant identification. One obtains therefore 
an arithmetic quotient picture for the classical F-theory moduli space as:
\begin{equation}
\m{M}_F \ \simeq \  \Gamma \backslash  O(2,18) / O(2) \times O(18).   
\end{equation}    
Along the lines of this description, the double-cover 
$ \m{M}_{K3} \rightarrow \m{M}_F $ can be seen as:
$$
\Gamma^+ \backslash  O(2,18) / O(2) \times O(18) \rightarrow 
\Gamma \backslash  O(2,18) / O(2) \times O(18).
$$ 

\subsection{Arithmetic of Compactification of $\m{M}_{K3}$}
\label{arith}
The moduli space $ \m{M}_{K3} $ is connected but not compact. There exists various arithmetic 
techniques aiming at compactifying $ \Gamma  \backslash  \Omega $. The simplest one is the 
Baily-Borel procedure \cite{bb} which we briefly review next. Later, we shall turn  
our attention to a particular case of Mumford's toroidal compactification \cite{ash} 
which plays a central role in the computation we undertake in this paper.  
\par The Baily-Borel procedure \cite{bb} introduces 
an auxiliary space $ \Omega^* $ with $ \Omega \subset \Omega^* \subset \Omega^{\vee} $.  
The topological boundary of $ \Omega \subset \Omega^{\vee} $ decomposes into a 
disjoint union of closed analytic subsets, called boundary components. There are 
two types of such components. Some are zero-dimensional and are represented by 
the points in in the real quadric 
$ \Omega^{\vee} \cap \mathbb{P} \left ( \mathbb{L}_o 
\otimes _{\Zee} \mathbb{R} \right ) $. The others are copies of $ \mathbb{P}^1 $ 
and are generated by the complexified images in  $ \mathbb{P} \left ( \mathbb{L}_o 
\otimes _{\Zee} \mathbb{C} \right ) $ of the 2-dimensional isotropic subspaces
of $ \mathbb{L}_o 
\otimes _{\Zee} \mathbb{R} $.   
Group theoretically, it can be seen that the stabilizer 
$$ {\rm Stab}(F) = \{ \ g \in O^{++}(2,18)  \ \vert \ gF=F \ \}  $$
of a boundary component $F$ is a maximal parabolic subgroup of $ O^{++}(2,18) $. A boundary
component $F$ is called then rational if its stabilizer $ {\rm Stab}(F) $ is defined over 
$ \mathbb{Q} $. The 
assignment $ P \to F_P $ with $ {\rm Stab}(F_P) = P $ determines a bijective 
correspondence between the set of proper maximal
parabolic subgroups of $  O^{++}(2,18) $ and the set of all rational boundary components. 
One defines then:
$$ \Omega^* \ = \ \Omega \ \cup \ \left ( \bigcup _P F_P  \right ) $$
where the right union is made over all proper maximal rational parabolics. The action of 
$ \Gamma $ extends naturally to $ \Omega ^* $. Moreover, one can
endow $ \Omega^* $ with the Satake topology, under which the $\Gamma$-action is
continuous. The Baily-Borel compactification appears then as:
\begin{equation}
\left ( \Gamma \backslash \Omega \right )^* \ \stackrel{{\rm def}}{=} \ 
\Gamma \backslash \Omega^*. 
\end{equation} 
The main features of this new quotient space are as follows (see \cite{bb} for details). 
The space $ \left ( \Gamma \backslash \Omega \right )^* $ is Hausdorff, compact, connected and 
can be given a structure of complex algebraic space. The quotient 
$ \Gamma \backslash \Omega $  is embedded in 
$ \left ( \Gamma \backslash \Omega \right )^* $ as a Zariski open subset. If 
$ I_i(  \mathbb{L}_o ) $, $ i \in \{1,2\} $ represents the set of primitive
isotropic sub-lattices of rank $i$ in  $ \mathbb{L}_o $ then the complement 
$$  \left ( \Gamma \backslash \Omega \right )^*  - \Gamma \backslash \Omega $$
consists of $ \vert \Gamma \backslash I_1(  \mathbb{L}_o ) \vert  $ points and   
$ \vert \Gamma \backslash I_2(  \mathbb{L}_o ) \vert  $ copies of 
$ {\rm PSL}(2, \Zee) \backslash \mathbb{H}. $
%Each curve is isomorphic 
%to a copy of:
%$$ {\rm PSL}(2,\Zee) \backslash \mathbb{H}. $$
\noindent Let us note that the complex conjugation involution on $ \Omega $ extends 
to $ \Omega ^* $. On boundary, it preserves the points and induces complex conjugation on the 
one-dimensional $ \mathbb{P}^1 $'s. The procedure provides therefore a compactification 
for the classical F-theory moduli space $ \m{M}_F $: 
$$  \hat{\Gamma}  \backslash \Omega \ \subset \ 
\left ( \hat{\Gamma} \backslash \Omega \right )^*  $$
with boundary strata given by points and copies of 
$ \left ( {\rm PSL}(2,\Zee) \rtimes \Zee_2 \right )   \backslash \mathbb{H} $ 
with $ \Zee_2 $ generated by $ \tau \to - \bar{\tau} $. 
\par It is known that Baily-Borel construction gives the minimal geometrically meaningful 
compactification of $ \Gamma \backslash \Omega $ in the sense that it is dominated by any 
other geometric
compactification. However, the disadvantage of the method is that the boundary has
large codimension (it consists of only points and curves) and contains only partial
geometrical information. One avoids these inconveniences by using a blow-up of
the Baily-Borel construction, the toroidal compactification of Mumford \cite{ash}.
This compactification, although not canonical in general, 
gives divisors as boundary components and carries significantly more 
geometrical information. The main arguments describing the construction, 
as presented in \cite{friedman} and 
\cite{frthesis}, are as follows. 
%For our purposes only parts of the new Mumford boundary components 
%will be relevant. We shall restrict therefore the review to the fundamentals of 
%those components as they are presented in \cite{friedman} and appendix of \cite{frthesis}. 
\par The Mumford boundary components associated to $ \Gamma \backslash \Omega $ involve again 
the maximal rational parabolic subgroups of $  O(2,18) $. These are stabilizers of non-trivial 
isotropic subspaces $ V_{\mathbb{Q}} \subset \mathbb{L}_o \otimes _{\Zee} \mathbb{Q} $. 
The lattice $ \mathbb{L}_o $ has signature $(2,18) $,
and hence, if $ V_{\mathbb{Q}} $ is isotropic then its dimension is either $2$ or $1$.
If $ {\rm dim} (V_{\mathbb{Q}}) = 1 $, then the associated Baily-Borel 
rational boundary component
$F$ is represented by just a point. Such a component is called of Type III. For 
$ {\rm dim} (V_{\mathbb{Q}}) = 2 $, the corresponding boundary component $F$ is 1-dimensional.
In this case $ F $ is said to be of Type II. Each rational Baily-Borel component $F$ will
determine a Mumford boundary component $ \m{B}(F) $. We shall be concerned here only with describing
the components of Type II for which the construction is canonical.    
\par Let $ V_{\mathbb{Q}} $ be a rank-two isotropic lattice and $F$ the associated Baily-Borel
component. We denote: 
$$ P(F) \ = \ {\rm Stab} \left ( V_{\mathbb{R}} \right ) \subset O(2,18) $$
$$ W(F) \ = \ {\rm the \ unipotent \ radical \ of} \ P(F) $$   
$$ U(F) \ = \ {\rm the \ center \ of} \ W(F) . $$ 
It turns out that $ U(F) $ is 1-dimensional (also definable over $\mathbb{Q}$) and the Lie algebra
of its real form can be described as:
\begin{equation}
u(F) \ = \ \{ N \in {\rm Hom} \left ( (\mathbb{L}_o)_{\mathbb{R}},  
(\mathbb{L}_o)_{\mathbb{R}} \right ) \ \vert \ {\rm Im}( N) \subset V_{\mathbb{R} } \ 
{\rm and} \  (Na,b) + (a,Nb) = 0, \ \forall \ a,b \in (\mathbb{L}_o)_{\mathbb{R}} \ \}.   
\end{equation} 
One obtains that any $ N \in u(F) $ satisfies $ N^2=0 $, 
$ \ {\rm Im}( N) = V_{\mathbb{R} } $ and  
$ {\rm Ker} (N) =  V_{\mathbb{R} }  ^{\perp} $. There is then an associated
weight filtration:
\begin{equation}
\label{wf}
0 \ \subset \ V_{\mathbb{R} } \ \subset \  V_{\mathbb{R} } ^{\perp} \ \subset 
\ (\mathbb{L}_o)_{\mathbb{R}}.
\end{equation}
We pick a primitive integral endomorphism $ N \in u(F) $ 
%such that 
%$ (.,N.) $ is positive definite on 
%$$ (\mathbb{L}_o)_{\mathbb{R}} / \ 
%\left ( U_{\mathbb{R} } \right ) ^{\perp} $$ 
and consider the groups:
$$ {\rm U}(N)_{\Cee} \ = 
\ \{ {\rm exp} \left ( \lambda N \right ) \ \vert \ \lambda \in \Cee \ \} $$
$$ {\rm U}(N)_{\Zee} \ = 
\ \{ {\rm exp} \left ( \lambda N \right ) \ \vert \ \lambda \in \Zee \ \} \ = 
\ {\rm U}(N)_{\Cee} \cap O^{++}(2,18;\Zee).  $$
The group $ {\rm U}(N)_{\Cee} $ acts upon the extended period domain 
$$ \Omega^{\vee} \ = \{ \ [z] \in \mathbb{P} 
\left ( \mathbb{L}_o \otimes _{\Zee} \Cee \right ) \ \vert \ (z,z)=0 \ \} .$$ 
providing an intermediate filtration $ \Omega \subset \Omega(F) \subset \Omega^{\vee} $ 
where $ \Omega(F) = {\rm U}(N)_{\Cee} \cdot \Omega $. 
\par One defines then the Mumford boundary component associated to $F$ as the space of
nilpotent orbits:
\begin{equation}
\m{B}(F) \ = \ \Omega(F)  /  {\rm U}(N)_{\Cee}. 
\end{equation}
In this setting, 
\begin{equation}
\label{princ1}
\Omega(F)  /  {\rm U}(N)_{\Zee} \rightarrow \m{B}(F) 
\end{equation}
is a holomorphic principal bundle with structure group $ {\rm U}(N)_{\Cee}  /  
{\rm U}(N)_{\Zee} \simeq \Cee^* $.
The inclusion
$$ \Omega / {\rm U}(N)_{\Zee} \ \hookrightarrow 
\ \Omega(F) /  {\rm U}(N)_{\Zee} $$
realizes $ \Omega  /   {\rm U}(N)_{\Zee} $ as an open subset in 
the total space of $ (\ref{princ1}) $. Let then:
\begin{equation} 
\overline {\Omega(F) /   {\rm U}(N)_{\Zee}    } \ = \ 
\left (   \Omega(F) / {\rm U}(N)_{\Zee}         \right )   
\times_{\Cee^*} \Cee. 
\end{equation}  
This amounts to gluing in the zero section in the $ \Cee^* $-fibration $ (\ref{princ1}) $.
One defines then:
$$ \overline {\Omega /   {\rm U}(N)_{\Zee}    } \ \stackrel{{\rm def}}{=} \  {\rm interior \ of \ the \ 
closure \ of} \ \Omega /   {\rm U}(N)_{\Zee} \  {\rm in} \ 
\overline {\Omega(F) /   {\rm U}(N)_{\Zee}    }. $$
Set-theoretically, one has:
$$ \overline {\Omega /   {\rm U}(N)_{\Zee}    } \ = \ \Omega / {\rm U}(N)_{\Zee} \ \sqcup \ 
\m{B}(F). $$
Finally:
\begin{equation}
\label{omegabar}
\overline{\Omega} \ \stackrel{{\rm def}}{=} \ \bigcup_F \ \overline {\Omega /   
{\rm U}(N)_{\Zee}    } \ = \ 
\Omega \ \sqcup \ \left ( \bigsqcup _F \m{B}(F)\right )  
\end{equation}
the union being performed over all rational Baily-Borel boundary components of Type II. 
This space inherits
a topology. The arithmetic action of $ \Gamma $ induces a closed discrete 
equivalence relation
on $ (\ref{omegabar}) $. The quotient space, denoted by $ \overline{\Gamma \backslash \Omega} $, 
enjoys the following properties (see \cite{ash}, \cite{friedman} for details): 
\begin{teo} \ 
\begin{itemize}
\item $ \overline{\Gamma \backslash \Omega} $ is a quasi-projective analytic variety.
\item $ \overline{\Gamma \backslash \Omega} $ contains  $ \Gamma \backslash \Omega $ as a 
Zariski open dense subset. 
\item The complement $ \overline{\Gamma \backslash \Omega} \ - \ \Gamma \backslash \Omega $ 
consists of two irreducible divisors. These divisors are quotients of smooth spaces 
by finite group actions.   
\end{itemize}
\end{teo}  
\noindent We shall denote the two divisors by $\m{D}_{E_8 \oplus E_8} $ and $ \m{D}_{\Gamma_{16}} $. 
The reason for this terminology is the following. The two Type II divisors in question correspond 
to the two distinct orbits in 
$ \Gamma \backslash I_2(\mathbb{L}_o ) $ where 
$ I_2(\mathbb{L}_o) $ is the set of primitive isotropic rank-two sub-lattices in $\mathbb{L}_o $. 
On can identify
the orbit to which a certain isotropic sub-lattice belongs using the following recipe. 
Let $ V \in I_2(\mathbb{L}_o) $. The quotient lattice $ V^{\perp} / 
V $ is even, unimodular, negative-definite and has rank $16$. It is known that, up 
to isomorphism, there exists only two lattices
of this type: $ -(E_8 \oplus E_8) $ and $ - \Gamma_{16} $. 
The two isomorphism classes perfectly differentiate the two orbits in $ \Gamma \backslash I_2(\mathbb{L}_o ) $.  
There are therefore only two distinct Baily-Borel boundary curves in 
$$ \left ( \Gamma \backslash \Omega \right )^* - \Gamma \backslash \Omega  $$ 
and, accordingly, there are two Type II components in Mumford's compactification.  
\par In fact, for each isotropic sub-lattice $V$ there is a natural projection:
\begin{equation}
\label{ps}
\m{B}(F) \rightarrow F 
\end{equation} 
defined by assigning to a nilpotent orbit $ \{ {\rm U}(N)_{\Cee} \cdot \omega \} $ the
complex line $ \{\omega\} ^{\perp} \cap V_{\Cee} \subset V_{\Cee} $. We shall see the 
geometrical significance of $ (\ref{ps}) $ in the next section. At this point, we just
note that these projections descend to maps from the Type II Mumford divisors to the
two Baily-Borel boundary curves under  
$$ \overline{\Gamma \backslash \Omega} \ 
\rightarrow \  
\left ( \Gamma \backslash \Omega \right )^*. $$   
\par As mentioned earlier, the main goal of this paper is to describe explicitly
the structure of $ \Gamma \backslash \Omega $ in a neighborhood of the two Type II divisors 
$\m{D}_{E_8 \oplus E_8} $ and $ \m{D}_{\Gamma_{16}} $. Our description will go along the 
following direction. Let $ F $ be a Type II Baily-Borel
component and denote by $ \Gamma_F = P(F) \cap \Gamma $
the stabilizer of the associated isotropic sub-lattice $V$. As subgroup of $ \Gamma $, 
the group $ \Gamma_F $ induces an equivalence relation on 
$ \overline{\Omega} $ dominating the $ \Gamma $-one. One obtains therefore the 
following sequence of analytic projections:
\begin{equation}
\overline{\Omega} \ \rightarrow  \ \Gamma_F \backslash \overline{\Omega}  \ \rightarrow \  
\overline{\Gamma \backslash \Omega}.
\end{equation} 
Then, as explained in Chapter 5 of \cite{ash}:
\begin{lem}
There exists an open subset 
$$ \m{U}_F \ \subset \ \overline {\Omega /   {\rm U}(N)_{\Zee}    } \ \subset \ \overline{\Omega}, $$
tubular neighborhood of the Mumford boundary component $ \m{B}(F) \subset \overline{\Omega} $ 
such that on $ \Gamma_F \cdot \m{U}_F $, the $ \Gamma $-equivalence reduces 
to $ \Gamma_F $-equivalence.    
\end{lem} 
\noindent In the light of this lemma, the analytic projection:
\begin{equation}
\Gamma_F \backslash \left( \Gamma_F \cdot \m{U}_F \right ) \ \rightarrow \ 
\Gamma \backslash \left( \Gamma_F \cdot \m{U}_F \right ) 
\end{equation}
is an isomorphism. One has therefore an analytic identification between an open 
neighborhood of the Mumford divisor associated to $F$ in $ \overline{\Gamma \backslash \Omega} $ and 
$$ \m{V}_F \stackrel{{\rm def}}{=} \  
\Gamma_F \backslash \left( \Gamma_F \cdot \m{U}_F \right ) \ \subset \ 
\Gamma_F \backslash \ \overline{\Omega/  {\rm U}(N)_{\Zee} } \ \ \subset \ 
\Gamma_F \backslash \ \overline{\Omega(F) /  {\rm U}(N)_{\Zee} }.   
$$  
But, as observed earlier, 
\begin{equation}
\label{parabcover1}
\w{\Theta} \colon 
\overline{\Omega(F) /  {\rm U}(N)_{\Zee} } \ \rightarrow \ \m{B}(F) 
\end{equation} 
is a holomorphic line bundle. After factoring out the action of $ \Gamma_F $, 
one obtains a holomorphic $ \Cee $-fibration:  
\begin{equation}
\label{parabcover} 
\Theta \colon \Gamma_F \backslash  \overline{\Omega(F) /  {\rm U}(N)_{\Zee} } \ \rightarrow \ 
\Gamma_F \backslash \m{B}(F). 
\end{equation}
It is easy to see that $ \m{V}_F $ is a tubular neighborhood of the zero-section
in $ (\ref{parabcover}) $. 
\par Based on the above arguments, one concludes that an open subset of the
period domain $ \Gamma \backslash \Omega $ which is a neighborhood of one
of the two possible Type II divisors can be identified with an open 
neighborhood of the zero-section in the parabolic fibration $ (\ref{parabcover}) $. 
Therefore, in order to describe the structure of $ \m{M}_{K3} $ in the vicinity
of one of the two Type II divisors $\m{D}_{E_8 \oplus E_8} $ and $ \m{D}_{\Gamma_{16}} $, 
it is essential to explicitly describe
the holomorphic type of $ (\ref{parabcover}) $. We accomplish this task in 
section $ \ref{main} $. 
\par We finish this section with a note on the behavior of complex conjugation
within the framework of the above construction. The complex conjugation on 
$ \Omega $ extends naturally to an involution of $ \overline{\Omega} $ giving 
producing complex conjugations on each Type II Mumford boundary component $ \m{B}(F) $. 
One can perform therefore a similar partial 
compactification:
$$ \m{M}_F =  \hat{\Gamma} \backslash \Omega \ \subset \ \overline{\hat{\Gamma} 
\backslash \Omega} $$
with $ \overline{\hat{\Gamma} \backslash \Omega} \ - \ \hat{\Gamma} \backslash \Omega $ 
consisting of two boundary divisors (obtained as quotients of the two Type II divisors
of $\Gamma \backslash \Omega$ by complex conjugation). Open neighborhoods of $\m{M}_F$ 
near the boundary divisors are still described by open neighborhoods of
the zero-section in the total space of the parabolic 
cover $ (\ref{parabcover}) $.    
\subsection{Boundary Components and Hodge Structures}
\label{p}
One can give a Hodge theoretic interpretation for the boundary component $\m{B}(F)$. A period  
$ \omega \in \Omega $ determines automatically a polarized Hodge structure of weight two 
on $ \mathbb{L} \otimes_{\Zee} \Cee $, corresponding geometrically to a marked elliptic $ K3 $ surface
with section. Taking orthogonals with respect to the fixed hyperbolic 
sub-lattice $ H \subset \mathbb{L} $ (which by construction 
consist of $(1,1)$-cycles and is therefore orthogonal to the period line),
one obtains a polarized Hodge structure of weight two on $ \mathbb{L}_o \otimes_{\Zee} \Cee $, 
\begin{equation}
\label{filtr1}
0 \ \subset \ \{\omega \} \ \subset \ \{\omega \} ^{\perp} \ \subset \ \mathbb{L}_o
\otimes_{\Zee} \Cee. 
\end{equation}
Let then $ V \subset \mathbb{L}_o $ be the primitive isotropic rank-two sub-lattice corresponding to the 
Type II Baily-Borel boundary component $F$. There is an induced weight filtration:
\begin{equation}
\label{weightf1}
0   \ \subset \ {\rm V}_{\Cee} \ \subset \  \left ( {\rm V}_{\Cee}  \right ) ^{\perp} \   \subset \ 
 \mathbb{L}_o \otimes \Cee. 
\end{equation}
Together, filtrations $ (\ref{filtr1}) $ and $(\ref{weightf1}) $ yield a mixed Hodge structure on 
$ \mathbb{L}_o \otimes_{\Zee} \Cee $. Taking this point of view, one can regard the 
domain $ \Omega(F) = U(N)_{\Cee} \cdot \Omega $
as the space of mixed Hodge structures on the weight filtration 
$ (\ref{weightf1}) $. These structures are acted 
upon by the group $ U(N)_{\Cee} $. The Type II Mumford boundary component  
$$ \m{B}(F) = \Omega(F) / U(N)_{\Cee} $$
appears then as the space of nilpotent orbits of such mixed Hodge structures.  
\par There are three $U(N)_{\Cee}$-invariant graded pure Hodge structures associated 
to each nilpotent orbit in $ \m{B}(F) $:
\begin{equation}
\label{ho1}
0 \ \subset \ \{\omega\}^{\perp} \cap V_{\Cee} \ \subset \ V_{\Cee}  
\end{equation}
\begin{equation}
\label{ho2}
0 \ \subset \  \left ( \{ \omega \} \cap V_{\Cee}^{\perp} + V_{\Cee} \right )  /  V_{\Cee} \ 
\subset \ \left ( \{ \omega \} ^{\perp} \cap V_{\Cee}^{\perp} + V_{\Cee} \right )  / 
V_{\Cee} \ \subset \ V_{\Cee}^{\perp}  /  V_{\Cee} \    
\end{equation}
\begin{equation}
\label{ho3}
0 \ \subset \  \left ( \{ \omega \} + V_{\Cee}^{\perp}  \right )  /  V_{\Cee}^{\perp} \ 
\subset \ \left ( \{ \omega \} ^{\perp} + V_{\Cee}^{\perp} \right )  / 
V_{\Cee}^{\perp} \ \subset \ (\mathbb{L}_o)_{\Cee}  /  V_{\Cee}^{\perp} \    
\end{equation}
\noindent The first one, which we denote by $\m{H} $, is a pure Hodge structure of weight one 
induced on  $V_{\Cee}$ and is polarized with respect to a certain
non-degenerate skew-symmetric bilinear form $ ( \cdot, \cdot)_1$ on $V$. 
Let $(\cdot, \cdot)_3 $ 
be the bilinear form on $ \mathbb{L}_o / V^{\perp} $ given by $ (x,y)_1= (x,Ny)$
and $ ( \cdot, \cdot)_1$ be the form on $V$ under which the
isomorphism:
\begin{equation}
\label{form1}
N \colon \mathbb{L}_o / V^{\perp} \rightarrow V  
\end{equation}
becomes an isometry. One has then $ (x,y)_1 = ( \w{x}, y) $ for $x,y \in V$, where 
$ \w{x} $ is a lift of $x$ to $\mathbb{L}_o$. The bilinear form $ ( \cdot, \cdot)_1$
is non-degenerate and skew-symmetric. The space $\m{B}(F) $ of nilpotent $ U(N)_{\Cee} $-orbits has two
connected components and one can check that the Hodge structure
$ (\ref{ho1}) $ is polarized with respect to $ ( \cdot, \cdot)_1$ or 
$-( \cdot, \cdot)_1$ depending on the component the nilpotent orbit is part of. We agree
to denote by $\m{B}^+(F) $ the component for which $ (\ref{ho1}) $ is polarized with respect to 
$( \cdot, \cdot)_1$. Then 
$$ \m{B}(F) = \m{B}^+(F) \sqcup \overline{\m{B}^+(F)}. $$          
\par The second graded Hodge structure, described in $ (\ref{ho2}) $, has pure weight two and can be 
seen to be of type $ (1,1) $. Finally, the third Hodge structure $ (\ref{ho3}) $ has weight three, 
but one can check that, under the isomorphism $ (\ref{form1}) $, filtration $ (\ref{ho3}) $ is just 
the $(1,1)$-shift of Hodge structure $ (\ref{ho1}) $.    
\par A mixed Hodge structure contains considerably more than the sum of its graded pieces. The 
first two graded parts are glued together by the extension of mixed Hodge structures:
\begin{equation}
\label{ddss}
\{0\} \rightarrow V \rightarrow V^{\perp} \rightarrow V^{\perp}/V \rightarrow  \{0\}. 
\end{equation}
In fact, one can check that the Hodge structure $ (\ref{ho1}) $ together with the 
extension $ (\ref{ddss}) $ completely determines the nilpotent orbit of $\omega $.  
This gives a natural isomorphism between $ \m{B}(F) $ and the space of equivalence 
classes of extensions of type $(\ref{ddss}) $. Such extensions of mixed Hodge structures 
have been studied by Carlson in \cite{carl}. They are classified, 
up to isomorphism, by an abelian group homomorphism:
\begin{equation}
\label{extmap1}
\psi \colon 
%\left ( {\rm V}^{\perp} /   {\rm V}   \right )_{\Zee} 
\Lambda \rightarrow 
%{\rm Jac}(E) \simeq E.  
J^1(\m{H})
\end{equation}
where $\Lambda $ is the lattice $ \left ( {\rm V}^{\perp} /   {\rm V}   \right )_{\Zee}$ 
and $J^1(\m{H} )\ = \ V_{\Cee}  /  
\left(  \{\omega\}^{\perp} \cap V_{\Cee} + V_{\Zee} \right ) $ 
is the generalized Jacobian associated to the pure Hodge structure $ \m{H} $ described in 
$ (\ref{ho1}) $. As mentioned before, $\Lambda$ has to be unimodular, even, negative-definite
and of rank $16$.  
%\par The extension map $ (\ref{extmap1}) $ can be described explicitly in our case 
%as follows. Let $ \gamma \in {\rm V}^{\perp}_{\Zee} $. Since 
%$ \omega \notin {\rm V}_{\Cee}^{\perp} $, one can always find an element 
%$ z \in {\rm V}_{\Cee} $ 
%such that $ (\gamma + z , \omega ) = 0 $. Given another $ z' \in  {\rm V}_{\Cee} $
%satisfying the same propriety, we have $ (z'-z, \omega ) =0 $ and therefore,
%$ z'-z \in \{ \omega \}^{\perp} \cap {\rm V}_{\Cee} $. Moreover, if $ \gamma $ is
%varied by an integral element, $ \gamma' = \gamma + u $ , $u \in {\rm V}_{\Zee} $ 
%and $ z'' \in {\rm V}_{\Cee} $ is chosen with $ (\gamma' + z'' , \omega) =0 $ then 
%$ u + z''-z \in \{ \omega \}^{\perp} \cap {\rm V}_{\Cee} $ and therefore
%$ z'' - z \in \{ \omega \}^{\perp} \cap V_{\Cee} + U_{\Zee} $. This recipe 
%provides a well-defined homomorphism:
%\begin{equation}
%\label{extmap44}
%\left ( {\rm V}^{\perp}/  {\rm V} \right ) _{\Zee} \ \rightarrow \   
%{\rm V}_{\Cee} / \left ( 
%\{ \omega \}^{\perp} \cap {\rm V}_{\Cee} + {\rm V}_{\Zee} 
%\right ) \ , \ \ \ 
%\end{equation}
%$$ 
%\gamma \ \left( {\rm mod} \ {\rm V}_{\Zee} \right ) \ \ \mapsto \ 
%z \ \left ( {\rm mod}              \ \{ \omega \}^{\perp} \cap {\rm V}_{\Cee} + {\rm V} _{\Zee}
%\right ) 
%$$   
%which provides the model for $ (\ref{extmap1}) $. 
\par One obtains then that, for a given Type II Baily-Borel component $F$, 
the Mumford boundary points lying in $ \m{B}(F) $ can be identified with pairs 
$ (\m{H}, \psi ) $ 
consisting of polarized Hodge structures $\m{H}$ of weight one on $ {\rm V}_{\Cee} $ together
with homomorphisms $ \psi \colon \Lambda \rightarrow J^1(E) $ where 
$\Lambda = \left ( {\rm V}^{\perp} /   {\rm V}   \right )_{\Zee}$. 
The projection to the $\m{H}$-component $ (\m{H}, \psi ) \rightarrow \m{H}$ recovers exactly the
projection:
$$ \m{B}(F) \ \rightarrow \ F $$ 
mentioned in the arithmetic discussion of previous section.    
\section{Stable $K3$ Surfaces}
\label{stablek3}
To this point we have described the partial compactification:
\begin{equation}
\label{procarith}
% \m{M}_{K3} \simeq  
\Gamma \backslash \Omega \ \subset \overline{\Gamma \backslash \Omega} 
\end{equation} 
from a purely arithmetic point of view. In this section, we claim that the above
compactification also has a geometrical interpretation. Namely, under the period 
map identification $\m{M}_{K3} = \Gamma \backslash \Omega $,  $ (\ref{procarith}) $ 
amounts to enlarging the moduli space $ \m{M}_{K3} $ by allowing certain explicit degenerations 
of elliptic $K3$ surfaces with section.
\par Let $ \Lambda_1 = \left (E_8 \oplus E_8 \right ) $ and $ \Lambda_2 = \Gamma_{16} $
be the two possible equivalence classes of unimodular, even, positive-definite lattices
of rank $16$. We claim that there is an identification:       
$$
\xymatrix{
\left \{ \txt{points on the Mumford\\ 
boundary divizor $ \m{D}_{\Lambda_i} $}
\right \} 
\ar @{<->}  [r] &     
\left \{ 
\txt{elliptic Type II stable \\ 
$K3$ surfaces with section \\ 
in  $\Lambda_i$-category} 
\right \}
}
$$
and furthermore, the above correspondence can be regarded as a natural extension of the period map
to the boundary. 
\subsection{Definition and Examples}
\noindent Let us start by reviewing the notion of a Type II stable $K3$ surface 
(following \cite{friedman1} \cite{friedman}) and the reason why these objects are 
natural geometrical candidates to be associated with the arithmetic 
Type II Mumford boundary points. 
\begin{dfn} (\cite{friedman}) A ${\bf Type \ II \ stable \ K3 \ surface} $ is a surface with normal 
crossings  
\begin{equation}
Z_o \ = X_1 \cup X_2 
\end{equation}   
satisfying the properties:
\begin{itemize}
\item $X_1 $ and $X_2 $ are smooth rational surfaces.
\item $ X_1 $ and $ X_2 $ intersect with normal crossings and $ D = X_1 \cap X_2 $
is a smooth elliptic curve. 
\item $ D \in \vert - K_{X_i} \vert $ for $ i=1,2$.
\item $ N_{D/X_1} \otimes N_{D/X_2} = \m{O}_D $ (d-semi-stability).    
\end{itemize}
\label{dddd33}
\end{dfn}
\noindent Let us note that the above conditions imply 
that $ \omega_{Z_o} \simeq \m{O}_{Z_o} $, where $\omega_{Z_o} $ is the dualizing sheaf.  
Specializing the above definition, we say that a Type II stable $K3$ surface is endowed with an 
$ {\bf elliptic \ structure \ with \ section}$ if $Z_o$ is in one of the following categories:
\begin{itemize}
\item [(a)] Both smooth rational surfaces $X_i$ are endowed with elliptic fibrations 
$ X_i \rightarrow \mathbb{P}^1 $ with sections $S_i \subset X_i $. The double curve $D$ is 
a smooth elliptic fiber on both sides. The two sections $S_1$ and $S_2 $ meet $D$
at the same point. 
\item [(b)] Both smooth 
rational surfaces $X_i$ carry rulings defining maps $ X_i \rightarrow \mathbb{P}^1 $. The two
restrictions on the double curve $D$ agree, providing the same branched double-cover
$D \rightarrow \mathbb{P}^1$. In addition $X_1 $ is endowed with a fixed section of 
the ruling, denoted $S_o$, disjoint from $D$.         
\end{itemize} 
\noindent In short, a stable surface $Z_o$ is, in the case (a), the total space of 
an elliptic fibration $X_1 \cup X_2 \rightarrow \mathbb{P}^1 \cup \mathbb{P}^1$ with 
a fixed section given by $ S_o = S_1 \cup S_2 $. In the case (b), $Z_o$ is the total 
space of a fibration $ X_1 \cup X_2 \rightarrow \mathbb{P}^1 $ whose generic fiber 
is a union of two smooth rational curves meeting at two points. The fixed rational 
curve $S_0 \subset X_1 - D $ is a section for the fibration. For reasons to be 
clarified shortly, we shall sometime refer to (a) and (b) as 
$ E_8 \oplus E_8 $ and $ \Gamma_{16} $ categories, respectively. 
\par Two elliptic Type II stable $K3$ surfaces with section 
$Z_o$ and $Z_o'$ are said to be isomorphic if there exists an isomorphism of analytic 
varieties $ f \colon Z_o \rightarrow Z_o'$ entering a commutative diagram (depending on the category):
\begin{equation}
\label{iiizzzooo}
\xymatrix{
X_1 \cup X_2 \ar [dr] \ar @/^/ [rr] _f &                                & X_1' \cup X_2' \ar [dl] \\
                              & \mathbb{P}^1 \cup \mathbb{P}^1 \ar @/^/ ^{S_o} [ul] \ar @/_/ _{S_o'} [ur] & 
\\             
} \ \ \ \ \ \ \ 
\xymatrix{
X_1 \cup X_2 \ar [dr] \ar @/^/ [rr] _f &                                & X_1' \cup X_2' \ar [dl] \\
                              & \mathbb{P}^1 \ar @/^/ ^{S_o} [ul] \ar @/_/ _{S_o'} [ur] &. \\ 
}
\end{equation}
\par The reasons why one considers the configurations in (a)-(b) as elliptic structures with section on
a stable $K3$ surface will be explained in section $ \ref{sdegen} $. Let us next 
describe explicit examples of 
such configurations. Our construction pattern is as follows. Let $E $ be a smooth elliptic curve. 
Consider $p_0, q_0 \in E$ and let 
$$ E \stackrel{\varphi_1}{\hookrightarrow} \mathbb{P}^2\ \ \ \ \ \ 
E \stackrel{\varphi_2}{\hookrightarrow} \mathbb{P}^2,  $$
be the projective embeddings determined by the linear systems $ \vert 3 p_0 \vert $ and 
$ \vert 3 q_0 \vert $. Pick $18$ more points  $ p_1, p_2, \cdots p_{18} $ 
(not necessarily distinct) on $E$ and partition them into two ordered subsets 
$$ \{ p_1, p_2 , \cdots p_t \} \cup \{ p_{t+1} , p_{t+2} , \cdots p_{18} \}. $$ 
Blow up the first copy of $\mathbb{P}^2 $ 
at $ p_1, p_2 , \cdots p_t $ (in the given order) and perform the same blow-up procedure on 
the second copy of 
$\mathbb{P}^2 $ using the points  $p_{t+1} , p_{t+2} , \cdots p_{18}$. Let $X_1$ and $X_2$
be the resulting surfaces. A surface $Z_o$ with normal crossings is 
obtained by gluing $X_1$ and $X_2$ 
together along the proper transforms of $ \varphi_1(E) $ and $ \varphi_2(E) $ using, as 
gluing map, the isomorphism $ (\varphi_2)^{-1} \circ \varphi_1 $.    
\begin{dfn}
\label{defex}
A collection $ \{ 3p_0;p_1,p_2, \cdots p_t; 3q_0; p_{t+1} \cdots p_{18} \} $, with $3p_0 $ and
$3q_0$ considered as divisor classes in $ {\rm Pic}(E)$, is called 
a ${\bf special \ family} $ 
if one of the following sets of conditions holds:
\begin{itemize}
\item [(a)] $t=9$, $p_9 = p_{18} $ and
$$ \m{O}_E(p_1+p_2+ \cdots +p_8 +p_9) = \m{O}_E(9p_0), \ \ 
\m{O}_E(p_{10} + p_{11} + \cdots + p_{18}) = \m{O}_E( 9 q_0).   $$  
 \item [(b)] $ 2\leq t \leq 17$, $p_1 = p_2 $ and
$$ \m{O}_E(p_1+p_2+ \cdots +p_{18}) = \m{O}_E( 9p_0 +9 q_0), \ \ 
\m{O}_E(3p_0 -p_1) = \m{O}_E(3q_0 - p_{t+1}).   $$  
\end{itemize} 
\end{dfn}
\noindent Let 
$ \{ 3p_0;p_1,p_2, \cdots p_t; 3q_0; p_{t+1} \cdots p_{18} \} $ 
be a special family on $E$. Denote by 
$$Z_o\left (E;3p_0;p_1,p_2, \cdots p_t; 3q_0; p_{t+1} \cdots p_{18} \right ) $$
the surface with normal crossings constructed by the pattern described earlier. 
\begin{teo}
\label{teoexample}
The surface:
$$Z_o\left (E;3p_0;p_1,p_2, \cdots p_t; 3q_0; p_{t+1} \cdots p_{18} \right ) $$ 
is an elliptic Type  II stable $K3$ surface with section. Moreover the surface falls 
in category (a) when the special family satisfies condition (a), and in 
category (b) when the special family satisfies condition (b).
\end{teo}
\begin{proof}
Assume that  $\{ 3p_0;p_1,p_2, \cdots p_9; 3q_0; p_{10} \cdots p_{18} \} $ is 
a special family of type (a). Then, 
the double curve $D$ of $Z_o$ is smooth elliptic and satisfies  
$ D \in \vert -K_{X_i} \vert $, $ D^2 = 0 $. A computation involving Riemann-Roch theorem 
leads to  $ h^0(X_i, D ) = 2 $. The linear system $ \vert D \vert $ is a base-point 
free pencil on each $X_i$ 
and induces elliptic fibrations $ X_i \rightarrow \mathbb{P}^1 $. The exceptional 
curves $E_{9} $ and $E_{18} $ corresponding to $p_{9}$ and $p_{18}$ 
are sections in the two fibrations and they meet the double
curve $D$ at the same point. The d-stability condition on 
$Z_o $ is satisfied as both normal bundles 
$ N_{D/X_i}$ are holomorphically trivial. We have therefore an explicit model 
$$  Z_o \left ( E; 3p_0;p_1,p_2, \cdots p_9; 3q_0; p_{10} \cdots p_{18} \right ) = 
X_1 \cup X_2 \rightarrow \mathbb{P}^1 \cup \mathbb{P}^1 $$
for an elliptic Type II stable $K3$ surface with section in the (a)-category.
\par We treat now the case when 
$\{ 3p_0;p_1,p_2, \cdots p_t; 3q_0; p_{t+1} \cdots p_{18} \} $ is a special 
family of type (b). Let $H_1$, $H_2$ be the hyper-plane divisors of the two copies of $\mathbb{P}^2$ 
and denote by $E_i$ the exceptional curve corresponding to $p_i$. The linear systems 
$ \vert H_1 - E_1 \vert $ and $ \vert H_2 - E_{t+1} \vert $ are base-point free pencils 
inducing rulings $ X_i \rightarrow \mathbb{P}^1 $. The restrictions of the two rulings 
agree on the double curve $D$,
recovering the branched double cover $ E \rightarrow \mathbb{P}^1 $
associated to the pencil $ \vert 3p_0 - p_1 \vert = \vert 3q_0 - p_{t+1} \vert$. Moreover, 
if one denotes by $S_o$ the proper 
transform of $ E_{1} $ in $X_1$, then $S_o$ is a smooth rational curve, 
with self-intersection $-2$, disjoint from $D$, 
and realizing a section of the ruling
$X_1 \rightarrow \mathbb{P}^1$. The d-semi-stability condition on 
$Z_o\left (E;3p_0;p_1,p_2, \cdots p_t; 3q_0; p_{t+1} \cdots p_{18} \right ) $ 
is satisfied since the line bundle 
$ N_{D/X_1} \otimes N_{D/X_2} $ is represented on $E$ by the principal
divisor $ 9p_0 + 9q_0 - p_1- p_2 - \cdots - p_{18} $. We obtain therefore 
an elliptic Type II stable $K3$ surface with section
$$ Z_o\left (E;3p_0;p_1,p_2, \cdots p_t; 3q_0; p_{t+1} \cdots p_{18} \right )
=X_1 \cup X_2 \rightarrow \mathbb{P}^1 $$
in the (b)-category.    
\end{proof}
\noindent The surfaces of Theorem $ \ref{teoexample} $ represent quite a large set of examples of 
elliptic Type II stable $K3$ surfaces with section. In fact, one can see that, up to certain 
explicit transformations, these surfaces actually exhaust all possibilities. 
\begin{dfn}
\label{dedede}
Let $Z_o$ be an elliptic Type II stable $K3$ surface with section. A ${\bf blowdown}$ 
$$ \rho \colon Z_0 \rightarrow \mathbb{P}^2 \cup \mathbb{P}^2 $$ consists of two sequences
of applications:
\begin{equation}
\label{bd1}
\rho_1 \colon X_1 = X_1^{(n)} \rightarrow X_1^{(n-1)} \rightarrow \cdots \rightarrow 
X_1^{(1)} \rightarrow X_1^{(0)} 
\end{equation}
\begin{equation}
\label{bd2}
\rho_2 \colon X_2 = X_2^{(m)} \rightarrow X_2^{(m-1)} \rightarrow \cdots \rightarrow 
X_2^{(1)} \rightarrow X_2^{(0)} 
\end{equation}  
such that:
\begin{enumerate}
\item The surfaces $X_1^{(0)} $ and $X_2^{(0)} $ are copies of $\mathbb{P}^2 $.   
\item Each map $ X_i^{(l)} \rightarrow X_i^{(l-1)} $ is a contraction of an exceptional curve 
in  $ X_i^{(l)} $.
\item If $Z_o$ is of type (a) then $S_1$, $S_2$ are the exceptional curves associated to 
$ X_1^{(n)} \rightarrow X_1^{(n-1)} $ and $ X_2^{(m)} \rightarrow X_2^{(m-1)} $.
\item If $Z_o$ is of type (b) then the exceptional curve associated to 
$ X_i^{(l)} \rightarrow X_i^{(l-1)} $, $ l \geq 2 $, is a component of a reducible 
fiber of the ruling. Moreover, for $ i=1$, $ l \geq 3 $ this exceptional curve is disjoint 
from $S_o$.    
\end{enumerate} 
\end{dfn} 
\noindent Due to their specific construction pattern, the special surfaces 
$ Z_o\left (E;3p_0;p_1,p_2, \cdots p_t; 3q_0; p_{t+1} \cdots p_{18} \right ) $ carry 
a canonical blow-down. Furthermore, if a stable surface $Z_o$ admits a blow-down, 
then $Z_o$ is isomorphic to a special surface of Theorem $ \ref{teoexample} $. Indeed, 
let us assume a choice of blow-down 
$ \rho \colon Z_o \rightarrow \mathbb{P}^2 \cup \mathbb{P}^2 $. Choose 
$ p_0, q_0 \in D$ such that $ 3p_0 $ and $3q_0$ 
are hyper-plane section divisors for the embeddings:
$$ D \hookrightarrow X_i \stackrel{\rho_i}{\rightarrow} X_i^{(0)} , \ \ i = 1,2. $$
Let 
$ x_1, x_2, \cdots , x_n, y_1, y_2, \cdots y_m $ be the points of intersection between 
the exceptional curves of $ X_i^{(l)} \rightarrow X_i^{(l-1)} $ and the double curve $D$. 
A cohomology calculation shows that $ m+n = 18$. Then, one can see that 
$ \{3p_0; x_1, x_2, \cdots, x_n; 3q_0; y_1, y_2, \cdots y_m\} $ is a special family
and there is a canonical isomorphism:
$$ Z_o \ \simeq \ Z_o ( D;3p_0; x_1, x_2, \cdots, x_n; 3q_0; y_1, y_2, \cdots y_m) $$
restricting to identity over $D$.       
\begin{prop}
\label{paa}
For any elliptic Type II stable $K3$ surface with section $Z_o$ in category (a), there 
exists a blow-down $ \rho \colon Z_o \rightarrow \mathbb{P}^2 \cup \mathbb{P}^2 $. 
\end{prop}
\begin{proof}
As rational surfaces, both $X_1 $ and $X_2$ have to dominate one of the geometrically 
ruled rational 
surfaces $ \mathbb{F}_n $, $ n \geq 0 $. Since $ D $ meets all exceptional curves, the double 
curve has to be the proper transform of an effective anti-canonical divisor in $ \mathbb{F}_n $. 
Such divisors exist only if $ n \leq 2 $. But $ \mathbb{F}_1 $ dominates $ \mathbb{P}^2 $ 
and $ \mathbb{F}_0 $, $ \mathbb{F}_2 $ also dominate $ \mathbb{P}^2 $ after blowing up a 
point on an anti-canonical curve. Therefore, if $X_i$ is neither 
$ \mathbb{F}_o $ nor $ \mathbb{F}_2 $ (which is the case here since $X_1$, $X_2$ are elliptic), 
one can always find blow-up sequences as in
$ (\ref{bd1}) $ and $ (\ref{bd2}) $.  Since $D^2=0$ on each $X_i$, 
it has to be that $n=m=9$. Moreover, one can always choose the section components $S_i$ as
the first exceptional curves to be contracted on each side. We have therefore a blow-down
$ \rho \colon Z_0 \rightarrow \mathbb{P}^2 \cup \mathbb{P}^2 $ as in definition $ (\ref{dedede}) $.
%The linear system $ \vert D \vert $ on $X_1$ (respective $X_2$ ) corresponds to the family of 
%cubics in $ X_1^{(0)} $ (respective $ X_2^{(0)} $) passing through 
%$\rho_1(x_1), \rho_1(x_2), \cdots \rho_1(x_9)$ (respective $\rho_2(y_1), \rho_2(y_2), \cdots \rho_2(y_9)$). 
%One cubic in the family is $D_i = \rho_i(D) \subset X_i^{(0)} $. However, $ \vert D \vert $ arpencils on 
%each side and therefore the nine points $\rho_1(x_1), \rho_1(x_2), \cdots \rho_1(x_9)$ 
%(respective $\rho_2(y_1), \rho_2(y_2), \cdots \rho_2(y_9) $ have to 
%lie on an extra cubic, distinct of $ D_i $. 
%One has therefore that $ x_1 + x_2 + \cdots x_9 - 9 p_0 $ and $ y_1 + y_2 + \cdots y_9 - 9 q_0 $ 
%vanish in $ {\rm Pic}^o(D) $.  
\end{proof}
\noindent Not all stable surfaces of category (b) admit blow-downs in the sense of Definition
$ \ref{dedede} $. $X_2 $ may be $ \mathbb{F}_0 $ or $ \mathbb{F}_2 $ and $X_1 $ may be $ \mathbb{F}_2 $.
None of these surfaces dominate $ \mathbb{P}^2 $. However, it can be shown that any $Z_o$ of category
(b) can be transformed, using certain explicit modifications, to a surface that admits blow-downs.    
\par An ${\bf elementary \ modification} $ of an elliptic Type II stable $K3$ surface with section 
$Z_o$ of category (b) consists of the blow-down of an exceptional curve $C$ lying inside a fiber of the 
ruling $X_i \rightarrow \mathbb{P}^1 $ and disjoint from $S_o$, 
followed by the blow-up of the resulting point on the opposite rational surface. The resulting 
$Z'_o$ is still an elliptic Type II $K3$ surface with section in category (b). 
\begin{prop}
\label{teooob}
Any elliptic Type II stable $K3$ surface $Z_o$ with section of category (b) can be transformed, 
using elementary modifications, to a surface which admits a blow-down.
\end{prop}
\begin{proof}  We claim that, using elementary modifications, one can transform
$Z_o$ to a new stable surface $Z'_o$ such that $ X'_2 = \mathbb{F}_1$. Indeed, using an argument
mentioned during the proof of Proposition $ \ref{paa} $, $X_2$ is either $ \mathbb{F}_0 $, 
or $ \mathbb{F}_2 $ or dominates $\mathbb{F}_1 $. If there is a blow-down $ X_2 \rightarrow \mathbb{F}_1 $ 
then perform elementary transformations consisting of flipping successively to $X_1 $ 
the exceptional curves involved in the blow-down. The new $X'_2$ is clearly $\mathbb{F}_1$. 
If $X_2$ is rather a copy of $\mathbb{F}_0$ or $ \mathbb{F}_2 $ then, choose an exceptional curve
$C$ sitting inside a fiber of the ruling on $X_1$ ($X_1$ and $X_2$ cannot be simultaneously geometrically
ruled). Let $p$ be the point where $C$ meets the double curve. Perform the elementary transform that takes 
$C$ to $X_2$ and flip back to $X_1$ the proper transform of the initial rational fiber through $p$ in
$X_2$. The resulting $X'_1 $ is then a copy of $ \mathbb{F}_1 $. Contracting the unique section of
negative self-intersection one obtains:
\begin{equation}
\label{bd22}
X_2 \rightarrow X_2^{(0)} 
\end{equation} 
with $ X_2^{(0)} $ isomorphic to a projective space $ \mathbb{P}^2 $.  
\par Assuming $ X_2 = \mathbb{F}_1 $, one has that $ X_1$ is ruled but not
geometrically ruled. Let us then describe the blow-down process 
\begin{equation}
\label{bd11}
X_1 = X_1^{(n)} \rightarrow X_1^{(n-1)} \rightarrow \cdots \rightarrow 
X_1^{(1)} \rightarrow X_1^{(0)} .
\end{equation}
We contract successively exceptional curves inside the reducible fibers
of $X_1$, making sure that the exceptional curves in question do not intersect 
$S_o$. One can use this procedure to reduce $X_1$ to a new ruled surface $X^{(2)}_1 $
which has a unique reducible fiber $F$ consisting of a union $C_1 \cup C_2 $ of two
smooth exceptional curves. Pick the curve, among $ C_1, C_2$, which
intersects $S_o$ and contract it. One obtains in this manner a projection 
$X^{(2)}_1 \rightarrow X^{(1)}_1 $ with $  X^{(1)}_1 $ 
geometrically ruled of type $\mathbb{F}_1 $. After contracting the image of 
$S_o$, we are left with $X^{(0)}_1$ which is a copy of $ \mathbb{P}^2 $. 
Sequences $ (\ref{bd11}) $ and $ (\ref{bd22}) $ determine a blowdown $ \rho \colon Z_o \rightarrow 
\mathbb{P}^2 \cup \mathbb{P}^2 $.     
\end{proof} 
\noindent Summarizing the facts, every elliptic Type II stable $K3$ surface with section of 
category (a) is isomorphic to a surface 
$ Z_o ( D;3p_0; x_p, p_2, \cdots, p_9; 3q_0; p_{10}, p_{11}, \cdots p_{18}) $ with  
$ \{ 3p_0; x_p, p_2, \cdots, p_9; 3q_0; p_{10}, p_{11}, \cdots p_{18} \} $ a special family on $D$. 
Every stable surface of category (b) can be transformed, after elementary modifications
to a surface $ Z_o ( D;3p_0; p_1, p_2, \cdots, p_{17}; 3q_0; p_{18} )$  associated to a special
family $ \{ 3p_0; p_1, p_2, \cdots, p_{17}; 3q_0; p_{18} \} $.

\subsection{Stable Periods and Torelli Theorem}
\label{sstableper}
\noindent We are now in position to provide the formal connection between Type II elliptic 
stable $K3$ surfaces
with section and Type II boundary points in the arithmetic partial compactification of 
$\Gamma \backslash \Omega  $. 
This correspondence will be later justified geometrically as an extended period map, 
using the theory of $K3$ degenerations. 
\begin{teo}
\label{t113}
Let $Z_o$ be an elliptic Type II $K3$ surface with section. Denote by $D$ the double curve. 
One can naturally associate to $Z_o$ a 
rank-sixteen unimodular even negative-definite lattice $\Lambda_{Z_o}$ together
with an abelian group homomorphism:
$$ \psi_{Z_o} \colon \Lambda_{Z_o} \rightarrow {\rm Jac}(D). $$
Moreover, $\Lambda_{Z_o}$ is a lattice of type $ -(E_8 \oplus E_8) $, for surfaces 
$Z_o$ in the (a)-category, and is of type $-\Gamma_{16} $, for $Z_o$ in the (b)-category.    
\end{teo}
\begin{proof}
\par We shall use a few known facts (see \cite{friedman} and \cite{friedman3}) concerning the Hodge 
theory of a Type II stable $K3$ surface. 
The rank of $ H^2(Z_o, \Zee ) $ is $21$. The complex cohomology group
$ H^2(Z_o, \Cee) $ carries a canonical mixed Hodge structure of weight filtration: 
%which one can read from 
%the Mayer-Vietoris exact sequence
%$$ 
%\cdots \rightarrow H^1(D, \Cee) \rightarrow H^2(Z_o, \Cee) \rightarrow H^2(X_1, \Cee) 
%\oplus H^2(X_2, \Cee) \rightarrow 
%H^2(D, \Cee) \rightarrow \cdots \ . 
%$$
%The weight filtration is:
\begin{equation}
\label{hodge666}
 0 \ \subset \ {\rm W}_1 \ \subset {\rm W}_2 \ = \  H^2(Z_o).
\end{equation} 
The two associated graded Hodge structures involved satisfy: 
$$ {\rm W}_1 \simeq H^1(D) \ \ {\rm (isomorphism \ of \ Hodge \ structures)}, $$
$$  {\rm W}_2 /   {\rm W}_1 \ \simeq \ {\rm Ker} \left ( 
H^2(X_1) \oplus H^2(X_2) \rightarrow H^2(D) 
\right ). $$
One deduces that $ {\rm W}_2 /   {\rm W}_1 $ has rank $19$ and carries a pure Hodge structure 
of type $ (1,1) $.   
\par The mixed structure on $ H^2(Z_o) $ produces an extension of mixed Hodge 
structures:
\begin{equation}
\label{zz1}
0 \rightarrow {\rm W}_1 \rightarrow {\rm W}_2 \rightarrow {\rm W}_2 /   {\rm W}_1 
\rightarrow 0 
\end{equation}
which, according to Carlson \cite{carl}, is classified by the associated abelian group homomorphism:
\begin{equation}
\label{extmap}
\w{\psi}_{Z_o} \colon \left ( {\rm W}_2 /   {\rm W}_1 \right )_{\Zee} \rightarrow 
%{\rm Jac}(E) \simeq E.  
J^1({\rm W}_1 ).
\end{equation}
Here  $J^1({\rm W}_1 ) =  W_1 / \ F^1W^1 + \left ( W_1 \right )_{\Zee} $ 
is the generalized Jacobian associated to the Hodge structure
on $ W_1 $. There is a purely geometrical description for $(\ref{extmap}) $. Since
the Hodge structure on $W_1$ is isomorphic to the geometrical weight-one Hodge structure of the
double curve $D$, one has a natural identification:
$$ J^1({\rm W}_1 ) \simeq {\rm Jac}(D) = {\rm Pic}^o(D). $$
Moreover, since the two surfaces $X_1$, $X_2$ are rational, any given cohomology class 
$$ [L] \in \left ( {\rm W}_2 /  {\rm W}_1 \right )_{\Zee} = {\rm Ker} \left ( 
H^2(X_1, \Zee) \oplus H^2(X_2, \Zee) \rightarrow H^2(D, \Zee) 
\right ) $$ is uniquely 
represented by a pair of holomorphic line bundles 
$ L = ( L_1 , L_2 ) \in {\rm Pic}(X_1) \times {\rm Pic}(X_2) $
satisfying $ L_1 \cdot D = L_2 \cdot D $. The image of $[L]$ under $(\ref{extmap}) $ can be
then described as:
\begin{equation}
\label{eext}
\w{\psi}_{Z_o} \left ( [L] \right ) \ = \ \m{O}_D(L_1) \otimes \m{O}_D(-L_2) \in {\rm Pic}^o(D) = 
{\rm Jac}(D) .  
\end{equation} 
In particular $ \w{\psi}_{Z_o} \left ( [L] \right ) = 0 $ for any cohomology class $ [L] $ representing 
a Cartier divisor on $ Z_o $. 
\par The lattice $\left ( {\rm W}_2 /  {\rm W}_1 \right )_{\Zee}$ has rank $19$ and 
is indefinite. However, the elliptic structure with section on $Z_o$ induces a series of 
Cartier divisors producing special cohomology classes. Firstly, the section $S_o$, 
which in case (a) is 
represented by two rational curves in $X_1$ and $X_2$ meeting $D$ at the same point, while  
in case (b) is a unique rational 
curve in $X_1$ disjoint from $D$, determines a Cartier divisor $\m{S}_o$ on $Z_o$. Secondly, the fiber
on $Z_o$, which in case (a) consists of elliptic fibers merging at $D$, while in case (b) consists of
rulings on each $X_i $ 
agreeing over the double curve, determines a  Cartier divisor class $ \m{F}_o $. Thirdly, let:
$$ \xi_1 = \m{O}_{X_1}(-D) \in {\rm Pic}(X_1), \ \ \ \xi_2 = \m{O}_{X_2}(D) \in {\rm Pic}(X_2) . $$
The d-stability condition assures us that the two line bundles agree over the double curve and
therefore they can be seen to determine a line bundle $ \xi_o $ over $Z_o$. The three Cartier 
divisors $\m{S}_o, \m{F}_o, \xi_o $ on $Z_o$ determine integral cohomology classes:
$$ [\m{S}_o], [\m{F}_o], [\xi_o] \in \left ( {\rm W}_2 /  {\rm W}_1 \right )_{\Zee} 
$$
satisfying 
$ [\m{S}_o]^2 = -2$, $[\m{F}_o]^2=0$, $[\xi_o]^2=0$, $[\m{S}_o].[\m{F}_o]=1 $,
$ [\m{S}_o].[\xi_o]=0$, $ [\m{F}_o].[\xi_o]=0$.
\par Denote by $ \{ [\xi_o]\}^{\perp} $ the sub-lattice of 
$ \left ( {\rm W}_2 /  {\rm W}_1 \right )_{\Zee} $ 
orthogonal to the class $ [\xi_o]$. Clearly 
all three elements $ [\xi_o] $ $ [\m{S}_o]$ and $ [\m{F}_o]$ belong to 
$  \{ [\xi_o]\}^{\perp} $. Then, define:
$$ \Lambda_{Z_o} 
%\w{\left ( {\rm W}_2 / \  {\rm W}_1 \right )} _{\Zee} 
\ \subset \{ [\xi_o]\}^{\perp} / 
\left (\Zee \cdot [\xi_o]  \right ) $$ 
as the sub-lattice orthogonal to the equivalence classes induced by $ [\m{S}_o]$ and $ [\m{F}_o]$. 
A simple observation involving the Hodge index theorem on $X_1 $ and $X_2 $ allows one to conclude
that 
$\Lambda_{Z_o}$  
%$ \w{\left ( {\rm W}_2 / \  {\rm W}_1 \right )} _{\Zee} $ 
is even, unimodular,
negative-definite and of rank $16$. As mentioned earlier, 
the extension homomorphism $ (\ref{eext}) $ vanishes on
cohomology classes representing Cartier divisors. In particular $ \w{\psi}_{Z_o} $ vanishes an all    
$ [\m{S}_o], [\m{F}_o] $ and $ [\xi_o] $. Therefore, without losing geometrical information, 
one can descend $ (\ref{eext}) $ to an abelian group homomorphism:
\begin{equation}
\label{eeeeext1112}
\psi_{Z_o} \colon 
%\w{\left ( {\rm W}_2 /   {\rm W}_1 \right )}_{\Zee}  
\Lambda_{Z_o} \rightarrow {\rm Jac}(D). 
\end{equation}  
\par The isomorphism type of the lattice $ \Lambda_{Z_o} $ is characterized by the category 
to which the stable surface $Z_o$ belongs. Assume that $Z_o $ is a surface in the (a)-category. There is
then a natural splitting $ \Lambda_{Z_o} =\Lambda^1_{Z_o} \oplus \Lambda^2_{Z_o} $ where
$$ \Lambda^i_{Z_o} \ = \ \{ \gamma \in H^2(X_i, \Zee) \ \vert \ \gamma \cdot [D] = 0, \ 
\gamma \cdot [S_i]=0 \ \}. $$
Pick a blow-down $ \rho \colon Z_o \rightarrow \mathbb{P}^2 \cup \mathbb{P}^2 $ as 
in Definition $ \ref{dedede} $ 
and consider the associated classes: 
$$ \{ H_1, H_2, E_1, \cdots E_{18} \} \subset H^2(X_1, \Zee) \oplus  H^2(X_1, \Zee) $$ 
representing the proper transforms of a hyper-planes 
in $ \mathbb{P}^2 $ and the total transforms of the exceptional 
curves associated to the blow-ups $ X^{(1)}_1 \rightarrow  X^{(0)}_1 $, 
$ X^{(2)}_1 \rightarrow  X^{(1)}_1 $, $ \cdots $ $ X^{(9)}_1 \rightarrow  X^{(8)}_1 $, 
$ X^{(1)}_2 \rightarrow  X^{(0)}_2 $, $ X^{(2)}_2 \rightarrow  X^{(1)}_2 $, $ \cdots$ $ 
X^{(9)}_2 \rightarrow  X^{(8)}_2 $. Let $ \alpha_1, \alpha_2, \cdots \alpha_8, 
\beta_1, \beta_2, \cdots \beta_8  $ be the following sixteen elements in $\Lambda_{Z_o}$:
\begin{equation}
\label{rootsyse8}
\alpha_1= E_1-E_2, \  \alpha_2= E_2-E_3,  \ \cdots , \ \alpha_7=E_7-E_8,  \  
\alpha_8=H_1-E_1-E_2-E_3  
\end{equation}   
$$
\beta_1=E_{10}-E_{11}, \  \beta_2=E_{11}-E_{12}, \ \cdots , \ 
\beta_7= E_{16}-E_{17}, \  \beta_8=H_2-E_{10}-E_{11}-E_{12}.  
$$  
One verifies that  $ \{ \alpha_1, \alpha_2, \cdots \alpha_8 \} $ and $ \{ \beta_1, \beta_2, \cdots \beta_8 \} $
are basis for $ \Lambda^1_{Z_o} $ and $ \Lambda^2_{Z_o} $. Moreover, analyzing the intersection numbers, one 
finds out that, after changing the sign of the quadratic pairing, each of the two lines in 
$ (\ref{rootsyse8}) $ consists of a set of  $E_8$ simple roots.  
$$ 
\def\objectstyle{\scriptstyle}
\def\labelstyle{\scriptstyle}
\xymatrix @-1.2pc  {
\stackrel{\alpha_1}{\bullet} \ar @{-} [r] &
\stackrel{\alpha_2}{\bullet} \ar @{-} [r] &
\stackrel{\alpha_3}{\bullet} \ar @{-} [r] \ar @{-} [d] &
\stackrel{\alpha_4}{\bullet} \ar @{-} [r] &
\stackrel{\alpha_5}{\bullet} \ar @{-} [r] &
\stackrel{\alpha_6}{\bullet} \ar @{-} [r] &
\stackrel{\alpha_7}{\bullet}   \\
 & & \alpha_8 \bullet & & & & \\ 
} \ \ \ \ 
\xymatrix @-1.2pc  {
\stackrel{\beta_1}{\bullet} \ar @{-} [r] &
\stackrel{\beta_2}{\bullet} \ar @{-} [r] &
\stackrel{\beta_3}{\bullet} \ar @{-} [r] \ar @{-} [d] &
\stackrel{\beta_4}{\bullet} \ar @{-} [r] &
\stackrel{\beta_5}{\bullet} \ar @{-} [r] &
\stackrel{\beta_6}{\bullet} \ar @{-} [r] &
\stackrel{\beta_7}{\bullet}   \\
 & & \beta_8 \bullet & & & & \\ 
}.
$$
\noindent The lattice 
$ \Lambda_{Z_o} $ is therefore isomorphic to $ - \left ( E_8 \oplus E_8 \right ) $. 
%Moreover, in this setting
%one can easily describe the homomorphism  $ (\ref{eext}) $. Given the group structure on $E$ with $p_o$ 
%as origin and the identification $ E \simeq D$, there is a natural group isomorphism 
%$ {\rm Pic}^o(D) = E $. Under this identification, the behavior of $ (\ref{eext}) $ on the basis elements
%of $(\ref{rootsyse8}) $ is given by:
%\begin{equation}
%\psi_{Z_o} (\alpha_1)= p_1-p_2, \  \psi_{Z_o} (\alpha_2)=p_2-p_3,  \ \cdots , \ 
%\psi_{Z_o} (\alpha_7)=p_7-p_8,  \  \psi_{Z_o} (\alpha_8)=p_0-p_1-p_2-p_3 
%\end{equation}
\par One can do a similar analysis in the case when $ Z_o $ is in category (b). 
Note that the isomorphism class of the 
pair $ \left (\Lambda_{Z_o} , \psi_{Z_o} \right ) $ does not change under elementary modifications. 
Indeed, a modification that flips an exceptional 
curve $C$ from $X_1$ to $X_2$ induces an isometry:
\begin{equation}
\label{modifizo}
\xymatrix{
H^2(X_1, \Zee) \oplus H^2(X_2, \Zee) \ar [r] \ar [d] _{\simeq}  & H^2(X'_1 , \Zee) \oplus H^2(X'_2 , \Zee) 
\ar [d] ^{\simeq} \\
H^2(X'_1, \Zee) \oplus \Zee [C] \oplus H^2(X_2, \Zee ) \ar [r] 
%_{ ({\rm id}) \oplus (-{\rm id}) \oplus ({\rm id})} 
& H^2(X'_1, \Zee) \oplus \Zee [C] \oplus H^2(X_2, \Zee). 
}
\end{equation}  
This map sends $ [C] \in H^2(X_1, \Zee) $ to $-[C] \in H^2(X'_2, \Zee) $, 
$ [\m{F}_o] $ to $ [\m{F}'_o] $, $ [\m{S}_o] $ to $ [\m{F}'_o] $ and $ [\xi_o] $ to $ [\xi'_o] $.  
There is then an induced lattice isomorphism $ \Lambda_{Z_o} \simeq   \Lambda_{Z'_o} $ which clearly makes the
diagram:
$$ 
\xymatrix @-1.2pc  {
\Lambda_{Z_o}  \ar [rr] ^{\simeq} \ar [dr] _{\psi_{Z_o}} & & \Lambda_{Z'_o} \ar [dl] ^{\psi_{Z'_o}}\\
 & {\rm Jac}(D) & \\
}
$$
commutative.
\par According to Proposition $ \ref{teooob}$, $ Z_o$ can be transformed, using elementary modifications, 
such that the resulting surface $Z'_o$ admits a blow-down $ \rho \colon Z'_o \rightarrow \mathbb{P}^2 \cup \mathbb{P}^2 $ 
associated to a special family $ \{ 3p_0; p_1, p_2, \cdots p_{17}; 3q_0; p_{18} \} $ on $D$. In such conditions,
a basis $ \{ \gamma_1, \gamma_2, \cdots \gamma_{16} \} $ for $ \Lambda_{Z'_o} $ is given by: 
$$ \gamma_1 = H_1 -E_1-E_2-E_{3}, $$ 
\begin{equation}
\label{rootsysd16}
\gamma_2= E_{3}-E_{4}, \ \gamma_3 = E_4 -E_5 , \ \cdots  , 
\gamma_{14} = E_{15} -E_{16} 
\end{equation}
$$ \gamma_{15} = E_{16} -E_{17}, \  
\gamma_{16} = H_2-E_{18}+E_{16}+E_{17}. $$
One verifies that, after reversing the sign of the pairing, $ (\ref{rootsysd16}) $ is a root system of type $D_{16} $. 
$$ 
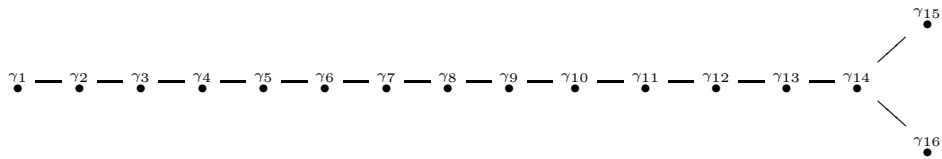
\def\objectstyle{\scriptstyle}
\def\labelstyle{\scriptstyle}
\xymatrix @-1.2pc  {
 & & & & & & & & & & & & & & \stackrel{\gamma_{15}}{\bullet} \ar @{-} [dl] \\
\stackrel{\gamma_{1}}{\bullet} \ar @{-} [r] &
\stackrel{\gamma_{2}}{\bullet} \ar @{-} [r] &
\stackrel{\gamma_{3}}{\bullet} \ar @{-} [r] &
\stackrel{\gamma_{4}}{\bullet} \ar @{-} [r] &
\stackrel{\gamma_{5}}{\bullet} \ar @{-} [r] &
\stackrel{\gamma_{6}}{\bullet} \ar @{-} [r] &
\stackrel{\gamma_7}{\bullet} \ar @{-} [r] &
\stackrel{\gamma_8}{\bullet} \ar @{-} [r] &
\stackrel{\gamma_9}{\bullet} \ar @{-} [r] &
\stackrel{\gamma_{10}}{\bullet} \ar @{-} [r] &
\stackrel{\gamma_{11}}{\bullet} \ar @{-} [r] &
\stackrel{\gamma_{12}}{\bullet} \ar @{-} [r] &
\stackrel{\gamma_{13}}{\bullet} \ar @{-} [r] &
\stackrel{\gamma_{14}}{\bullet} \ar @{-} [dr] & \\
 & & & & & & & & & & &  & & & \stackrel{\gamma_{16}}{\bullet} \\
}
$$
The lattice $ \Lambda_{Z_o}$ is therefore isomorphic to $ - \Gamma_{16} $.      
\end{proof} 
\noindent We make now the connection with the arithmetic Mumford boundary points. 
In the notation of $\ref{arith} $, assume that $F$ is a Type II Bailey-Borel component for 
$ \Gamma \backslash \Omega $, 
corresponding to the isotropic rank-two sub-lattice $ V \subset \mathbb{L}_o $, and $\m{B}(F)$ is the
associated Type II Mumford boundary divisor. 
Let $ \Lambda $ be the rank-sixteen lattice $ V^{\perp}/V$.
Recall from $ \ref{p} $ 
that $\m{B}(F)$ decomposes into two connected components 
$$\m{B}^+(F) \sqcup \overline{\m{B}^+(F)} $$ 
and there is a bijective identification between boundary points in $ \m{B}^+(F) $ 
and pairs $ (\m{H}, \psi) $ 
consisting of weight-one Hodge structures on $V_{\Cee} $ polarized with respect to the 
skew-symmetric form
$(\cdot, \cdot )_1 $ together with 
abelian group homomorphisms $ \psi \colon \Lambda \rightarrow J^1(\m{H}) $ . 
\par Let $Z_o$ be an elliptic Type II stable $K3$ surface with section as in definition 
$ \ref{dddd33} $ (a)-(b). Attach 
to $Z_o$ a set of markings $ \phi_1, \ \phi_2 $ consisting of isometries     
\begin{equation}
\label{stabmark}
\phi_1 \colon H^1(D, \Zee) \rightarrow V, \ \ \ \ \phi_2 \colon \Lambda_{Z_o} \rightarrow \Lambda.    
\end{equation}
The marking $ \phi_1$ can be used to transport
the geometrical weight-one Hodge structure $W_1$ of $D$ to a formal weight-one 
polarized Hodge structure $\m{H}$ on $V$. There is then an induced isomorphism of abelian groups 
$ {\rm Jac}(D) \simeq J^1(\m{H}) $. This isomorphism, together with the marking $\phi_2$, allows one 
to transport the homomorphism $ \psi_{Z_o} $ of Theorem $ \ref{t113} $ to a formal homomorphism 
$ \psi \colon \Lambda \rightarrow J^1(\m{H}) $. In the light of the arguments in previous paragraph, 
this procedure can be regarded as a period correspondence, associating to every 
marked elliptic Type II stable $K3$ surface with section $(Z_o, \phi_1, \phi_2) $ 
in $\Lambda$-category a ${\bf marked \ stable \ period} $ in the form of a pair 
$ (\m{H}, \psi) \in \m{B}^+(F) $.    
\par One can further refine this correspondence by removing the markings and
considering the pairs $ (\m{H}, \psi) $ modulo the isometries of $ \Lambda $ and $V$.
Let us denote by $ {\rm SL}_2(\Zee) $ the group of automorphisms of $ V$ preserving the
skew-symmetric pairing $ (\cdot, \cdot)_1 $. This group acts naturally on the set of
weight-one Hodge structures on $V$ polarized with respect to $(\cdot, \cdot)_1 $. Consider
${\rm Aut}(\Lambda)$ to be the group of isometries of lattice $ \Lambda $. The product
group $ \m{G} = {\rm Aut}(\Lambda) \times {\rm SL}_2(\Zee) $ acts then on the set of pairs 
$ (\m{H}, \psi) $ as:
$$ (f,\alpha) . (\m{H}, \psi) \ = \ \left ( 
\alpha(\m{H}), \w{\alpha} \circ \psi \circ f^{-1}
\right ) $$
where $ \w{\alpha} \colon J^1(\m{H}) \rightarrow J^1(\alpha(\m{H})) $ is the natural 
isomorphism induced by $\alpha $. It is clear that, 
given a marked triplet $ (Z_o, \phi_1, \phi_2) $ inducing a marked pair $(\m{H}, \psi) $, 
a variation of markings $ \phi_1, \phi_2 $ or a change 
of $Z_o$ under an isomorphisms as in $ (\ref{iiizzzooo}) $ leaves $(\m{H}, \psi) $ 
within the same $ \m{G} $-orbit. Therefore, one can associate to any elliptic Type II
stable $K3$ surface with section a well-defined $ {\bf stable \ period} $ in 
$ \m{G} \backslash \m{B}^+(F) $.    
\begin{dfn}
Two elliptic Type II stable $K3$ surfaces with section in the same category:
$$ Z_o = X_1 \cup X_2 \ \ {\rm and} \ \  Z'_o = X'_1 \cup X'_2 $$
are said to be $ {\bf equivalent} $ if one of the following holds:
\begin{enumerate}
\item $ Z_o $ and $Z'_o$ are isomorphic (as in $ (\ref{iiizzzooo}) $). 
\item $ Z_o $ and $Z'_o$ are both of category (a) and $Z_o$ is isomorphic to $X'_2 \cup X'_1$.
\item $ Z_o $ and $Z'_o$ are both of category (b) and can be made to be isomorphic by 
transforming each of them using a finite sequence of elementary modifications.
\end{enumerate} 
\end{dfn} 
\noindent Let then 
$ \m{M}^{{\rm stable}}_{\Lambda} $, $ \Lambda = E_8 \oplus E_8 $ or $ \Gamma_{16} $, be 
the coarse moduli spaces of equivalence classes in category (a), respective (b). It can be
easily seen that the stable period of a surface $Z_o=X_1 \cup X_2$ does not change 
when $Z_o$ gets replaced by $ X_2 \cup X_1 $ (if $Z_o$ is of category (a) ) or when
$Z_o$ gets transformed by an elementary modification. One has therefore a well-defined 
period map:
\begin{equation}
\label{per22}
{\rm per}_{\Lambda} \colon \m{M}^{{\rm stable}}_{\Lambda} \rightarrow \m{G} \backslash \m{B}^+(F).  
\end{equation} 
Furthermore, as we shall see from the analysis in section $ \ref{ellipticflat} $, there exists 
a natural group isomorphism $ \m{G} \simeq \Gamma^+_F = P(F) \cup \Gamma^+ $ (recall that $P(F)$ 
is the rational parabolic subgroup associated to the Baily-Borel boundary component $F$). Moreover,
under this isomorphism, the action of $ \m{G} $ on $ \m{B}^+(F) $ reduces to the standard arithmetic 
action of $P(F) \cup \Gamma^+$. This produces a natural identification:
$$ \m{G} \backslash \m{B}^+(F) \simeq \Gamma^+ \ \m{B}^+(F) = \Gamma_F \backslash \m{B}(F) = 
\m{D}_{\Lambda}. $$
One can therefore interpret the stable period of an elliptic Type II stable $K3$ surface with section 
as a point of the arithmetic Mumford divisor $\m{D}_{\Lambda} $ and hence regard $ (\ref{per22}) $ as 
a map 
$ {\rm per}_{\Lambda} \colon \m{M}^{{\rm stable}}_{\Lambda} \rightarrow \m{D}_{\Lambda} $.        
\begin{teo}
\label{late}
The map $ (\ref{per22}) $ is an isomorphism.
\end{teo}
\noindent We prove this statement in two steps. To begin with, let us show that $ (\ref{per22}) $ is injective. 
\begin{teo}
\label{ainject}
Two stable surfaces $Z_o$ and $Z_o$ of category (a), which have the same stable period, 
are equivalent.
\end{teo} 
\begin{proof}
This follows from standard results concerning $E_8$ del Pezzo surfaces 
(see \cite{delpezzo} \cite{demazure} and \cite{manin} for details). If 
$Z_o = X_1 \cup X_2 $ and $Z_o' = X_1' \cup X_2' $ are stable surfaces 
of category (a), then, after contracting the sections, one obtains four $E_8$ del Pezzo surfaces 
$ \w{X}_1 , \ \w{X}_2, \ \w{X}'_1, \ \w{X}'_2 $. Moreover, one has isomorphisms:
$$ \Lambda_{Z_o} \ = \ \Lambda^1_{Z_o} \oplus \Lambda^2_{Z_o} \ \simeq \ 
[K_{\w{X}_1} ] ^{\perp} \oplus [K_{\w{X}_2} ] ^{\perp} \subset  
H^2(\w{X}_1, \Zee) \oplus H^2(\w{X}_2, \Zee)$$    
$$ \Lambda_{Z'_o} \ = \ \Lambda^1_{Z'_o} \oplus \Lambda^2_{Z'_o} \ \simeq \ 
[K_{\w{X}'_1} ] ^{\perp} \oplus [K_{\w{X}'_2} ] ^{\perp} \subset  
H^2(\w{X}_1', \Zee) \oplus H^2(\w{X}_2', \Zee). $$  
It was proved in \cite{delpezzo} that the isomorphism class of a pair $ (\w{X}, D) $, 
consisting of an $E_8$ del Pezzo surface $\w{X}$ with an embedded smooth elliptic curve 
$D$, is determined by the map $ [K_{\w{X}} ] ^{\perp} \rightarrow {\rm Jac}(D) $ 
modulo Weyl equivalence. Based on this argument, assuming that $Z_o$ and $Z_o'$ determine the same stable 
period in $\m{G} \backslash  \m{B}^+(F) $, it follows that there is an isomorphism of 
elliptic curves $D \simeq D'$ which extends to 
an isomorphism of stable surfaces of either $ X_1 \cup X_2 \simeq X'_1 \cup X'_2 $ form 
or $ X_1 \cup X_2 \simeq X'_2 \cup X'_1 $ form.            
\end{proof} 
\noindent We use different arguments for justifying the analog of Theorem $ \ref{ainject} $ for 
stable surfaces of category (b). As shown earlier, given an elliptic Type II stable $K3$ surface with section 
$Z_o= X_1 \cup X_2 $ of category (b), one can always transform $Z_o$, by performing 
elementary modifications, to a new stable surface $Z'_o $ such that $ X'_2 \simeq \mathbb{F}_1 $. 
In this setting, there exists blowdowns 
$ \rho \colon Z'_o \rightarrow \mathbb{P}^2 \cup \mathbb{P}^2 $ and
each choice of such blowdown induces a $D_{16} $ simple root system 
$$ \{ \gamma_1, \ \gamma_2, \ \cdots \gamma_{16} \}$$
for $ \Lambda_{Z'_o} $, as described in $ (\ref{rootsysd16}) $. The model $ Z'_o $ and 
the blowdown $ \rho$ are far from being unique. One can 
further transform $Z_o$, using sequences of elementary modifications, to new surfaces $Z''_o $, 
satisfying $ X''_2 \simeq \mathbb{F}_1 $, but not isomorphic to $Z'_o$. However, any modification 
from $Z_o$ to $ Z'_o $ induces a canonical isomorphism $ \Upsilon \colon 
\Lambda_{Z_o} \rightarrow  \Lambda_{Z'_o} $ 
(see $(\ref{modifizo})$) entering the commutative diagram: 
$$ 
\xymatrix @-1.2pc  {
\Lambda_{Z_o}  \ar [rr] ^{\Upsilon} \ar [dr] _{\psi_{Z_o}} & & \Lambda_{Z'_o} \ar [dl] ^{\psi_{Z'_o}}\\
 & {\rm Jac}(D). & \\
}
$$
\begin{lem}
\label{techforb}
Let $Z_o$ be an elliptic Type II stable surface with section, of category (b). For any basis of 
of simple roots $ \mathbb{S} \subset \Lambda_{Z_o} $, there exists a sequence of elementary modifications 
transforming $Z_o$ to a new stable surface $Z'_o=X'_1 \cup X'_2 $ with $ X'_2 \simeq \mathbb{F}_1 $ 
and a blowdown $ \rho \colon Z'_o \rightarrow \mathbb{P}^2 \cup \mathbb{P}^2 $ such that the simple root 
system associated to $ \rho $ is $ \Upsilon \left ( \mathbb{S} \right ) $.      
\end{lem}
\begin{proof}
Any two sets of $D_{16} $ simple roots can be transformed one into the other using a Weyl transformation. 
It suffices then to show that, given $Z_o $ with $X_2 \simeq \mathbb{F}_1 $ and fixing a blowdown 
$ \rho_0 \colon Z_o \rightarrow \mathbb{P}^2 \cup \mathbb{P}^2 $ with associated set of simple roots
$ \mathbb{S}_0 $, for any Weyl transformation $ w \in W(\Lambda_{Z_o}) $, there exists a sequence 
of elementary modifications transforming $Z_o$ to $Z'_o$ with $X'_2 \simeq \mathbb{F}_1 $ and a blowdown 
$ \rho \colon Z'_o \rightarrow \mathbb{P}^2 \cup \mathbb{P}^2 $, such that the simple root set associated
to $\rho $ is $ \Upsilon \left ( w \cdot \mathbb{S}_0 \right ) $.
\par Let $ H_1, E_1, \cdots E_{17} $ and $H_2 , E_{18} $ be the hyper-plane sections and the total transforms
of the exceptional curves associated to $ \rho_0 $. The simple root set $ \mathbb{S}_0 $ is:        
%$$ \gamma_1 = H_1 -E_1-E_2-E_{3}, $$ 
\begin{equation}
\label{rootsysd16d}
\gamma_1 = H_1 -E_1-E_2-E_{3}, \ \gamma_2= E_{3}-E_{4}, \ \gamma_3 = E_4 -E_5 , \ \cdots  , 
\gamma_{14} = E_{15} -E_{16} 
\end{equation}
$$ \gamma_{15} = E_{16} -E_{17}, \  
\gamma_{16} = H_2-E_{18}+E_{16}+E_{17}. $$
\noindent We define $ \varepsilon_1 , \varepsilon_2 \cdots , \varepsilon _{16} $, elements of 
$ \Lambda_{Z_o} \otimes \mathbb{Q} $ given by:
$$ \varepsilon_1 \ = \frac{1}{2} \left ( H_2 - E_{18} \right ) + H_1 -E_1 -E_2, \ \ 
\varepsilon_l = \frac{1}{2} \left ( H_2 - E_{18} \right ) + E_{l+1}, \ \ 2 \leq l \leq 16 . $$
%$$ \varepsilon_1 \ = \frac{1}{2} \left ( H_2 - E_{18} \right ) + H_1 -E_1 -E_2.   $$ 
The set 
$ \{ \varepsilon_1 , \varepsilon_2 \cdots , \varepsilon _{16} \} $ forms an orthonormal basis for
$ \Lambda^1_{Z_o} \otimes \mathbb{Q} $ (when changing the sign of the quadratic form) and the
roots in $ \mathbb{S}_0 $ appear as:
$$ \gamma_l = \varepsilon_l - \varepsilon_{l+1} , \ {\rm for} \  1 \leq l \leq 15 , \ {\rm and} \ 
\gamma_{16} = \varepsilon_{15} + \varepsilon_{16}. $$ 
%$$ \gamma_{16} = \varepsilon_{15} + \varepsilon_{16} .$$  
In this setting, it is known that the Weyl 
group $ W(\Lambda_{Z_o}) $ is generated by permutations of 
$ \varepsilon_1 , \varepsilon_2 \cdots , \varepsilon _{16} $ and transformations 
$ t_{ij} $, $ 1 \leq i < j \leq 16 $ 
taking $ \varepsilon_i, \varepsilon_j $ to $ -\varepsilon_i, -\varepsilon_j $ and leaving 
all other $ \varepsilon_l $ unchanged. In what follows, we shall indicate the 
elementary modifications and the change in blowdown sequence generating, at the 
level of roots, transpositions $ ( \varepsilon_i, \varepsilon_j) $. One can use a similar 
technique to treat the transformations $ t_{ij} $.
\par We shall denote by $ X_1^{(l)} $ the surfaces obtained from $X_1 $ during the blowdown $\rho_0 $, 
by $E_l $ the corresponding contracting curves, and by $p_l $ the intersection points $D \cap E_l $. 
Start with $X_1 = X_1^{(17) } $ and contract 
successively $ E_{17}, \cdots E_{j+2} $. The resulting surface is $X_1^{(j+1)} $. The total 
transform of $E_{i+1} $
on $X_1^{(j+1)} $ is a chain $ C_1 \cup C_2 \cup \cdots \cup C_k $ of smooth rational curves with 
self-intersection $-2$, with the exception of $C_1$ which is exceptional. One has intersecting numbers 
$ C_l \cdot C_{l+1} = 1 $ for  
$ 1 \leq l  \leq k-1 $ and $C_l \cdot C_l' =0 $ otherwise. Flip $ C_1, C_2, \cdots C_{k-1} $, 
successively, to $ X_2 $ and then contract $ C_k $. Then flip $ C_{k-1}, C_{k-2}, \cdots C_2 $ back. 
Flip $ E_{j+1} $ to the right. Denote the resulting 
stable surface by $ \w{X}^{j}_1 \cup \w{X}_2 $. Next, if $ i \geq 2 $ then contract successively the 
curves $E_{j}, E_{j-1} , \cdots E_{i+2} $ on  $\w{X}^{j}_1 $. Let the resulting surface be denoted 
$\w{X}^{i+1}_1 $. Flip $E_{j+1}$ back from $\w{X}_2$ and contract it on  $\w{X}^{i+1}_1 $. 
The resulting surface is exactly $ X^{(i)}_1 \cup X_2 $. Keep then the rest of the blowdown intact 
and construct the upper part of the new blowdown $\rho $ by retracing the steps and blowing up 
successively the points $ p_{j+1}, p_{i+2}, p_{i+3} , \cdots , p_{j} , p_{i+1} , p_{j+2}, \cdots p_{17} $, 
on  $X^{(i)}_1$.
$$
\def\objectstyle{\scriptstyle}
\def\labelstyle{\scriptstyle}
\xymatrix @-0.5pc  {
%\xymatrix{
X_1 = X^{(17)}_1 \ar [r] & \cdots \ar [r] & X^{(j+1)}_1 \ar [d] & \\
                         &                & \w{X}^{(j+1)}_1 \ar [r] & X^{(j)}_1 \ar [r] & \cdots \ar [r] & 
X^{(i+1)}_1 \ar [d] \\ 
X_1'= X'^{(17)}_1  \ar [r]        &  \cdots \ar [r]  & X'^{(j+1)}_1  \ar [r]  &  X'^{(j)}_1 \ar [r]  &
    \cdots \ar [r]         & \w{X} ^{(i+1)}_1 \ar [r] &  X ^{(i)}_1 \ar [r] &   X ^{(i-1)}_1 \ar [r] & 
\cdots \ar [r] & X ^{(0)}_1 \simeq \mathbb{P}^2 \\
}
$$
The new stable surface $ Z_o' = X_1' \cup X_2 $ is obtained from $Z_o $ through vcb a sequence of elementary 
modifications and the simple root system associated to the blowdown $ \rho $ 
is $ \Upsilon \left ( w \cdot \mathbb{S}_0 \right ) $.  
\par The case i = 1 requires a slight modification of the above procedure. After obtaining $ \w{X}^{j+1}_1$ 
continue by contracting the proper transform of the line passing 
through $p_1$ with multiplicity two. Flip $ E_{j+1} $ to $X_2$. Denote the resulting 
stable surface by $ \w{X}^{j}_1 \cup \w{X}_2 $. Contract successively $E_{j}, E_{j-1} , \cdots E_{3} $ on 
$\w{X}^{j}_1$. %The resulting surface coincides with the $\mathbb{F}_2 $ surface that one obtains after 
%contracting in $ X_1^{(2)} $ the proper transform of the line passing through $p_1$ with multiplicity two. 
Flip $ E_{j+1} $ back to the left. Contract the image of the proper transform of the line passing through 
$p_1$ and $ p_{l+1} $ and then contract the image of the proper transform of the line passing through 
$p_1$ and $ p_{l+1} $. The resulting 
surface is a copy of $ \mathbb{P}^2 $.
\end{proof} 
\noindent We are then in position to justify the injectivity of the stable period map $ (\ref{per22}) $ for
category (b) surfaces. 
\begin{teo}
\label{binject}
Two stable surfaces $ Z_o$ and $Z'_o$ of category (b), which have the same stable period, are equivalent. 
\end{teo}
\begin{proof}
Since elementary modifications do not vary the stable period, we can assume that both $X_2$ and $X'_2$ 
are copies of $ \mathbb{F}_1 $. Choose a blowdown 
$ \rho \colon Z_o \rightarrow \mathbb{P}^2 \cup \mathbb{P}^2 $ and denote by 
$ \mathbb{S} \subset \Lambda_{Z_o} $ the associated
basis of simple roots. 
\par Let $ (\phi_1, \phi_2) $ , $ ( \phi'_1, \phi_2 ) $ markings for $Z_o$, $Z'_o$ as in
$ (\ref{stabmark}) $. Denote by $ (\m{H}, \psi) $, $ (\m{H}', \psi') $ the induced marked
periods. Since the stable periods of the two surfaces are identical, there must exist 
isometries $ \alpha $ and $ f $ for $ V$ and $ \Lambda$, respectively, such that:
$$ \m{H}' = \w{\alpha}(\m{H} ) , \ \ \ \psi' = \w{\alpha} \circ \psi \circ f^{-1}. $$
Let $ \mathbb{S}' = \left ( (\phi'_2)^{-1} \circ f \circ \phi_2 \right )  (\mathbb{S}) $. Then $ \mathbb{S}' $
is a basis of simple roots in $ \Lambda_{Z'_o} $ and, according to Lemma $ \ref{techforb} $, 
there exists a new stable surface $Z''_o $, obtained from $Z'_o $ through a sequence of elementary 
modifications, which admits a blowdown $ \rho'' \colon Z''_o \rightarrow \mathbb{P}^2 \cup \mathbb{P}^2 $
such that the simple root basis $ \mathbb{S}'' $ associated to $ \rho''$ satisfies 
$ \mathbb{S}'' = \Upsilon \left ( \mathbb{S}' \right ) $.
%$$
%\xymatrix{
%\Lambda_{Z_o} \ar [d] _{(\phi_2)^{-1} \circ f \circ \phi_1} \ar [dr] ^{\phi_1} & & J^1(\m{H}) \ar [dd] ^{\w{\alpha}} \\  
%\Lambda_{Z'_o} \ar [r] ^{\phi'_2} \ar [d] _{\Upsilon} & \Lambda \ar [ur] ^{\psi} \ar@/_/ [dr] _{\psi''} 
%\ar@/^/ [dr] ^{\psi'} & \\ 
%\Lambda_{Z'_o} \ar [ur] %^{\phi_2 \circ \Upsilon^{-1} } 
%& & J^1(\m{H}') \\
%}
%$$       
\noindent Let then $ \left (  3p_0;p_1, p_2, p_3, \cdots, p_{16}, p_{17}; 3q_0; p_{18} \right ) $ and 
$\left (  3p''_0;p''_1, p''_2, p''_3, \cdots, p''_{16}, p''_{17}; 3q''_0; p''_{18} \right )  $ be 
the two special families on $D$ and $D'$ induced by the blowdowns $ \rho $ and $ \rho''$, respectively.
%$$ H_1, E_1 , E_2, \cdots E_{17}, H_2, E_{18}, \ \ \ \  H''_1, E''_1 , E''_2, \cdots E''_{17}, H''_2, E''_{18} $$
%be the hyper-plane sections and exceptional curves associated to the two blow-downs $ \rho $ and $ \rho'' $
%and denote by 
%\begin{equation}
%\label{fam11}
%\left (  3p_0;p_1, p_2, p_3, \cdots, p_{16}, p_{17}; 3q_0; p_{18} \right ) 
%\end{equation}           
%\begin{equation}
%\label{fam22}
%\left (  3p''_0;p''_1, p''_2, p''_3, \cdots, p''_{16}, p''_{17}; 3q''_0; p''_{18} \right ) 
%\end{equation}  
%be the two special families on $D$ and $D'$ induced by the blowdowns $ \rho $ and $ \rho''$, respectively.. 
Fix base points on $D$ and $D'$ and consider the induced identifications:
$$ D \simeq {\rm Jac}(D) = J^1(\m{H}) , \ \ \ D' \simeq {\rm Jac}(D') = J^1(\m{H}'). $$
Use these identifications to define an abelian group isomorphism:
\begin{equation}
\xymatrix{
\w{\eta} \colon \  D \ar [r] ^{\simeq} %\ar@/_/ [rrr] _{\w{\eta}} 
&  J^1(\m{H}) \ar [r] ^{\w{\alpha}} &  J^1(\m{H}') \ar [r] ^{\simeq} & D' . 
}
\end{equation} 
and then construct $ \eta \colon D \rightarrow D' $ with $ \eta(p) = \w{\eta}(p) - \w{\eta}(p_1) + p''_1 $. 
It turns out then that the isomorphism $ \eta $ transports the special family 
$ \left (  3p_0;p_1, p_2, p_3, \cdots, p_{16}, p_{17}; 3q_0; p_{18} \right ) $ to
$\left (  3p''_0;p''_1, p''_2, p''_3, \cdots, p''_{16}, p''_{17}; 3q''_0; p''_{18} \right )  $.
%$ (\ref{fam11}) $ to $ (\ref{fam22}) $. 
This implies that $Z_o$ and $Z''_o$ are isomorphic which, in turn, implies 
that $Z_o$ and $Z'_o$ are equivalent. 
\end{proof} 
\noindent One concludes from Theorems $ \ref{ainject} $ and $ \ref{binject} $ that the stable 
period map:
\begin{equation}
\label{per3334}
{\rm per}_{\Lambda} \colon \m{M}^{{\rm stable}}_{\Lambda} \rightarrow \m{G} \backslash \m{B}^+(F) 
\end{equation}  
is injective. Let us then complete the proof of Theorem $ \ref{late} $:
\begin{teo}
\label{p331}
The period map $ (\ref{per3334}) $ is surjective. 
\end{teo}
\begin{proof}
Let $ (\m{H}, \psi) $ be a pair in $ \m{B}^+(F) $. We show that there exists a marked surface 
$(Z_o, \phi_1,\phi_2)$ with stable marked period $ (\m{H}, \psi)$. The Hodge-theoretic Jacobian 
$ J^1(\m{H}) $ is itself a pointed elliptic curve endowed with a natural group 
structure. We agree to call it $ (E, p_0)$ and denote 
by $ \phi_1 \colon H^1(E, \Zee) \simeq V $ a marking that sends the geometrical Hodge structure 
of $E$ to $\m{H} $. In particular, $ \phi_2$  induces a group isomorphism 
\begin{equation}
\label{asdsa}
E \simeq {\rm Jac}(E) \simeq J^1(\m{H}) . 
\end{equation}     
\par If $ \Lambda = E_8 \oplus E_8 $, pick a basis for 
$ \{a_1, \cdots a_8, b_1, \cdots b_8 \}$ for $\Lambda $ such that 
$  \{a_1, \cdots a_8\}$ and $ \{b_1, \cdots b_8 \}$ are $E_8$ systems of simple roots. In what follows, we
construct $19$ points on $E$, denoted $ q_0, x_1, x_2 \cdots x_9 , y_1, y_2 \cdots y_{9} $. Choose
$ x_1 \in E $ such that:
$$ 3x_1 = 2 \psi (a_1) + \psi(a_2) - \psi(a_8). $$ 
Then construct $ p_2, \cdots p_9 $, recursively, by the rule:
$$ x_l = x_{l-1} - \psi(a_{l-1}), \ {\rm for} \ 2 \leq l \leq 8 $$
$$ x_9 =  - ( x_1 + x_2 + \cdots + x_8 ). $$ 
Then set $ y_9 = x_9 $ and:
$$ y_1 = y_9  + \left ( 7\psi(b_1)+ 6\psi(b_2) + 5\psi(b_3) + \cdots + 2 \psi(b_6) +\psi(b_7) \right ) \ - \ 
3 \left (2 \psi (b_1) + \psi(b_2) - \psi(b_8) \right ).    $$ 
Construct then recursively $ y_l = y_{l-1} + \psi(b_{l-1}) $ for $ 2 \leq l \leq 8 $ and then pick 
$ q_0 \in E $ such that: 
$$ 3q_0 = 2 \psi (b_1) + \psi(b_2) - \psi(b_8) + 3 y_1. $$
One verifies that $ \left (  3p_0; x_1, x_2, \cdots x_{9}; 3q_0; y_{1}, y_{2}, \cdots y_{9} \right )  $ 
is a special family of category (a) on $ E $. The stable surface 
\begin{equation}
\label{exxa}
Z_o \ = \ Z_o \left ( E; 3p_0; x_1, x_2, \cdots x_{9}; 3q_0; y_{1}, y_{2}, \cdots y_{9} \right ) 
\end{equation}
is then an elliptic Type II stable $ K3$ surface with section, of category (a). Moreover, 
$Z_o$  comes endowed with a natural blow-down. Let 
$ \{ \alpha_1, \alpha_2, \cdots \alpha_8, \beta_1, \beta_2, \cdots \beta_8 \}  \subset  \Lambda_{Z_o}  $
be the ordered set of simple roots associated to the respective blow-down.
Then, under the isomorphism $ ( \ref{asdsa}) $, 
$ \psi_{Z_o} (\alpha_i) = \psi (a_i), \ \psi_{Z_o} (\beta_i) = \psi (b_i) $ for $ 1 \leq i \leq 8 $. 
In other words, if $ \phi_2 \colon \Lambda_{Z_o} \rightarrow \Lambda $ is the marking
sending $ \{ \alpha_1, \alpha_2, \cdots \alpha_8, \beta_1, \beta_2, \cdots \beta_8 \} $ to 
$ \{a_1, \cdots a_8, b_1, \cdots b_8 \}$, then the diagram: $$
\xymatrix{
\Lambda_{Z_o} \ar [d] _{\phi_1} \ar [r] ^{\psi_{Z_o}} & {\rm Jac}(E) \ar [d] ^{\simeq} \\
\Lambda \ar [r] ^{\psi} & J^1(\m{H} ) \\
} $$
%\end{equation}           
is commutative. we conclude that the marked stable period of the marked triplet 
$(Z_o, \phi_1, \phi_2) $ coincides with the pair $ (\m{H}, \psi) $.
\par A similar procedure can be used if $ \Lambda = \Gamma_{16} $. Fix $ \{ c_1, c_2, \cdots c_{16} \} $
a basis of $ D_{16} $ simple roots in $ \Lambda $. We shall construct a set of $19$ points 
$ q_0, p_1, p_2, \cdots p_{18}$ in $E$. To begin with pick $ p_1 \in E $ such that:
$$ 3p_1 \ = \ - \left ( 2 \psi(c_1) + 2 \psi(c_2) + \cdots + 2 \psi(c_{13})+ 2 \psi(c_{14}) + 
\psi(c_{15}) + \psi(c_{16}) \right ). $$
Define then $ p_2 = p_1 $ and construct recursively 
$ p_l \ = \ p_{l-1} - \psi(c_{l-2}) $ for $ 3 \leq l \leq 17  $. Pick then $ q_0\in E$ such that:
$$ 6q_0 \ = \ 2p_1 + p_2 + p_3 + \cdots p_{17} .$$
Finally, set $ p_{18} = p_1 + 3q_0 $.  It follows then that 
$ \left (  3p_0; p_1, p_2, \cdots p_{17}; 3q_0; p_{18} \right )  $
is a special family of category (b) on $ E $. Then,  
\begin{equation}
\label{exxb}
Z_o \ = \ Z_o \left ( E; 3p_0; p_1, p_2, \cdots p_{17}; 3q_0; p_{18} \right ) 
\end{equation}
is an elliptic Type II stable $ K3$ surface with section endowed with a canonical blow-down. If $ 
\{ \gamma_1, \gamma_2, \cdots \gamma_{16}  \} \ \subset \ \Lambda_{Z_o}  $
is the ordered set of simple roots associated to the respective blow-down, then, 
under the isomorphism $ ( \ref{asdsa}) $, one has.
$ \psi_{Z_o} (\gamma_l) = \psi (c_l) $,  $ 1 \leq i \leq 16 $. Therefore, if one sets a marking 
$ \phi_2 \colon \Lambda_{Z_o} \rightarrow \Lambda $ such that the ordered basis 
$ \{ \gamma_1, \gamma_2, \cdots \gamma_{16}  \} $ is sent to $ \{ c_1, c_2, \cdots c_{16} \} $, then 
the marked stable period associated to $(Z_o, \phi_1, \phi_2) $ is $ (\m{H}, \psi) $.         
\end{proof} 

\subsection{Stable Surfaces as K3 Degenerations}
\label{sdegen} 
We have seen that the two Type II Mumford boundary divisors 
$ \m{D}_{E_8 \oplus E_8 } $ involved in the partial compactification
of $ \Gamma \backslash \Omega $ can be regarded as moduli spaces of 
periods for elliptic Type II stable $K3$ surfaces with section in the 
$ E_8 \oplus E_8 $ and $ \Gamma_{16} $ category, respectively. In this
section we justify the presence of such surfaces from a geometrical 
point of view, as they appear naturally as central fibers for certain
degenerations of $K3$ surfaces. 
\begin{dfn}
\label{defdegen} 
\par A one-variable degeneration of elliptically fibered $K3$
surfaces with section consists of a commutative diagram of analytic maps:
\begin{equation}
\label{degen1}
\xymatrix{
Z \ar [r] ^{\pi} \ar [d] & \Delta \\ 
S \ar [ur] &  \\
}
\end{equation}
where $Z$ is a smooth three-fold, $S$ is a smooth surface, 
$ \Delta $ is the unit disk. In addition, the structure requires the presence of 
an analytic section $ s \colon S \rightarrow Z $ and of a line bundle
$ \m{F} $ on $ Z $, such that:
\begin{itemize}
\item For every $t \in \Delta^* $, $Z_t$ is a smooth $K3$ surface, 
$ S_t $ is a smooth rational curve and the projection $ Z_t \rightarrow S_t $ 
is an elliptic fibration with section $ s_t \colon S_t \rightarrow Z_t $.    
\item The restriction of $\m{F} $ on $Z_t $ coincides with the line bundle
associated to the elliptic fiber in $ Z_t \rightarrow S_t$.    
\end{itemize} 
\end{dfn}  
\noindent Two degenerations $ Z \stackrel{\pi}{\rightarrow} \Delta $ and 
$ Z' \stackrel{\pi'}{\rightarrow} \Delta $ 
as in $ (\ref{degen1}) $ are said to be equivalent if one has a birational map 
$ \alpha \colon Z \rightarrow  Z' $ entering the commutative diagram:
\begin{equation}
\xymatrix{
Z \ar [r] _{\pi} \ar [d]  \ar@/^/ [rr] ^{\alpha} & \Delta & Z' \ar [l] ^{\pi'} \ar [d] \\
S \ar@/^/ [u] ^{s}\ar [ur] \ar@/^/ _{\overline{\alpha}} [rr]  & & S' \ar [ul] \ar@/_/ [u] _{s'}   
}
\end{equation}
and satisfying $  \m{F} = \alpha ^* \m{F}' $. 
\par One can think of an equivalence class of $K3$ degenerations as a 
punctured disc embedded in the moduli space $ \m{M}_{K3} $. 
Intuitively, one can then regard the degenerated central fibers $ Z_o $ as 
geometrical representatives for boundary points in a compactification of $ \m{M}_{K3} $. 
A major difficulty appears here due to the fact that equivalent degenerations can
have quite different central fibers. One tries to surmount this obstacle by restricting 
to more distinguished degenerations in the hope of obtaining a canonical model of central
fiber for each degenerating equivalence class, which is a requirement for any attempt 
of geometrical partial compactification. Along this reasoning line (see \cite{frthesis} 
\cite{pinkham} \cite{friedman} \cite{kul} for details), we restrict ourself to
degenerations $ Z \stackrel{\pi}{\rightarrow} \Delta $ as in $ (\ref{degen1}) $ which
are semi-stable (meaning that the central fiber $ Z_o $ is a surface 
with normal crossings) and satisfy $ K_Z = \m{O}_Z $. We shall call these Kulikov degenerations.  
\par Since $ \pi_1(\Delta^*) \simeq \Zee $, it is not necessarily possible to attach a 
consistent set of markings to the surfaces in such a Kulikov family. Attached to
each degeneration, there is a monodromy operator:
\begin{equation}
{\rm T } \in {\rm Aut} \left ( H^2(Z_t, \Zee) \right ) 
\end{equation} 
which can be described explicitly as the Picard-Lefschetz transformation obtained 
by transporting cycles around origin $ t=0$ in $ \Delta $ while preserving the
classes representing the elliptic structure and section. 
The operator $ T $ is unipotent, meaning $ \left ( T - I \right ) ^3 = 0 $ 
which is equivalent to saying that its
logarithm
\begin{equation}
{\rm N} \ = \ \left (T-I \right ) - \frac{1}{2} \left ( T-I \right )^2  
\end{equation}       
is a nilpotent endomorphism of $ H^2(Z_t, \mathbb{Q}) $ satisfying $ N^3 = 0 $. Complexifying
the picture, one obtains a monodromy weight filtration:
\begin{equation}
\label{lll1}
\{0\} \ \subset \ {\rm Im} \left ( N^2 \right ) \ \subset 
{\rm Im} \left ( N \right ) \cap {\rm Ker} \left ( N \right ) \ \subset \  
{\rm Im} \left ( N \right ) +  {\rm Ker} \left ( N \right ) \
\subset \ {\rm Ker} \left (N^2  \right ) \ \subset \ H^2(Z_t, \Cee) .  
\end{equation}
Moreover, as explained in \cite{schmid}, the degeneration data produces
a mixed Hodge structure on $ H^2(Z_t, \Cee) $ with weight filtration $ (\ref{lll1}) $, 
the limiting mixed Hodge structure. With respect to this structure, the nilpotent 
endomorphism $N$ becomes a morphism of mixed Hodge structures of type $(-1,-1) $.
\par Kulikov degenerations fall into three categories, denoted Type I, Type II 
and Type III, depending on whether $N=0 $, $ N^2 =0 $ but $N \neq 0 $, or $ N^3 = 0 $ 
but $ N^2 \neq 0 $. We shall restrict our attention here only to Type II families 
($ N \neq 0 $ but $ N^2 =0 $) as that will turn to be the case relevant to our prior
discussion. In this case $N$ is always an integral endomorphism. The elliptic 
Type II stable $K3$ surfaces with section appear then naturally as central
fibers for Type II degenerations with primitive $N$.   
\begin{prop} \
\begin{enumerate}
\item Let $ Z \stackrel{\pi}{\rightarrow} \Delta $ be a degeneration as in $ (\ref{degen1}) $ 
which is Kulikov of Type II with primitive endomorphism $N$. The central fiber $Z_o$ is
then an elliptic Type II stable $K3$ surface with section. 
\item For every elliptic Type II stable $K3$ surface with section $Z_o$, there exists
a Type II Kulikov degeneration $ Z \stackrel{\pi}{\rightarrow} \Delta $ 
as in $ (\ref{degen1}) $ with central fiber 
$Z_o$.     
\end{enumerate}
\label{pppp2121}
\end{prop} 
\begin{proof}
Both statements can be deduced easily from standard results on $K3$ degenerations 
(see \cite{frthesis} \cite{friedman1} \cite{friedman} and \cite{kul}). Indeed, to prove
the first part of the proposition, assume that $Z_o$ is a central fiber of a
degeneration:  
\begin{equation}
\label{degen1234}
\xymatrix{
Z \ar [r] ^{\pi} \ar [d] & \Delta \\ 
S \ar [ur] \ar@/^/ [u] ^{s} &  \\
}
\end{equation}
as in $ (\ref{degen1}) $, which is Kulikov of Type II with $N$ primitive. Recall the following 
fact from \cite{kul}:
\begin{teo} 
\label{kk}
The central fiber of a Type II Kulikov degeneration of $K3$ surfaces is always a chain of 
smooth rational surfaces:
$$ Z_0 \ = \ X_1 \cup X_2 \cup \cdots \cup X_{r+1}. $$ 
The surfaces $ V_2$, $V_3$, $ \cdots $ $ V_{r} $ are smooth elliptic ruled. 
The chain contains only double curves and all double curves are smooth and elliptic. 
\end{teo} 
\noindent The integer $r$ can be detected by from the arithmetic of the degeneration by 
writing $ N=rN_o $ with $ N_o $ primitive. Since we expect that $N$ itself is primitive, 
it has to be that $r=1$ and therefore the central fiber $Z_o$ is a union $X_1 \cup X_2 $ 
of two rational surfaces glued along an elliptic curve $D$. Let us analyze then the 
central fiber configuration: 
\begin{equation}
\label{uuu}
X_1 \cup X_2 \rightarrow S_o. 
\end{equation}
Since $S \rightarrow \Delta $ is a degeneration of smooth rational curves, $S_o$ has to 
be a chain of rational curves. The map $ (\ref{uuu}) $ is proper and its domain is a 
union of two irreducible varieties. Therefore, $S_o $ cannot have more than two 
irreducible components. We divide then our discussion into two cases:
\begin{enumerate}
\item $S_o $ is a union $S_1 \cup S_2 $ of two copies of 
$ \mathbb{P}^1 $ meeting at one point.
\item $S_o$ is a smooth rational curve.
\end{enumerate}
In the first case, $S_i$ represents the image of $X_i$ through $ (\ref{uuu}) $ and $D$ is
the fiber above the common point. We have therefore two elliptic fibrations  
$ X_i \rightarrow S_i $ agreeing over the double curve. The section $s_o \colon S_o 
\rightarrow Z_o$ allows us to regard $S_1$ and $S_2 $ as two smooth rational curves 
embedded in $X_1 $ and $X_2 $, respectively. The two curves $S_1$ and $S_2 $ meet $D$
at the same point. This is exactly the configuration required for $Z_o $ to be an 
elliptic Type II stable $K3$ surface with section of category (a). 
\par In the second case, the section $ S_o $ lies entirely inside one of two
surfaces $V_i $. Assume $ S_o \subset V_1 $. Since $S_o $ corresponds to a Cartier divisor
on $ Z_o $, it cannot intersect the double curve $D$. Therefore $ S_o^2=-2$. The projection        
$ (\ref{uuu}) $ restricts to rulings:
$$ V_i \rightarrow \mathbb{P}^1 $$
with the double curve $D$ playing the role of a bi-section on each side. We obtain therefore
that $Z_o$ is an elliptic Type II stable $K3$ surface with section of category (b).  
\par In order to prove the second assumption, we recall some facts pertaining to the 
deformation theory for stable $K3$ surfaces.
\begin{teo} (\cite{friedman1})
Let $Z_o$ be a Type II stable $K3$ surface.
\begin{enumerate}
\item $Z_o$ is smoothable and appears as central fiber in a Kulikov semi-stable degeneration.
\item The space of first-order deformations of $Z_o$:   
$$ {\bf T}^1_{Z_o} = {\rm Ext}^1 \left ( \Omega^1_{Z_o}, \m{O}_{Z_o} \right ) $$
is $21$-dimensional.
\item The versal deformation space of $Z_o$ looks like $V_1 \cup V_2 \subset {\bf T}^1_{Z_o}$ 
where $V_1$ and $V_2 $ are two smooth divisors meeting normally. The points of $V_1 $ 
corresponds to locally trivial deformations of $Z_o$. The points of $V_2 \backslash  V_1 $
represent deformations of $K3$ surfaces. $V_1 \cup V_2 $ corresponds to locally trivial
deformations of $Z_o$ which remain d-semi-stable.     
\end{enumerate} 
\end{teo}
\noindent Therefore, given an elliptic Type II stable $K3$ surface with section $Z_o$, there 
are always plenty of smoothings of $Z_o$. We just have to show that, on some of these
deformations, the two Cartier divisors $ \m{S}_o $ and $ \m{F}_o $ can be extended on 
the three-fold. The obstruction to extending a Cartier divisor of the central fiber is 
measured by the Yoneda pairing \cite{siu}:
\begin{equation}
\label{yoneda} 
< \cdot, \cdot > \colon \ {\rm Ext}^1 \left ( \Omega^1_{Z_o}, \m{O}_{Z_o} \right ) \otimes 
H^1 \left ( Z_o, \Omega^1_{Z_o} \right ) \ \rightarrow \  
{\rm Ext}^2 \left ( \m{O}_{Z_o}, \m{O}_{Z_o} \right ) = H^2 \left ( Z_o, \m{O}_{Z_o} \right ) 
\end{equation} 
which is non-degenerate for stable $K3$ surfaces. The Zariski tangent space to the smoothing 
component $V_2$ is given by the hyper-plane (see \cite{friedman}):
$$ \left \{ \sigma \in {\bf T}^1_{Z_o} \ \vert \ < \sigma , [\xi_o] > = 0 \ \right \} \ 
\subset \ {\bf T}^1_{Z_o} $$ 
where $ [\xi_o] $ is the class in $ H^1 ( Z_o, \Omega^1_{Z_o} ) $ associated to the 
Cartier divisor $ \xi_o $. The formal Zariski tangent space to the space of smoothings 
extending the elliptic structure and section is then given by:
\begin{equation}
\label{zariski}
\left \{ \sigma \in {\bf T}^1_{Z_o} \ \vert \ 
< \sigma , [\xi_o] > = < \sigma , [\m{F}_o] > = < \sigma , [\m{S}_o] >=0  \ \right \}. 
\end{equation}     
Since Yoneda pairing is non-degenerate and $ [\xi_o] $, $ [\m{F}_o] $ and $ [\m{S}_o] $ 
are independent in $ H^1 ( Z_o, \Omega^1_{Z_o} ) $, $ (\ref{zariski}) $ is $18$-dimensional. 
The space of versal deformations extending the elliptic structure and section has then a unique 
smoothing component $V'_2 $ of dimension $18$. The points $V'_2$ away from the discriminant 
locus correspond to deformations as in Definition $ \ref{defdegen} $.
\end{proof}
\noindent Let us then present the stable period map $(\ref{per22})$ as a natural extension 
of the $K3$ period correspondence $ \m{M}_{K3} \simeq  \Gamma \backslash \Omega. $ 
Assume that $ Z \rightarrow \Delta $ is a degeneration of elliptic $K3$ surfaces with section, as
in Definition $ \ref{defdegen} $, which is Kulikov, Type II semi-stable, and has 
primitive endomorphism $N$. There is then a corresponding Griffiths' period map (see \cite{griff}):
\begin{equation}
\label{holopath}
\Phi \colon \Delta^* \ \rightarrow \ \Gamma \backslash \Omega. 
\end{equation}
Following results of Mumford \cite{ash} and Schmid \cite{schmid}, one sees that $ (\ref{holopath}) $ 
extends to a holomorphic map:
$$  \w{\Phi} \colon \Delta \ \rightarrow \ \overline{\Gamma \backslash \Omega}. $$
Recall the construction of the boundary point $ \w{\Phi}(0) $. Choose a 
compatible marking $ H^2(Z_t, \Zee) \simeq \mathbb{L} $ as in section $\ref{per1}$. 
The endomorphism $N$ is integral and vanishes on both cohomology classes 
$ [\m{F}_t], [\m{S}_t] \in H^2(Z_t, \Zee)  $ corresponding to the elliptic fiber and section in 
$ Z_t \rightarrow S_t $. Therefore, it defines an isotropic rank-two sub-lattice 
$ V \subset \mathbb{L}_o $. Define $ \Lambda = V^{\perp} / V $ and let $ F $ be the Baily-Borel
component associated to $V$. The monodromy weight filtration $ (\ref{lll1}) $ associated to the 
degeneration $ Z \rightarrow \Delta $ is just the complexification of:
\begin{equation}
\label{llll}
\{0\} \ \subset \ {\rm Im} \left ( N \right ) \ \subset \  
{\rm Ker} \left ( N \right ) \ \subset \ H^2(Z_t, \Zee) .  
\end{equation}   
Taking the orthogonal part to the fiber and section classes $ [\m{F}_t] $ and $ [\m{S}_t] $, 
reduces $ (\ref{llll}) $ to: 
\begin{equation}
\label{llll2}
\begin{array}{ccccccc}
\{0\} & \subset &  {\rm Im} \left ( N \right ) & \subset &  
{\rm Ker} \left ( N \right ) \cap \{ [\m{F}_t], [\m{S}_t] \}^{\perp} &  \subset &  H^2(Z_t, \Zee) 
\cap \{ [\m{F}_t], [\m{S}_t] \}^{\perp} \\
\end{array}
\end{equation}
which corresponds under the marking to:
\begin{equation}
\label{llll3}
\{0\} \ \subset \ V \ \subset \  
V^{\perp} \ \subset \ \mathbb{L}_o .  
\end{equation}  
By the classical construction of Schmid \cite{schmid}, the family $ Z \rightarrow \Delta $ induces 
a nilpotent orbit of limiting mixed Hodge structures with weight filtration $ (\ref{llll}) $. These
structures descend, under the marking, to give a nilpotent orbit of polarized 
mixed Hodge structures on 
$ (\ref{llll3}) $. The resulting $ U(N)_{\Cee}$-orbit consists essentially of the 
decreasing filtrations:
\begin{equation}
\mathbb{L}_o \otimes \Cee \ \supset \ \{ {\rm exp}(zN) \cdot \omega_t \}^{\perp} \ \supset  \ 
\{ {\rm exp}(zN) \cdot\omega_t \} \ \supset \{0\} , \ \ z \in \Cee
\end{equation}  
where $ \{ \omega_t \} \subset \mathbb{L}_o \otimes \Cee $ is the marked period line of $Z_t $. 
But, as explained in section $ \ref{p} $, such a nilpotent orbit of Hodge structures is 
equivalent to a point on the Type II Mumford component $ \m{B}^+(F) $. $ \w{\Phi}(0) $ is the class
of this point on the quotient boundary divisor 
$ \m{D}_{\Lambda} \subset \overline{\Gamma \backslash \Omega} - \Gamma \backslash \Omega  $.          
\par Let then $Z_o$ be the central fiber of $ Z \rightarrow \Delta $. According to Proposition 
$\ref{pppp2121} $, $Z_o$ is an elliptic Type II stable $K3$ surface with section. The Clemens-Schmid 
exact sequence \cite{gs}:
\begin{equation}
\label{ccss1}
\{0\} \rightarrow H^0(Z_t) \stackrel{(-2,-2)}{\longrightarrow} H_4(Z_o) 
\stackrel{(3,3)}{\longrightarrow} H^2(Z_o) \stackrel{(0,0)}{\longrightarrow}
H^2(Z_t) \stackrel{N}{\longrightarrow} H^2(Z_t) \stackrel{(-2,-2)}{\longrightarrow}
H_2(Z_o) \stackrel{(3,3)}{\longrightarrow} H^4(Z_t) \cdots \ \ 
\end{equation}
allows one to relate the geometric mixed Hodge structure of $Z_o$ with 
the limiting mixed Hodge structure associated to the degeneration 
$ Z \rightarrow \Delta $. A careful analysis of $ (\ref{ccss1}) $ reveals that:  
\begin{teo}
The boundary point $ \w{\Phi}(0) \in \m{D}_{\Lambda}$ is the stable period of $Z_o$, 
as defined in section $ \ref{sstableper} $. 
\end{teo}
\section{Boundary Components and Flat Bundles}
\label{lapp}
There exists a second geometric interpretation, more relevant from the point of view of 
heterotic/F-theory duality, for the boundary points on the two Type II 
Mumford divisors $ \m{D}_{\Lambda} $ with $ \Lambda = E_8 \oplus E_8 $ or $ \Lambda = \Gamma_{16} $. 
Recall that, given a Baily-Borel component $F$, $ \m{D}_{\Lambda} = \Gamma^+_F \backslash \m{B}^+(F) $ 
and the points in $ \m{B}^+(F) $ are in one-to-one correspondence to pairs $ (\m{H}, \psi ) $ of
polarized weight-one mixed Hodge structures $ \m{H} $ on $V$ together with abelian group
homomorphisms $ \psi \colon \Lambda = V^{\perp}/V \rightarrow J^1(\m{H} ). $ Such a pair is known
to determine a flat $G$-connection over the elliptic curve $E = J^1(\m{H}) $. The Lie group
$G$ is $(E_8 \times E_8) \rtimes \Zee_2 $ if 
$ \Lambda =  E_8 \oplus E_8$ and $ {\rm Spin}(32)/\Zee_2 $ 
if $ \Lambda = \Gamma_{16} $.
\par Let us briefly review the connection. For explicit details, see \cite{mf1}. It is a 
standard fact that, given a compact Lie group $G$ and a smooth two-torus $E$, there is a 
bijective correspondence between 
the equivalence classes of flat $G$-connections on $E$ and their associated  
holonomy morphisms $ \pi_1(E) \rightarrow G $, up to conjugation. One can
therefore formally identify a flat connection with a commuting pair of elements in $G$, 
up to simultaneous conjugation. Fix a maximal torus $T \subset G $. If $G$ is simply connected, 
(in particular for $G = (E_8 \times E_8) \rtimes \Zee_2 $), it was shown that any given pair 
of commuting elements in $G$ can be simultaneously conjugated in $T$. The same statement is true for
$ G =  {\rm Spin}(32)/\Zee_2 $, providing that one considers only connections which can
be lifted to $ {\rm Spin}(32)$-connections. In this way, a flat $G$-connection on $E$ can be formally 
understood as an element in:
$$ {\rm Hom}\left ( \pi_1(E) , T \right ) / W  $$
where $W$ is the Weyl group of $G$. The lattice $ \Lambda $ plays the role of the lattice of 
the maximal torus $T$. In this framework:
$$ {\rm Hom}\left ( \pi_1(E) , T \right ) \ \simeq \ 
{\rm Hom} \left ( \pi_1(E) , U(1) \otimes \Lambda \right ) \ \simeq \ 
{\rm Hom} \left ( \pi_1(E) , U(1) \right ) \otimes \Lambda . $$
The first factor of the last term above represents the set of gauge equivalence 
classes of flat hermitian line bundles over $E$. In the presence of a complex structure
on $E$, one can identify $ {\rm Hom} \left ( \pi_1(E) , U(1) \right ) $ to $ {\rm Pic}^o(E) $ which, 
in turn, is a complex torus isomorphic to $E$. There exists then a bijective correspondence between 
flat $G$-connections and points of the analytic quotient:
\begin{equation}
\label{l60} 
E \otimes \Lambda / W
\end{equation}    
which one can see as the moduli space of flat $G$-bundles over $E$. 
\par Due to the unimodularity $ \Lambda ^* \simeq \Lambda $, each element in
$ (\ref{l60}) $ can be regarded as a class of a morphism $ \Lambda \rightarrow E $.
Along this idea, one can associate to any Type II boundary point of $ \m{B}^+(F) $, a smooth 
elliptic curve $E=J^1(\m{H})$ and a flat $G$-bundle. In section $ \ref{ellipticflat} $ we shall
show that all possible flat $G$-bundles are realized\footnote{For $ G =  {\rm Spin}(32)/\Zee_2 $, we only look 
at flat bundles liftable to $ {\rm Spin}(32)$.} and that two points in $ \m{B}^+(F) $
determine equivalent pairs of elliptic curves and flat bundles exactly when they belong to the
same $ \Gamma^+_F $-orbit. This will lead to a holomorphic identification between the boundary
divisors $ \m{D}_{E_8 \oplus E_8} $, $ \m{D}_{\Gamma_{16}} $ and the moduli spaces 
$ \m{M}_{E, E_8 \times E_8 \rtimes \Zee_2} $ and  $ \m{M}_{E, {\rm Spin}(32)/\Zee_2} $ of
equivalence classes of 
pairs of elliptic curves and flat bundles, respectively.                
\section{Explicit Description of the Parabolic Cover}
\label{main} 
In this section we give an explicit description of the two Type II boundary divisors
$ \m{D}_{\Lambda} = \Gamma_F \backslash \m{B}(F)$ and identify precisely the holomorphic 
type of the parabolic fibrations given in $(\ref{parabcover})$:
\begin{equation}
\label{pparacover}
\Theta \colon \Gamma_F \backslash \overline{\Omega(F)/U(N)_{\Zee}} \ \rightarrow \ \Gamma_F \backslash \m{B}(F). 
\end{equation}
This leads, following $ \ref{arith} $, to a description of the structure of $ \m{M}_{K3} $ 
in a neighborhood of $ \m{D}_{\Lambda} $. 
\subsection{Fixing the Parabolic Group} 
Let $F$ be a fixed Type II Baily-Borel boundary component for $ \Gamma \backslash \Omega $. 
Denote by $  V $ the associated primitive isotropic rank-two sub-lattice of $ \mathbb{L}_o $ 
and set, as in $ \ref{arith} $:  
$$ \Lambda = - \left ( V^{\perp}/V \right ) $$
$$ {\rm P(F)} \ = \ {\rm Stab} \left ( V_{\mathbb{R}} \right ) \subset O^{++}(2,18) $$
$$ {\rm W}(F) \ = \ {\rm the \ unipotent \ radical \ of \ P(F)} $$   
$$ {\rm U}(F) \ = \ {\rm the \ center \ of \ W(F) }. $$
It follows then that ${\rm U}(F) $ is a one-dimensional Lie group with Lie algebra:
$$ {\rm u}(F) \ = \ \{ {\rm N} \in {\rm End}_{\mathbb{R}} \left ( \mathbb{L}_o \otimes 
{\mathbb{R}} \right ) 
 \ \vert \  {\rm Im} (N) = V_{\mathbb{R}} \ {\rm and} \ (Nx,y) + (x,Ny) = 0 \}. $$
\begin{lem}
\label{l007}
There exists a basis $ \{ A, B \} $ for $V $ such that the endomorphism 
$N \colon \mathbb{L}_o \rightarrow \mathbb{L}_o$ defined by:
\begin{equation}
\label{morphism} 
N(x) = (x,B) A - (x,A) B 
\end{equation}  
is primitive, integral and belongs to $ {\rm u}(F) $. 
\end{lem}  
\begin{proof}
Let $A$ be a primitive element of $V$. Due to unimodularity, there exists 
$ A' \in \mathbb{L}_o $ satisfying $ (A, A') = 1 $. Pick $ B \in V $, primitive, 
such that $ (B, A') = 0 $. It follows that 
$ \{ A, B \} $ forms a basis for $ V $.
\par Let then $N$ be the endomorphism defined in $ (\ref{morphism}) $. Pick 
$ C \in \mathbb{L}_o $ such that $ (B,C) = 1 $ and define:
$$ B' \ = \ C - (C,A') A - (C,A)A' + (C,A) (A,A') A \ \in V.  $$
One has verifies that $ N(A') =-B $ and $ N(B') = A $. Therefore,  $N$ is primitive 
and $ {\rm Im}(N) = V $. Moreover, since:
$$ (Nx,y) \ = \ (x, B)(y,A) - (x,A)(y,B) \ = \ - (x, Ny) $$
for any $x,y \in \mathbb{L}_o $, the endomorphism $N$ belongs to $ {\rm u}(F) $.             
\end{proof} 
\noindent In order to facilitate future computations, we shall introduce a 
special coordinate system on $\mathbb{L}_o$. The linearly independent family 
$ \{ A', B', A, B \} $, can be seen to provide a decomposition: 
\begin{equation}
\label{lat}
\begin{array}{ccccc}
\mathbb{L}_o  \ \simeq \ & \underbrace{\left ( \Zee \cdot A' \oplus \Zee \cdot B'\right ) } 
\ \oplus & \underbrace{\left ( \Zee \cdot A \oplus \Zee \cdot B\right ) } \ \oplus & 
 \Lambda . \\
&      \Zee^2                                &     \Zee^2                                & & \\
\end{array}
\end{equation} 
In this light, any element $ \mathbb{L}_{o} $ 
(or $ \mathbb{L}_{o} \otimes \Cee$) can be written uniquely as:  
$$ x_1 A' + x_2 B' + y_1 A + y_2 B + z.    $$ 
We convene therefore to regard the elements of $ \mathbb{L}_{o} $ 
as a triplets $ (x,y,z) $ with 
$ x=(x_1,x_2) \in \Zee^2$, $ y = (y_1, y_2) \in \Zee^2  $ 
and $ z \in \Lambda $. The quadratic pairing on $ \mathbb{L}_o $ is 
recovered as:    
$$
\left ( (x, y, z), \ (x', y', z') \right ) = x.y'+x'.y - (z,z') $$
where the first two dot-pairings on the left represent the standard 
Euclidean pairing on $\Zee^2$ and $(\cdot, \cdot) $ is the pairing of $\Lambda$. 
Under this rule, the isotropic lattice $ V $ corresponds to the space of 
triplets $ (0,y,0) $ and the integral endomorphism $N$ is 
given by $ N(x, y, z) \ = \ (0, Tx, 0) $ 
with $ T \colon \Zee^2 \rightarrow \Zee^2 $ is the standard skew-adjoint endomorphism  
$ T(x_1,x_2) = (x_2, -x_1) $. 
\par As in $ \ref{arith} $, we define the groups:
$$ {\rm U}(N)_{\Cee} \ \colon = 
\ \{ {\rm exp} \left ( \lambda N \right ) \ \vert \ \lambda \in \Cee \ \}  $$
$$ {\rm U}(N)_{\Zee} \ \colon = 
\ \{ {\rm exp} \left ( \lambda N \right ) \ \vert \ \lambda \in \Zee \ \}  $$
leading to the sequence of inclusions:
$$ \Omega \ \subset \ \Omega(F) = {\rm U}(N)_{\Cee} \cdot \Omega \ \subset \ \Omega^{\vee} . $$ 
We shall use the newly introduced coordinate system to analyze these inclusions. Let 
$r \colon \mathbb{L}_o \otimes \Cee \rightarrow \mathbb{R} $ be the function defined
by $ r(\omega) = - i (N\omega, \bar{\omega}) $. This function is invariant under the
action of $ U(N)_{\Cee} $. In fact, if $ \omega = (x,y,z) $ then $ r(\omega) = 
2 {\rm Im}(x_1 \overline{x_2}).$ Let 
$$ \Omega^{\vee} = \Omega^-(F) \ \cup \ \Omega^0(F)  \ \cup \ \Omega^+(F) $$ 
be the decomposition of  $\Omega^{\vee}$ in subsets for which $ r(\omega) $ is 
strictly negative, zero and strictly positive, respectively. 
\begin{prop} The following statements hold:   
\begin{enumerate}
\item $ \Omega(F) = \Omega^{-}(F) \ \cup \ \Omega^{+}(F) $.  
\item If $[\omega] = [a,b,c] \in \Omega^+(F) $, then:
\begin{equation}
\label{fractie}
- \frac{a_1}{a_2} 
\end{equation}
is a well-defined element of the upper-half plane.    
\end{enumerate}
\label{lem2}
\end{prop}
\begin{proof}
Let $ [\omega] \in \Omega(F) $. We show that $ r(\omega) \neq 0 $ by proving that 
the opposite statement leads to a contradiction. Indeed, assume that 
$ r(\omega) = 0 $. Since $ \omega \in \Omega(F) $, there exists 
$ \omega_o \in \Omega $ such that $ \omega = {\rm exp}(zN). \omega_o = \omega_o + z N \omega_o $
for some $ z \in \Cee $. Then:
\begin{equation}
\label{aaa}
(\omega, \bar{\omega} ) = (\omega_o, \bar{\omega_o} ) - 2 \ {\rm Im}(z) \cdot 
r (\omega_o) = (\omega_o, \bar{\omega_o} )>0. 
\end{equation} 
But $ r(\omega) = 0 $ also implies that $ \{ \omega, N\omega , N \bar{\omega} \} $ 
span an isotropic subspace of $ \mathbb{L}_o \otimes \Cee $. Clearly, $ \omega $ and $ N\omega $ 
are independent (otherwise $ \omega \in V_{\Cee} $, contradicting $ (\ref{aaa}) $).
Since the largest isotropic subspace in  $ \mathbb{L}_o \otimes \Cee $ is two-dimensional, 
it has to happen that 
$ N \bar{\omega} $ is generated by $ \omega $ and $ N\omega $. But that also 
implies $ \omega \in V_{\Cee} $, leading to a contradiction. 
\par This shows that:
$$ \Omega(F) \subset \Omega^{-}(F) \ \cup \ \Omega^{+}(F). $$
The reverse inclusion is straightforward.
\par Turning to the second statement, one has $ r(\omega) = 2 {\rm Im}(a_1 \overline{a_2})> 0$. 
The denominator of $(\ref{fractie}) $ is therefore non-zero. Moreover, the same formula
leads to:
$$ {\rm Im} \left ( - \frac{a_1}{a_2} \right ) = \frac{r(\omega)}{2 \vert a_2 \vert ^2 }$$
which assures us that $(\ref{fractie}) $ is an element of the upper half-plane.       
\end{proof}     
\noindent We are now in position to write explicit formulas for the geometric assignment, 
described in $\ref{p} $, that associates to a nilpotent orbit in
$$ \m{B}(F) \ = \ \Omega(F) / {\rm U}(N)_{\Cee}, $$
a pair $ (\m{H}, \psi) $ consisting of a weight-one Hodge structure $ \m{H} $ on $V$ 
and a homomorphism $ \psi \colon \Lambda \rightarrow J^1(\m{H})  $.
\par Under the identification $  V \simeq \Zee^2 $, provided by the 
basis $ \{A,B \} $, the skew-symmetric bilinear form $ (\cdot, \cdot)_1 $ is transported
to $ (x,y)_1 = x . Ty = x_1y_2-x_2y_1$. The Hodge structures of weight one on $V$ which
are polarized with respect to $ (\cdot, \cdot)_1 $ are then indexed by purely imaginary
complex numbers $\tau$ belonging to the upper half-plane $ \mathbb{H} $. Every such $ \tau $
induces the polarized weight-one Hodge structure:
\begin{equation}
\label{tauhodge}
0 \ \subset \  \{ A + \tau B \}  \ \subset \ V_{\Cee} 
\end{equation}   
and the correspondence is one-to-one.  
\par Let then $ [\omega] = [a,b,c] $ be an element in $ \Omega(F) $. As described in $ \ref{p} $, 
the Hodge structure $ \m{H} $ associated to the nilpotent orbit of $ [\omega] $ in 
$$ \m{B}(F) \ = \ \Omega(F) / \ {\rm U}(N)_{\Cee} $$
is given by the filtration:
\begin{equation}
\label{uhodge}
0 \ \subset \ \{ [\omega] \} ^{\perp} \ \cap \ V_{\Cee}   \ \subset \  
V_{\Cee}.
\end{equation}
Using the coordinate framework, the middle space in $ (\ref{uhodge}) $ is:
$$ \{ [\omega] \} ^{\perp} \ \cap \ V_{\Cee} \ = \ \{ (0,y,0) \in V_{\Cee} \ \vert \ a.y = 0 \ \}. $$
An identification of the two filtrations $ (\ref{tauhodge}) $ and $ (\ref{uhodge}) $ leads 
one to:
\begin{equation}
\label{tauform}
\tau \ = \ - \frac{a_1}{a_2}.  
\end{equation}
Connecting $ (\ref{tauform}) $ to Proposition $ \ref{lem2} $ we see that the decomposition
$$ \m{B}(F) = \Omega(F)/ U(N)_{\Cee} = 
\Omega^+(F)/ U(N)_{\Cee} \ \cup \ \Omega^-(F)/ U(N)_{\Cee}$$
corresponds to the decomposition 
$ \m{B}(F) = \m{B}^+(F) \cup \overline{\m{B}(F)} $ of section $ \ref{p} $. 
The Hodge structure $ \m{H} $ is polarized with respect to $ (\cdot, \cdot)_1 $ if and only if
$ [\omega] \in \Omega^+(F) $. 
\par The coordinate framework can also be used to give a straightforward procedure 
constructing the extension homomorphism: 
\begin{equation}
\psi \colon \Lambda \rightarrow J^1(\m{H}) = 
V_{\Cee} / \left ( \{ \omega \} ^{\perp} \cap \ V_{\Cee} \right ) + V  
\end{equation} 
associated to $ [\omega] = [a,b,c] $. If $ \gamma \in \Lambda $, 
choose a lifting $ \w{\gamma} = ( 0, \beta, \gamma ) \in V^{\perp}_{\Cee} $ 
such that $ \left ( \w{\gamma}, \omega \right ) = 0 $. This amounts to choosing 
$ \beta \in V_{\Cee} $ with $ \beta. a = \gamma.c $. Clearly, such
a $ \beta $ is not unique but two different choices always differ by an element 
in $ \{ [\omega] \} ^{\perp} \cap \ V_{\Cee} $. Moreover, if one denotes:
\begin{equation}
\label{r1}
z \ = \ \frac{c}{a_2} \ = \ z_2 - \tau z_1, \ \ z_1, z_2 \in \Lambda_{\mathbb{R}} 
\end{equation}      
then the homomorphism $ \psi $ can be described as assigning:
$$ \gamma \ \longmapsto \  \left ( ( \gamma, z_1), \ 
(\gamma, z_2) \right ) \in J^1(\m{H}) .$$ 
The element $z \in \Lambda_{\Cee} $, defined as in $(\ref{r1}) $, totally 
controls the homomorphism $ \psi $.
\par We have reached therefore the following conclusion:
\begin{teo}
\label{tteeoo1}
The geometric correspondence of section $ \ref{p} $ which associates to boundary points 
in $ \m{B}^+(F) $ pairs
$ (\m{H}, \psi) $ of polarized weight-one Hodge structures on $V_{\Cee} $
and extension homomorphisms $ \psi \colon \Lambda \rightarrow J^1(\m{H}) $ induces 
an identification:
\begin{equation}
\label{a12}
\m{B}^+(F) \ = \ \Omega^+(F) /  {\rm U}(N)_{\Cee} \ \simeq \ \mathbb{H} \times \Lambda_{\Cee}. 
\end{equation}
Under this identification, the holomorphic $\Cee$-fibration of $ (\ref{parabcover1}) $ 
is described by the the map: 
\begin{equation}
\label{th}
\w{\Theta} \colon \Omega^+(F) \rightarrow \mathbb{H} \times \Lambda_{\Cee}, \ \ 
\w{\Theta} \left ( [a,b,c] \right ) \ = \ \left ( -\frac{a_1}{a_2}, \ \frac{c}{a_2} 
\right ). 
\end{equation}   
\end{teo} 
\noindent One immediately verifies in $(\ref{th})$ the main features of $ (\ref{parabcover1})$,  
namely :
\begin{itemize} 
\item $ \w{\Theta} $ is an onto holomorphic map.
\item $ \w{\Theta} $ is invariant under the action 
of  $ {\rm U}(N)_{\Cee} $ on $ \Omega^+(F)$ and the fibers of $\w{\Theta} $ coincide with 
the orbits of the ${\rm U}(N)_{\Cee} $-action. $ \Theta $ is therefore a holomorphic 
${\rm U}(N)_{\Cee} $-principal bundle. 
\end{itemize}
%\noindent Moreover, after factoring out the action of ${\rm U}(N)_{\Zee} $ in the fibers 
%of $ \w{\Theta} $, one obtains a holomorphic $ \Cee^* $-bundle:
%\begin{equation}
%\label{paracover11}
%\Omega^+(F) / {\rm U}(N)_{\Zee} \ 
%\rightarrow \ \m{B}^+(F)  \simeq \mathbb{H} \times \Lambda_{\Cee}.
%\end{equation}
\noindent At this point, recall that one obtains the 
parabolic cover $ (\ref{pparacover}) $ by further taking the 
quotient with respect to the action of the parabolic group of integral isometries 
$ \Gamma_F^+ = P(F) \cap \Gamma^+ $. It is important 
therefore to understand the group 
$ \Gamma_F^+$ and its action on $ \Omega^+(F) / {\rm U}(N)_{\Zee} $ 
and $ \mathbb{H} \times \Lambda_{\Cee} $. 

\subsection{Description of $ \Gamma^+_F$ and its action on $\m{B}^+(F)$}  
\label{ellipticflat}
%We continue to work within the framework of model $ (\ref{lat}) $:
%$$ \mathbb{L}_o \ = \ \Zee^2 \oplus \Zee^2 \oplus \left ( - \Lambda \right ) $$
%\begin{equation}
%\label{inpr} 
% (x,y,z).(x',y',z') \ = \ x.y' + x'.y - z.z' 
%\end{equation}
%with
%$$ N \colon \mathbb{L}_o \rightarrow \mathbb{L}_o,  \ \ N(x,y,z) = (0,Tx,0). $$
%Then:
The integral isometries $ \Gamma^+_F $ can be given a matrix description using the
coordinate framework $ (\ref{lat}) $.
\begin{lem}
A transformation in $ \Gamma $, stabilizing the isotropic sub-lattice 
$V$ is of the form: 
\begin{equation}
\label{type}
g( m, Q, R, F ) \ = \ \left ( \ 
\begin{array}{ccc}
m & 0 & 0 \\
R & \w{m} & Q f \\
Q^t m & 0 & f \\
\end{array}
 \ \right ) 
\end{equation}\
where:
\begin{enumerate}
\item $ m \in {\rm GL}_2(\Zee) $.
\item $ \w{m} = \left ( m^t \right )^{-1} . $ 
\item $ Q \in {\rm Hom} (\Lambda, \Zee^2 ) $, $ R \in {\rm End} (\Zee^2) $ 
satisfying $ R^t m + m^t R \ = \ m^t Q Q^t m. $ 
\item $ f $ is an isometry of $ \Lambda $. 
\end{enumerate}
Moreover $ g( m, Q, R, f ) \in \Gamma_F^+ $ if and only if $ m \in {\rm SL}2(\Zee) $. 
\end{lem}    
%\begin{proof}
%Let $ \varphi \in  \Gamma $ represented by a transformation matrix:
%$$
% \left ( \ 
%\begin{array}{ccc}
%m & n & p \\
%R & h  & S \\
%D & E & f \\
%\end{array}
% \ \right ). 
%$$ 
%Assume that $ \varphi $ preserves $U$. This implies immediately that $ n = 0 $ and $ E = 0 $. 
%Moreover, $ \varphi $ has to preserve the bilinear form $ (\ref{inpr}) $. One obtains then:
%$$ h^t \cdot m = I $$
%$$ h^t \cdot p=0 $$ 
%$$ f^t \cdot f = I $$
%$$ S^t \cdot m = f^t \cdot D  $$
%$$ R^t \cdot m + m^t \cdot R = D^t \cdot D. $$
%We conclude that $ m \in {\rm GL}(\Zee^2) $, $ h = \left ( m^t \right )^{-1} $, 
%$ p = 0 $ and $ f $ is an isometry for $ \Lambda $. Furthermore, denoting 
%$ Q = S \cdot f^{-1} $, 
%we obtain $ D = Q^t \cdot m $ and 
%$ R^t \cdot m + m^t \cdot R \ = \ m^t \cdot Q \cdot Q^t \cdot m $. In other words, 
%$$ \varphi = g( m, Q, R, f ). $$ 
%Clearly, if $ \varphi \in \Gamma^+_F $ then $ {\rm det}(m) = 1 $ and 
%therefore $ m \in {\rm SL}(\Zee^2) $.     
%\end{proof}
%\vspace{.1in}
\noindent This gives a matrix characterization for $ \Gamma_F^+ $. The 
group multiplication law goes as follows:
\begin{equation}
g( m_1, Q_1, R_1, d_1 ) \cdot g( m_2, Q_2, R_2, d_2 ) \ = \
\end{equation} 
$$ \ = \ 
g( m_1 m_2, \ Q_1+ \w{m_1} Q_2 f_1^t , \ R_1 m_2 + 
\w{m_1} R_2 + Q_1 f_1 Q_2^t  m_2 , \ f_1 f_2  ). $$ 
In particular: 
$$ g( m, Q, R, f )^{-1} \ = \ g( m^{-1}, \ -m^t Q f, \ R^t, \ f^{-1} ) . $$   
\noindent We single out the following special subgroups of $ \Gamma_F^+ $ :
\begin{enumerate}
\item $ {\rm U}(N)_{\Zee} \ = \ \left \{ g(I,0,R,I) \ \vert \ R+R^t=0 \right \} \ \simeq \ \Zee $.
\item $ \m{S} \ = \ \left 
\{ g(m,0,0, I) \ \vert \ m \in {\rm SL}_2(\Zee) \right \} \ \simeq \ {\rm SL}_2(\Zee)$. 
\item $ \m{W} \ = \ \left \{ g(I,0,0,f) \ \vert \ f \in \m{O}(\Lambda) 
\right \} \ \simeq \ \m{O}(\Lambda) $.
\item $ \m{T} \ = \ \left \{ g(I,Q,R,I) \ \vert \ R+R^t = QQ^t \right \}$.
\end{enumerate}
\noindent It can be verified that:
\begin{itemize}
\item $ {\rm U}(N)_{\Zee} \subset Z \left ( \Gamma_F^+ \right ) $  
\item $ \m{T} \vartriangleleft \Gamma_F^+ .$
\item $ \m{S} \cap \m{W} = \{ \pm I \} .$ 
\item $ \Gamma_F^+ $ decomposes as a 
semi-direct product $ \m{T} \rtimes \left ( \m{W} \times_{ \{ \pm I \} }  \m{S} \right ) $.  
\end{itemize}
%
%\subsection{The $ \Gamma_F^+$ action on $\m{B}^+(F)$}
%\label{ellipticflat}
\noindent The parabolic subgroup $ \Gamma_F^+$ acts on the total space $ \Omega^+(F) $ of 
the holomorphic $ \Cee^*$-bundle:
\begin{equation}
\label{parab2233}
  \w{\Theta} \colon \Omega^+(F) \ \rightarrow \ \mathbb{H} \times \Lambda_{\Cee}, \ \ \ 
\w{\Theta} \left ( [a,b,c] \right ) \ = \ \left ( -\frac{a_1}{a_2}, \ \frac{c}{a_2} 
\right ).
\end{equation}
There is a compatible action on $ \mathbb{H} \times \Lambda_{\Cee} $ which 
carries an important geometric significance.  
\par Recall that a pair $ (\tau, z) \in \mathbb{H} \times \Lambda_{\Cee} $ determines
a polarized mixed Hodge structure on $ V_{\Cee} $ together with a homomorphism 
$ \psi \colon \Lambda \rightarrow J^1(\m{H}) $ given essentially by $  
\psi (\gamma) \ = \ \left ( (\gamma, z_1), (\gamma, z_2) \right ) $ where $ z = z_2 - \tau z_1 $.    
As mentioned earlier, the Jacobian $J^1(\m{H})$ can be regarded as 
an elliptic curve $ E_{\tau} \ = \ \Cee / \Zee \oplus \tau \Zee  $ and, in this 
setting, the morphism $ \psi $ determines a flat $ G $-connection over 
$ E_{\tau} $ (the Lie group $G$ is 
$ E_8 \times E_8 \rtimes \Zee_2 $ if $ \Lambda = E_8 \oplus E_8 $ lattice and 
$ G = {\rm Spin}(32)/\Zee_2 $ if $ \Lambda = \Gamma_{16} $). 
\par Denote by $ \pi \colon \mathbb{H} \times \Lambda_{\Cee} \rightarrow \mathbb{H} $ the projection
on the first coordinate. Taking then:
$$ 
\m{L}_{\tau} \colon = \{ \tau \} \times \Lambda \otimes \left ( \Zee \oplus \tau \Zee \right ) \subset 
\pi^{-1}(\tau) 
$$ 
one obtains a family of 32-dimensional lattices, parameterized by $\tau$, 
moving in the fibration $ \pi $. 
\begin{dfn}
Let $ \Pi $ be the group of holomorphic automorphisms of the fibration $\pi$ which preserve 
the lattice family $\m{L} $ and cover $ {\rm PSL}(2, \Zee) $ transformations on $ \mathbb{H} $. 
\end{dfn} 
\noindent It turns out that two elements $ (\tau, z) $ and $ (\tau', z') $ of 
$ \mathbb{H} \times \Lambda_{\Cee} $ determine isomorphic pairs of elliptic curves 
and flat $G$-connections if and only if they can be transformed one into another through an isomorphism 
in $ \Pi $. In this sense, the analytic space:
\begin{equation}
\m{M}_{E, G} \ = \ \Pi \ \backslash \left ( \mathbb{H} \times \Lambda_{\Cee} \right )   
\end{equation}     
can be seen as the moduli space of pairs of elliptic curves and flat 
$G$-bundles\footnote{Again, in the case $\Lambda = \Gamma_{16} $, one considers only 
${\rm Spin}(32)$-liftable connections}. 
\begin{teo}
There is a short exact sequence of groups:
\begin{equation}
\{1\} \ \rightarrow \ {\rm U}(N)_{\Zee} \ \rightarrow \ \Gamma_F^+ \ 
\stackrel{\alpha}{\rightarrow} \ \Pi \ \rightarrow \ \{1\}. 
\end{equation} 
with respect to which the analytic fibration $ \w{\Theta} $ of $ (\ref{parab2233}) $ 
is $ \alpha$-equivariant. This induces a holomorphic identification:
\begin{equation}
\label{identif}
\m{D}_{\Lambda} \ = \ \Gamma^+_F \backslash \m{B}^+(F) \ \simeq \ 
\Pi \backslash  \left ( \mathbb{H} \times \Lambda_{\Cee} \right ) \ = \ \m{M}_{E,G}
\end{equation}  
between the Type II boundary divisor $\m{D}_{\Lambda} $ corresponding to $F$ and 
the moduli space of pairs of elliptic curves and flat $G$-bundles. Moreover, 
under $ (\ref{identif}) $, the quotient map:
\begin{equation} 
\label{pcc33}
\Theta \colon \Gamma_F^+ \ \backslash \Omega^+(F) \ \rightarrow \ \Pi \backslash 
\left ( \mathbb{H} \times \Lambda_{\Cee} \right ).   
\end{equation}
is exactly the parabolic Seifert fibration $ (\ref{parabcover}) $ of section 
$ \ref{arith} $. 
\end{teo}  
\begin{proof}
For any $ \varphi \in \Gamma^+_F $, one can construct a well-defined 
automorphism of $ \mathbb{H} \times \Lambda_{\Cee} $ by taking 
$$ (\tau, z) \ \longmapsto \Theta\left ( \varphi(\omega) \right ) $$ 
where $ [\omega] $ is a lift (under $\w{\Theta}$) of 
$ (\tau,z) $ in $ \Omega^+(F) $. We claim that all such transformations
are elements of $ \Pi $. 
\par Let $ (\tau, z) \in \mathbb{H} \times \Lambda_{\Cee} $ and 
$ g( m, Q, R, f ) \in \Gamma_F^+ $ defined as in $ (\ref{type}) $. Choose
$ [\omega] = [ x, y, z] \in \w{\Theta}^{-1}(\tau,z) $. It can be assumed that 
$ x = (-\tau, 1) $ and $  z = z_2 - z_1 \tau $ with $z_1$, $z_2 \in 
\Lambda_{\mathbb{R}} $. 
\par If $ m \in {\rm SL}_2(\Zee) $ is has the matrix form:
$$  
\left ( 
\begin{array}{cc}
a & b \\
c & d \\ 
\end{array}
\right ) \ $$
then 
$$  \w{m} = 
\left ( 
\begin{array}{cc}
d & -c \\
-b & a \\ 
\end{array}
\right ) \ $$
and the action of $ g( m, Q, R, f) $ is just: 
$$ g( m, Q, R, f ) . [\omega] \ = \ \left [ m\cdot x, \  
R \cdot x + \w{m} \cdot y + Q\cdot f \cdot z ,   \ 
Q^t \cdot m \cdot x + f \cdot z 
\right ]. $$ 
An easy calculation shows that 
$ \w{\Theta} \left (  g( m, Q, R, f ).[\omega] \right ) = (\tau', z') $ 
with
$$ \tau' = \frac{a \tau - b }{-c \tau + d }, \ \ \ z' = Q^t ( - \tau', 1) + 
\left ( d f(z_2) - b f(z_1) \right ) + \left ( -c f(z_2) + a f(z_1) \right ) \tau'. $$  
It is clear then that the transformation: 
\begin{equation}
\label{nnn} 
(\tau, z) \ \rightarrow \ \w{\Theta} \left ( g( m, Q, R, f ) . [\omega] \right ) 
= ( \tau', z') 
\end{equation} 
covers a $ {\rm PSL}(2, \Zee) $ transformation on the first factor of 
$ \mathbb{H} \times \Lambda_{\Cee} $ corresponding 
to the matrix action of 
$$  \w{m} = 
\left ( 
\begin{array}{cc}
a & -b \\
-c & d \\ 
\end{array}
\right ).  $$
In addition, one notes that 
the transformation $ (\ref{nnn}) $ preserves the lattice family $\m{L}$. It is  
therefore with a well-defined transformation in $ \Pi $. 
\par The above assignment induces a group homomorphism 
$ \alpha \colon \Gamma_F^+ \rightarrow \Pi $. It can be easily seen that 
$ \alpha \left ( g( m, Q, R, f ) 
\right ) = 1 $ requires $ m = I $, $ f = I $ and $ Q = 0 $. This implies that  
$ {\rm Ker}(\alpha) = {\rm U}(N)_{\Zee} $. 
\par Let us check that $ \alpha $ is an onto morphism. For this purpose, we 
single out the following special subgroups of $ \Pi $.  
\begin{itemize}
\item $ \m{S}_{\Pi} \ = \  \left \{ \psi \in \pi \ \Big\vert \ 
\psi(\tau, z \otimes \lambda) \ = \ \left ( 
\frac{a \tau + b}{c \tau + d} , \ \frac{z}{c\tau+d} \otimes \lambda 
\right ) \ \right \} $
\item $ \m{T}_{\Pi}  \ = \ 
\left \{ \psi \in \Pi  \ \Big\vert \ 
\psi ( \tau, z \otimes \lambda) = ( \tau, z \otimes \lambda + 1 \otimes q_1 + 
\tau \otimes q_2 ) \ \ {\rm where} \ \ (q_1, q_2) \in 
\Lambda \oplus \Lambda   
\right \} $
\item $ \m{W}_{\Pi} = \ \left \{ \psi = {\rm id} \oplus 
( {\rm id} \otimes f ) \in \Pi
\ \Big\vert \ f \in O(\Lambda)  \right \} $.  
\end{itemize}          
The three subgroups $ \m{S}_{\Pi}$, $\m{T}_{\Pi}$ and $\m{W}_{\Pi} $ generate 
the entire $ \Pi $. In addition, note that:  
$$ \m{S}_{\Pi}  \cap \m{W}_{\Pi}  =  
\left \{ \psi \in \Pi  \ \Big\vert \ 
\psi(\tau, z \otimes \lambda) \ = \ \left ( 
\tau, \ \pm z \otimes \lambda 
\right ) \ \right \} = \{ \pm I \} $$
and if $ p \colon \Pi \rightarrow {\rm PSL}(2, \Zee) $ is the
projection to ${\rm PSL}(2, \Zee) $ then $ {\rm Ker}(p) $ is 
generated by $ \m{W}_{\Pi} $ and $ \m{T}_{\Pi}  $. 
One concludes from these facts that $ \Pi $ is a semi-direct product: 
\begin{equation}
\Pi \ = \ \m{T}_{\Pi}   \rtimes 
\left ( \m{W}_{\Pi}  \times_{\{ \pm 1 \}} \m{S}_{\Pi}  \right ).  
\end{equation}  
The above three subgroups are naturally related through the homomorphism 
$\alpha$ to the three particular subgroups of $ \Gamma^+_F $ described earlier.
\par When restricted to $ \m{S} \subset \Gamma_F^+ $, the morphism $ \alpha $ produces
an isomorphism $ \m{S} \simeq \m{S}_{\Pi} $ sending $g(m, 0, 0, I) $ with   
$$ m = 
\left ( 
\begin{array}{cc}
a & b \\
c & d \\ 
\end{array}
\right ) \ $$
to the automorphism in $ \m{S}_{\Pi} $ associated to the matrix:  
$$ 
\left ( 
\begin{array}{cc}
a & -b \\
-c & d \\ 
\end{array}
\right ). $$
\par When restricted to $ \m{W} \subset \Gamma_F^+ $ the morphism $ \alpha $ produces
an isomorphism $ \m{W} \simeq \m{W}_{\Pi} $, which sends  $ g(I, 0, 0, f)$ to 
the automorphism induced by $f$ in $\m{W}_{\Pi}$. 
\par Finally, when restricted to $ \m{T} \subset \Gamma_F^+ $, 
the morphism $ \alpha $ produces
a surjective morphism $ \m{T} \simeq \m{T}_{\Pi} $ with kernel $ {\rm U}(N)_{\Zee}$.
If $ Q \colon \Lambda \rightarrow \Zee^2 $ is given by  
$$ Q(\gamma) = \left ( (\gamma, q_1), \ (\gamma, q_2) \right ), \ \ q_1, q_2 \in \Lambda $$
then $ \alpha $ represents the assignment $ g(I,Q,R,I) \rightarrow (q_2, -q_1) $.    
\par All three subgroups, $\m{S}_{\Pi}  $, $ \m{T}_{\Pi} $ and
$ \m{W}_{\Pi} $ are therefore entirely covered by the image of $ \alpha $. Since 
they generate $ \Pi $, the morphism $ \alpha $ is surjective. 
\par One verifies immediately that the map: 
\begin{equation}
\w{\Theta} \colon \Omega^+(F) \rightarrow \mathbb{H} \times \Lambda_{\Cee} , \ \ 
\w{\Theta} \left ( [a,b,c] \right ) = 
\left ( - \frac{a_1}{a_2} , \frac{c}{a_2} \right )   
\end{equation} 
is equivariant with respect to $ \alpha $. This leads to the Seifert fibration:  
\begin{equation} 
\Theta \colon \Gamma_F^+ \backslash \Omega^+(F) \ \rightarrow \ \Pi \backslash   
\left ( \mathbb{H} \times \Lambda_{\Cee} \right ) = \m{M}_{E, G}   
\end{equation} 
whose fibers are isomorphic to 
$ {\rm U}(N)_{\Cee} / {\rm U}(N)_{\Zee} $ and therefore are copies 
of $\Cee^* $. 
\par The identification $ (\ref{identif}) $ follows from the arguments above.  
\end{proof} \

\subsection{Automorphy Factors for the Parabolic Cover}
\noindent Let us remark that, based on the above arguments, one obtains a canonical isomorphism
\begin{equation}
\m{D}_{\Lambda} = \Gamma^+_F \backslash \m{B}^+(F) = \Gamma^+_F \backslash 
\Omega^+(F) / U(N)_{\Cee} \ \simeq \ \Pi \ \backslash \left ( \mathbb{H} \times \Lambda_{\Cee} \right ) = 
\m{M}_{E,G}   
\end{equation} 
identifying the Type II Mumford divisor $ \m{D}_{\Lambda} $ with the moduli space $ \m{M}_{E,G} $ of 
pairs of elliptic curves and flat $G$-bundles. 
Under this isomorphism, the parabolic cover $ (\ref{pparacover}) $ becomes the 
induced holomorphic Seifert $ \Cee^* $- fibration:
\begin{equation}
\label{pcc}
\Gamma_F^+ \ \backslash \Omega^+(F) \ \rightarrow \ \Pi \ \backslash 
\left ( \mathbb{H} \times \Lambda_{\Cee} \right ). 
\end{equation}
Our task in this section is to analyze the holomorphic type 
of $ (\ref{pcc}) $.
\par We use the following strategy. The base space of $ (\ref{pcc}) $ is a complex orbifold
$ \Pi  \backslash V $ where $ V = \mathbb{H} \times \Lambda_{\Cee} $. One can describe 
holomorphic $ \Cee^* $-fibrations over such a space in terms of equivariant line 
bundles over the cover $V$. These equivariant objects are line bundles 
$ \m{L} \rightarrow V $ where the action of the group $ \Pi $ on the base is given a 
lift to the fibers. All holomorphic line bundles over $ V $ are trivial and,
choosing a trivializing section, one obtains a lift of the action to fibers through 
a set of automorphy factors
$ (\varphi_a)_{a \in \Pi} $ with 
$ \varphi_a \in H^0(V, \m{O}^*_{V}) $ satisfying:
$$ \varphi_{ab}(x) = \varphi_a(b \cdot x) \varphi_b(x). $$
Such a set determines a 1-cocycle $ \varphi $ 
in $ Z^1(\Pi, H^0(V_G, \m{O}^*_{V})) $.
Two automorphy factors provide isomorphic fibrations on 
$ \m{M}_{E, \ G} $ if and only if they determine the same group cohomology class
in $ H^1(\Pi, H^0(V, \m{O}^*_{V})) $. To state this rigorously, 
there is a canonical map $ \phi $ entering the following exact sequence:
\begin{equation}
\{1\} \rightarrow H^1(\Pi, H^0(V, \m{O}^*_{V})) 
\stackrel{\phi}{\rightarrow} H^1(\m{M}_{E, \ G}, \ \m{O}^*_{\m{M}_{E, \ G}}) 
\stackrel{p^*}{\rightarrow} H^1(V, \ \m{O}^*_{V}) \simeq \{1\} .
\end{equation}   
\noindent We are going to write down explicitly a set of automorphy factors for 
fibration $ (\ref{pcc}) $. Since the modular group $ \Pi $ is generated by the 
three subgroups $ \m{S}_{\Pi} $, $\m{W}_{\Pi}  $ and $ \m{T}_{\Pi}  $ it 
will suffice to describe automorphy factors for elements in those subgroups. 
\par The first step towards computing the automorphy factors of $ (\ref{pcc}) $ 
is defining a holomorphic trivialization of the covering $ \Cee $-bundle:
\begin{equation}
\label{f3}
\w{\Theta} \colon \Omega^+(F) \rightarrow \mathbb{H} \times \Lambda_{\Cee}, \ \ \ 
\w{\Theta} 
\left ( [a,b,c] \right ) = \left ( - \frac{a_1}{a_2} , \frac{c}{a_2} \right ).
\end{equation}
Recall that this map provides the arithmetic recipe through which one can obtain
out of a given $ K3 $ period an elliptic curve $ E_{\tau} $ and a morphism 
$ \psi \colon \Lambda \rightarrow E_{\tau} $ which carries the holonomy information
of a flat $ G$-connection. Building a trivializing section for $ (\ref{f3}) $ amounts
then intuitively to finding a way to recover a $K3$ period out of geometric data
given by an elliptic curve endowed with a flat $G$-connection.
\par Surprisingly, such a method arises in string theory, precisely in the Narain 
construction (see \cite{narain1} \cite{narain2} \cite{ginsparg}) of the lattice of momenta 
related with toroidal compactification of heterotic strings. This construction leads 
one to consider the following map (see Appendix $ \ref{appendix} $ for details):
\begin{equation}
%\label{narainmap2}
\sigma_n \colon \mathbb{H} \times \Lambda_{\Cee}  \times \Cee \ \rightarrow \ 
\Omega^+(F)   
\end{equation}
\begin{equation}
\label{narain33} 
\sigma_n(\tau, z, u ) \ = \ {\rm exp} \left ( u \cdot N \right ). \left [ 
\ (-\tau, 1) , \ 
\frac{1}{2} \left ( 
\frac{(z,z) - (z,\bar{z} ) }{\bar{\tau} - \tau }, \ 
\frac{\bar{\tau} ( z,z) - \tau ( z , \bar{z} ) }{\bar{\tau} - \tau }  
\right ) , \ z \ \right ]. \ \  
\end{equation}
A brief analysis of the above formula reveals the following:
\begin{rem} \
\begin{enumerate}
\item The image of $ \sigma_n $ indeed lies in the indicated space since for 
any triplet $ (\tau, z, u) $, 
$$ \left ( \sigma_n (\tau, z,u ), \ \sigma_n(\tau, z,u) \right ) =  0 \ \ {\rm and} \ \ 
-i ( \ N \sigma_n(\tau, z,u) , \ \overline{\sigma_n(\tau, z,u)} \ ) = 
{\rm Im} \tau > 0.   $$  \item One has:
$$ \left ( \sigma_n(\tau, z,u), \ \overline{\sigma_n(\tau, z,u)} \right ) =  \ 
{\rm Im}(u)   \ {\rm Im}(\tau)  $$
and therefore, $ \sigma_n(\tau, z,u) $ is a $K3$ period for any $u \in \Cee$ with 
strictly positive imaginary part.
\item The map 
\begin{equation}
\label{narains}
(\tau, z) \rightarrow \sigma_n(\tau, z,0) 
\end{equation} 
makes a smooth section for the line bundle $ (\ref{f3}) $. 
\item When one factors out the action of $ {\rm U}(N)_{\Zee} $, application $ (\ref{narain33}) $   
provides a smooth trivialization for the induced $ \Cee^* $-bundle:
\begin{equation}
\label{th22}
\w{\Theta} \colon \Omega^+(F) / {\rm U}(N)_{\Zee} \ 
\rightarrow \ \mathbb{H} \times \Lambda_{\Cee}.  
\end{equation}
\end{enumerate}  
\end{rem} 
\noindent The above Narain trivialization has a major drawback !  
It is not holomorphic. Nevertheless, one can get around this problem and obtain a 
holomorphic trivialization by perturbing slightly the map $ (\ref{narain33}) $.
\par Note that the middle term in expression $ (\ref{narain33}) $ can
be rewritten:
$$ 
\frac{1}{2} \left ( 
\frac{(z,z) - (z,\bar{z} ) }{\bar{\tau} - \tau }, \ 
\frac{\bar{\tau} ( z,z) - \tau ( z , \bar{z} ) }{\bar{\tau} - \tau }  
\right ) \ = \ 
\frac{1}{2} \left ( 
\frac{(z,z) - (z,\bar{z} ) }{\bar{\tau} - \tau }, \
(z , z) + \tau \frac{(z,z) - (z,\bar{z} ) }{\bar{\tau} - \tau }   
\right ) \ = 
$$ 
$$ 
\ = \frac{1}{2}  \left ( 0, (z,z) \right ) \ + \ 
\frac{1}{2} \frac{(z,z) - (z,\bar{z} ) }{\bar{\tau} - \tau } \left ( 
1,\tau   
\right ) \ \ = \ 
\frac{1}{2}  \left ( 0, (z,z) \right ) \ + \ 
\frac{1}{2} \frac{(z,z) - (z,\bar{z} ) }{\bar{\tau} - \tau } T \left ( 
-\tau, 1   
\right )  
$$
Following the above equality, one can see the Narain section $ (\ref{narain33}) $ as:
$$ \sigma_n(\tau, z, 0 ) \ = \ {\rm exp} 
\left( \frac{1}{2} \frac{(z,z) - (z,\bar{z} ) }{\bar{\tau} - \tau } N \right ). 
\left [ 
\ (-\tau, 1) , \ 
\frac{1}{2} \left ( 
0, (z,z)  
\right ) , \ z \ \right ] \  \in \Omega^+(F) . \ \  $$ 
The second factor in the right-hand side term is holomorphic. This suggest the 
following perturbation:
\begin{equation}
\label{narainmap2}
\sigma \colon \mathbb{H} \times \Lambda_{\Cee}  \ \rightarrow \ 
\Omega^+(F)   
\end{equation}
$$ \sigma(\tau, z) \ = \ {\rm exp} 
\left( - \frac{1}{2} \frac{(z,z) - (z,\bar{z} ) }{\bar{\tau} - \tau } N \right ). 
\sigma_n(\tau, z) \ = \ 
\left [ 
\ (-\tau, 1) , \ 
\frac{1}{2} \left ( 
0, (z,z)  
\right ) , \ z \ \right ].  $$
We call this the ${\bf perturbed \ Narain \ map}$. 
One can immediately check that:
\begin{teo}
The perturbed Narain map $ (\ref{narainmap2}) $ is a holomorphic section for 
the line bundle:
\begin{equation}
\label{th3333}
\w{\Theta} \colon \Omega^+(F) \rightarrow \mathbb{H} \times \Lambda_{\Cee}, \ \ 
\w{\Theta} \left ( [a,b,c] \right ) \ = \ \left ( -\frac{a_1}{a_2}, \ \frac{c}{a_2} 
\right ). 
\end{equation}
It descends to a holomorphic section for the $\Cee^* $-fibration $ (\ref{th22}) $, 
providing therefore a holomorphic trivialization:
$$ \mathbb{H} \times \Lambda_{\Cee} \times \Cee^* \simeq 
\Omega^+(F)/ {\rm U}(N)_{\Zee} , \ \ \  (\tau, z,u) 
\rightarrow u \cdot \sigma(\tau, z).   $$
The $ \Cee^* $-action in the right-hand side expression represents the action 
of $  {\rm U}(N)_{\Cee}  /   {\rm U}(N)_{\Zee} $ upon 
$ \Omega^+(F) / {\rm U}(N)_{\Zee}  $.
\end{teo} 
\noindent We are now in position to compute the automorphy factors of the parabolic cover map:
\begin{equation}
\Theta \colon \Gamma_F^+ \ \backslash \Omega^+(F) \ \rightarrow \ \Pi  \backslash 
\left ( \mathbb{H} \times \Lambda_{\Cee} \right ). 
\end{equation}
In order to obtain a set of factors, one needs to analyze the variation
of the perturbed Narain map $ \sigma $ under the action of the modular
group $ \Pi $. 
\vspace{.1in}
\begin{lem}
\label{ll1} 
Assume $ (q_1, q_2 ) \in \m{T}_{\Pi} $. Then:
\begin{equation}
\label{f5}
\sigma(\tau, z+ q_1 + \tau q_2 ) \ = \ e^{\pi i \left ( 
2(q_2, z)  +\tau (q_2, q_2) 
\right )} \cdot g(I,Q,R,I) . \sigma(\tau, z) 
\end{equation}  
where $Q \in {\rm Hom}_{\Zee} 
\left ( \Lambda , \Zee^2 \right ) $ is given by $ Q(\gamma) = \left ( 
- (\gamma, q_2), (\gamma, q_1)  \right ) $ and $ R \in {\rm End}(\Zee^2) $ 
with $ R+R^t = QQ^t $
\end{lem}
\begin{proof}
We perform the computations in $ \Omega^+(F) $. 
$$ \sigma(\tau, z+ q_1 + \tau q_2 ) \ = \ \left [ (-\tau,1), \ \frac{1}{2} 
\left ( 0, (z,z) + (q_1,q_1) + \tau^2 (q_2,q_2) + 2(z,q_2) + 2 \tau (z, q_2) + 
2 \tau (q_1, q_2) \right ) , \  
z + q_1 + \tau q_2 \right ]. $$
On the other hand, 
$$  g(I,Q,R,I).\sigma(\tau, z)  \ = \ \left [ 
(-\tau, 1), \ R(-\tau, 1) + \frac{1}{2} (0, (z,z)) + Qz , \ z+ Q^t(-\tau, 1) \right ]. $$
But $ Qz = ( -z q_2 , z q_1 ) $ and $ Q^t(-\tau, 1) = q_1 + \tau q_2 $. 
Moreover, one can 
see that: 
$$ R(-\tau, 1) = \left ( -(q_1, q_2) - \frac{1}{2} \tau (q_2,q_2) , \ 
\frac{1}{2} (q_1,q_1) 
\right ) +  A(-\tau, 1) $$
where $ A \in {\rm End}(\Zee^2) $ skew-symmetric. One obtains then the following 
equality in $ \Omega^+(F) $:
\begin{equation}
\sigma(\tau, z+ q_1 + \tau q_2 ) \ = \ {\rm exp} \left ( \left (  
(q_2, z) + \frac{1}{2} \tau (q_2, q_2) + \alpha \right ) N \right ) \ g(I,Q,R,I) . \sigma(\tau, z) 
\end{equation}
with  $ \alpha \in \Zee $. After 
factoring out the $ {\rm U}(N)_{\Zee} $-action, one is led to 
$ (\ref{f5}) $. 
\end{proof}  
\vspace{.1in} 
\begin{lem}
\label{ll2}
Assume $ 
\left ( 
\begin{array}{cc}
a & b \\
c & d \\ 
\end{array}
\right ) \ \in {\rm SL}_2(\Zee) $. Then:
\begin{equation}
\label{f10}
\sigma \left ( \frac{a \tau + b}{c \tau + d}, \ \frac{z}{c \tau + d } \right ) \ = \ 
e^{ 
\left ( 
- \frac{\pi i c(z,z)}{c \tau + d } 
\right )} 
\cdot g(m,0,0,I).\sigma(\tau, z),  \ \ {\rm where} \  m = 
\left ( 
\begin{array}{cc}
a & -b \\
-c & d \\ 
\end{array}
\right ).
\end{equation}
\end{lem}
\begin{proof} 
As in the previous lemma, we write the calculations in $ \Omega^+(F) $. One has:
\begin{equation}
\sigma \left ( \frac{a \tau + b}{c \tau + d}, \ \frac{z}{c \tau + d } \right ) \ = \ 
\left [ 
\left ( - \frac{a \tau + b}{c \tau + d}, \ 1 \right ), \
\frac{1}{2} \left ( 0, \ \frac{(z,z)}{(c \tau + d)^2 } \right ), \ 
\frac{z}{c \tau + d }   
\right ]. 
\end{equation}
%$$ 
%= \ \left [ 
%\left ( - ( a \tau + b ) , \ c\tau + d  \right ), \
%\frac{1}{2} \left ( 0, \ \frac{z^2}{(c \tau + d) } \right ), \ z   
%\right ] .
%$$
In the same time:
\begin{equation}
g(m,0,0,I) . \sigma(\tau, z) \ = \ 
\left [ 
m(-\tau, 1), \ \w{m} \left ( \frac{1}{2} (0, (z,z)) \right ) , \ 
z 
\right ] \ = \ 
\end{equation}
$$ 
= \ \left [ 
\left ( - ( a \tau + b ) , \ c\tau + d  \right ), \
\frac{1}{2} \left ( c(z,z), d(z,z) \right ), \ z   
\right ] \ = \ 
\left [ 
\left ( - \frac{a \tau + b}{c \tau + d}, \ 1 \right ), \
\frac{1}{2} \left ( 
\frac{c(z,z)}{c \tau + d} , \ \frac{d(z,z)}{c \tau + d } \right ), \ 
\frac{z}{c \tau + d }   
\right ].
$$
Comparing the two formulas we get the following identity in $ \Omega^+(F) $:
\begin{equation}
\label{f8}
\sigma \left ( \frac{a \tau + b}{c \tau + d}, \ \frac{z}{c \tau + d } \right ) \ = \ 
{\rm exp}  
\left ( 
\left (  - \frac{1}{2} \frac{c(z,z)}{c \tau + d } \right ) N \right )
 \ g(m,0,0,I) . \sigma(\tau, z).
\end{equation}
Factoring out the $ {\rm U}(N)_{\Zee} $-action, one obtains $ (\ref{f10}) $.
\end{proof} 
\vspace{.1in}
\noindent We can state then:
\begin{teo}
\label{mm}
Let $  \left ( \varphi_{g} \right ) _{g \in K} $ be the automorphy factors 
of parabolic cover $ \Cee^* $-fibration:
\begin{equation}
\label{f22}
\Gamma_F^+ \  \backslash \ \Omega^+(F) \ \rightarrow \ \Pi \ \backslash  
\left ( \mathbb{H} \times \Lambda_{\Cee} \right ) 
\end{equation} 
associated to the trivialization generated by $ \sigma $. Then:
\begin{enumerate}
\item $ \varphi_g(\tau, z) \ = \ e^{-\pi i \left ( 
2(q_2, z) +\tau (q_2, q_2) \right )} $ for 
$ g = (q_1, q_2)  \in \m{T}_{\Pi}  $.
\item $ \varphi_g(\tau, z) \ = \ e^{ 
\left ( 
\frac{\pi i c(z,z)}{c \tau + d } 
\right )} $ for $ g = \left ( 
\begin{array}{cc}
a & b \\
c & d \\ 
\end{array}
\right ) \in \m{S}_{\Pi}  $. 
\item $ \varphi_g(\tau, z) \ = \ 1 $ for $ g \in \m{W}_{\Pi}  $.
\end{enumerate}
\end{teo}
\begin{proof} 
The first two expressions are direct consequences of Lemmas $ \ref{ll1} $ and 
$ \ref{ll2} $. The fact that $ \varphi_g(\tau, z) \ = \ 1 $ for $ g \in \m{W}_{\Pi}  $.
follows from:
$$ \sigma ( \tau, f(z)) \ = \ \left [ 
\ (-\tau, 1) , \ 
\frac{1}{2} \left ( 
0, \left ( f(z),f(z) \right )  
\right ) , \ f(z) \ \right ] \ = \ $$ 
$$ = \ \left [ 
\ (-\tau, 1) , \ 
\frac{1}{2} \left ( 
0, (z,z)  
\right ) , \ f(z) \ \right ] \ = \ g(I,0,0,f). \sigma(\tau,z) $$
for any $ f \in O(\Lambda) $.  
\end{proof} 
\vspace{.1in}
\noindent The three subgroups 
$ \m{T}_{\Pi}  $, $ \m{S}_{\Pi}  $ and $\m{W}_{\Pi} $ generate
the entire modular group $ \Pi $. Therefore, the above automorphy factors are enough 
to characterize completely the holomorphic type of fibration $ (\ref{f22}) $.

\subsection{Theta Function Interpretation and Relation to Heterotic String Theory}  
Given the particular automorphy factor expressions computed in the previous section, 
one can provide for the parabolic cover $\Cee^* $-fibration  $ (\ref{f22}) $ 
a theta function interpretation.  
\par Let $ \mathbb{H} \times \Lambda_{\Cee} $ be the orbifold cover of 
$ \m{M}_{E,G} $. Since $  \Lambda $ is positive definite, unimodular and even, 
there is an associated holomorphic theta function 
(see \cite{kac} \cite{mumford2} for details) :
\begin{equation}
\Theta_{\Lambda} \colon \mathbb{H} \times \Lambda_{\Cee} \rightarrow \mathbb{C}, \ \ \ \ 
\Theta_{\Lambda}(\tau, \ z) = \ \sum_{\gamma \in \Lambda} \ 
e^{ \pi i ( 2 (z,\gamma) + \tau (\gamma, \gamma))}.  
\end{equation}
The pairing appearing above represents the bilinear complexification of the 
integral pairing on $ \Lambda $. The $\Lambda$-character 
function can be written then as a quotient of $ \Theta_{\Lambda} $:
\begin{equation}
\label{char}
B_{\Lambda} \colon \mathbb{H} \times \Lambda_{\Cee} \rightarrow \mathbb{C}, \ \ \ \ B_{\Lambda}
(\tau, \ z) = 
\frac{\Theta_{\Lambda}(\tau, \ z)}{\eta(\tau)^{16}}. 
\end{equation}
Here, $ \eta $ is Dedekind's eta function:
$$ \eta(\tau) = e^{\pi i \tau /12} \prod_{m=1}^{\infty} \ \left ( 
1- e^{2 \pi i m \tau} \right ),$$
which is an automorphic form of weight $ 1/2$ and multiplier system 
given by a group homomorphism 
$ \chi \colon {\rm SL}_2(\Zee) \rightarrow \Zee / 24 \Zee $, in the
sense that \cite{sig}:
$$ \eta ( \gamma \cdot \tau ) \ = \ 
\chi(\gamma) \  \sqrt{c \tau+d} \ \eta(\tau) \ \ {\rm for} \ \ 
\gamma = \left (
\begin{array}{cc}
a & b \\
c & d \\
\end{array}
\right ) \in {\rm SL}_2(\Zee). 
$$
The character terminology for $ ( \ref{char} ) $ is justified by its role in the
representation theory of infinite-dimensional Lie algebras. The function 
$ B_{\Lambda} $ is the
zero-character of the level $l=1$ basic highest weight representation of the Kac-Moody
algebra associated to $G$ (see \cite{kac} for details).   
\par According to \cite{kac}, the character function 
$ B_{\Lambda} $ obeys the following transformation properties: 
\begin{prop}  
\label{kac} 
Under the action of the modular group $ \Pi $, the character function
$ (\ref{char}) $ transforms as :
$$  B_{\Lambda} \left ( g \cdot ( \tau, z) \right )  \ = \ 
\varphi^{{\rm ch}}_g (\tau, z) \cdot B_{\Lambda}( \tau, z). $$
The factors $ \varphi^{{\rm ch}}_g, g \in \Pi $ can be described as:
\begin{itemize}
\item $ \varphi ^{{\rm ch}}_{(q_1, q_2)} (\tau, z) = 
e^{\pi i(-2 (q_2 z)-\tau (q_2 q_2) )}$  for $ (q_1, q_2) \in \m{T}_{\Pi} $.  
\item $ \varphi ^{{\rm ch}}_{m} = e^{\frac{\pi i c (z z)}{c \tau + d}} $ 
 for $ m = \left (
\begin{array}{cc}
a & b \\
c & d \\
\end{array} \right ) 
\in \m{S}_{\Pi}  $. 
\item $ \varphi^{{\rm ch}}_{w} = 1 $ for $ w \in \m{W}_{\Pi}  $.
\end{itemize} 
\end{prop}  
\noindent The holomorphic function $ B_{\Lambda} $ descends therefore to a section 
of a $ \mathbb{C} $-fibration:
\begin{equation}
\label{thetaline}
\m{Z} \rightarrow \Pi  \backslash  
\left ( \mathbb{H} \times \Lambda_{\Cee} \right ) = \m{M}_{E, \ G} 
\end{equation}
with automorphy factors  $ \varphi^{{\rm ch}}_g $ described above. We call this
the character fibration.
\par One can compare then the character fibration $ (\ref{thetaline}) $ with the
parabolic cover $ \Cee^* $-fibration:
\begin{equation}
\label{f2222}
\Theta \colon \Gamma_F^+ \backslash  \Omega^+(F) \ \rightarrow \ \Pi  \backslash  
\left ( \mathbb{H} \times \Lambda_{\Cee} \right )  = \m{M}_{E, \ G} 
\end{equation} 
analyzed in the previous section. A look at Theorem $ \ref{mm} $ and Proposition $ \ref{kac} $ 
is enough to convince us that the two holomorphic fibrations are defined through identical 
automorphy factors. Therefore: 
\begin{teo}
\label{mainth} 
The parabolic cover fibration  
\begin{equation}
\label{fibr66} 
\Theta \colon \Gamma_F ^+\backslash  \Omega^+(F) \ \rightarrow \ \Pi  \backslash 
\left ( \mathbb{H} \times \Lambda_{\Cee} \right ) = \m{M}_{E, \ G} 
\end{equation} 
is holomorphically isomorphic  to the character fibration $ (\ref{thetaline}) $ 
with the zero section removed.   
\end{teo}   
\par We conclude the section by placing the outcome of Theorem $ \ref{mainth} $ in 
connection with the parabolic compactification construction presented in 
section $ \ref{first} $, and comparing the resulted structure to the classical moduli 
spaces of eight dimensional heterotic string theory.   
\par Recall that, up to isomorphism, there exist only two even, positive-definite 
and unimodular 
lattices $\Lambda$ of rank $16$. To each choice of $ \Lambda $ one can associate 
a corresponding Lie group $G$. For $ \Lambda_1 = E_8 \times E_8 $ one sets
$ G_1 = \left ( E_8 \times E_8 \right )  \rtimes \Zee_2 $. If 
$ \Lambda_2 = \Gamma_{16} $ then $ G_2 = {\rm Spin}(32)/\Zee_2 $. The moduli space $ \m{M}_{K3} $ 
of $K3$ surfaces with section admits a partial compactification $ \overline{\m{M}_{K3} } $ 
obtained by adding two distinct divisors at infinity $ \m{D}_{\Lambda_i} $. Each point 
on $ \m{D} _{\Lambda_i} $ can be identified with an equivalence class of elliptic Type II 
stable $K3$ surfaces with section, in $\Lambda_i$-category, and with 
an isomorphism class of a pair $ (E,A) $ consisting of an 
elliptic curve $E$ and a flat $G_i$-connection $A$. The 
correspondence gives a natural holomorphic isomorphism:
\begin{equation}
\label{adas} 
\m{D}_{\Lambda_i} \ \simeq \  \m{M}_{E,G_i}  
\end{equation} 
where $ \m{M} _{E,G_i} $ is the $17$-dimensional moduli space of pairs of elliptic curves
and flat $G_i$-bundles\footnote{Again, if $ G = {\rm Spin}(32)/\Zee_2 $ only connection
liftable to $ {\rm Spin}(32) $-connection are considered}. 
\par As explained in $ \ref{arith} $, in each of the two cases, one has the parabolic cover 
\begin{equation} 
 \m{P}_{\Lambda_i} \stackrel{p}{\rightarrow} \m{M}_{K3}    
\end{equation}     
modeling the projection $ \Gamma_F \backslash \Omega \rightarrow 
\Gamma \backslash \Omega $ where $ \Gamma_F $ is the stabilizer in $ \Gamma $
of a rank two isotropic sub-lattice of $ \mathbb{L}_o $ determining $ \Lambda_i $. 
Moreover, the space $ \m{P}_{\Lambda_i} $ fibers holomorphically over the corresponding 
divisor:
\begin{equation}
\label{adas1} 
\m{P}_{\Lambda_i} \rightarrow \m{D}_{\Lambda_i} 
\end{equation}         
with all fibers being copies of $ \mathbb{C}^* $. Theorem $ \ref{mainth} $ shows that,
under identification $ (\ref{adas}) $, the above fibration is the character
fibration of $ \Lambda_i $ with the zero-section removed. That allows one to 
holomrphically identify  
$ \m{P}_{\Lambda_i} $ with the total space of the character $ \Cee^* $-fibration.
\par Turning our attention to the heterotic side of the duality, it was shown in 
\cite{clingher} (see Theorems $ \ref{het1} $ and $ \ref{het2} $ in 
section $ \ref{intro} $) that the moduli space $ \m{M}^{G_i}_{{\rm het}} $ of 
classical vacua for heterotic string theory compactified over the torus is 
holomorphically isomorphic to the same total space of the 
character $\Cee^*$-fibration corresponding to the lattice $ \Lambda_i$. Corroborating 
these facts to Theorem $ \ref{mainth} $, one obtains a holomorphic isomorphism of 
$\Cee^*$-fibrations:
\begin{equation}
\label{lastdiag}
\xymatrix{
\m{P}_{\Lambda_i} \ar [r] ^{\simeq} \ar [d] &   \m{M}^{G_i}_{{\rm het}} \ar [d] \\
 \m{D}_{\Lambda_i} \ar [r]^{\simeq} & \m{M}_{E,G_i} & \\ 
}
\end{equation}
which can be seen as an identification between the parabolic cover space 
$\m{P}_{\Lambda_i}$ and the classical moduli space of eight-dimensional 
heterotic string theory with group $G_i$.   
\par But, as described in section $ \ref{arith} $, there exists an open 
set $ \m{V} $, punctured tubular neighborhood of the divisor 
$ \m{D}_{\Lambda_i} $ in $ \m{M}_{K3} $ such that the pre-image 
$ p^{-1}(\m{V}) $ in $ \m{P}_{\Lambda_i} $ is a tubular neighborhood 
of the zero section in $ (\ref{adas1}) $ and the projection 
$ p^{-1}(\m{V}) \rightarrow \m{V} $ is an isomorphism. 
%This fact allows one to describe 
%the structure of $ \m{M}_{K3} $ 
%near the boundary divisor $ \m{D}_{\Lambda_i} $ by analyzing
%the parabolic cover $ \m{P}_{\Lambda_i} $ in the neighborhood of the zero-section
%of the Seifert fibration 
%$ (\ref{adas1}) $. But, by Theorem $ \ref{mainth} $, the fibration $ (\ref{adas1}) $ 
%is the character fibration. 
This fact allows us to conclude that:
\begin{teo} (F-Theory/Heterotic String Duality in Eight Dimensions) \
\label{duality}

\noindent There exists a holomorphic isomorphism between an open neighborhood of 
$ \m{M}_{K3} $ near the divisor $ \m{D}_{\Lambda_i} $ and an
open neighborhood of $ \m{M}^{G_i}_{{\rm het}} $ near the zero-section of 
the left fibration in $(\ref{lastdiag})$. 
\end{teo}
The open neighborhood of $ \m{M}^{G_i}_{{\rm het}} $ in the above statement 
corresponds to large volumes of the elliptic curve. Hence, the two regions identified by 
Theorem $ \ref{duality} $ are exactly the sectors that physics predicts should 
closely resemble each other.

\section{Appendix: Narain Construction} 
\label{appendix} 
The parameter fields for $8$-dimensional heterotic string theory are, 
after Narain \cite{ginsparg} \cite{narain1},  triplets 
$ (A,g,B) $ consisting of a flat $G$-connection, a flat metric and a constant 
anti-symmetric 2-tensor $B$, 
all defined over a two-torus $E$. The Lie group $G$ involved is either 
$ E_8 \times E_8 \rtimes \Zee_2 $ or $ {\rm Spin}(32)/\Zee_2 $. 
\par One usually describes a flat torus as a quotient:
$$ E = \ \ \mathbb{R}^2 / U  $$
of the Euclidean space $ \mathbb{R}^2 $ through a rank two lattice 
$ U = \Zee e_1 \oplus \Zee e_2 $. In this way, $E$ inherits a flat metric, which in 
turn generates a volume $ v \in \mathbb{R}^*_+ $ and a complex structure parameterized 
by $ \tau \in \mathbb{H} $. These parameters are obtained as: 
$$
v \ = \ \sqrt{g_{11}g_{22} - g_{12}^2 }    
$$                
$$
\tau \ = \ \frac{g_{12} + v \cdot i  }{g_{11} }. 
$$  
where $ g_{ij} = e_i \cdot e_j $. 
\par A flat $G$-connection on $E$ is, formally, a morphism 
$$ A \ \colon \ U \ \rightarrow \ \Lambda_{\mathbb{R} }. $$ 
The lattice $ \Lambda $ is the coroot lattice of $G$ if 
$ G = E_8 \times E_8 \rtimes \Zee_2 $ and the lattice of a maximal torus of
$G$ if $ G= {\rm Spin}(32)/\Zee_2 $. As is the standard procedure, one parameterizes
holomorphically these flat connections by taking:
$$ z= z_2 - \tau z_1 \ \in \Lambda_{\Cee} , \ \ {\rm where} \ \ 
z_i \colon = A(e_i) \ \in \Lambda_{\mathbb{R}}. $$
\par The last ingredient, the B-field is seen as a two-form $ B = b \ (e_1^* ) \wedge (e_2^*) $
with $ b \in \mathbb{R} $. The B-field holonomy along $ E $ is given by 
$$ {\rm exp} \left ( i \int_E B \right ) \ = \ {\rm exp} \left ( i b v \right ). $$
\par One considers then the space:
$$  \mathbb{R}^{2,18} \ = \ \mathbb{R}^2 \oplus \mathbb{R}^2 \oplus \Lambda_{\mathbb{R}} 
$$ 
endowed with the inner product:
$$ ( x,y,z) . (x',y',z') = x.x' -y.y' - z.z'. $$
The lattice of momenta \cite{ginsparg}, denoted $\m{L}_{(A,g,B)} $, 
associated to a heterotic triplet $ (A,g,B) $ 
is obtained as the image of the map: 
\begin{equation}
\label{mapnn} 
\varphi_{(A,g,B)} \ \colon U \oplus U^* \oplus \Lambda \ \rightarrow \ \mathbb{R}^2 \oplus \mathbb{R}^2 
\oplus \Lambda_{\mathbb{R}}
\end{equation}             
$$ 
\varphi_{(A,g,B)} (w, p, l)  \ = \ $$
$$
= \ \left ( 
\frac{1}{2} p - bTw -\frac{1}{2} A^t l  -\frac{1}{4} A^tAw - w, \ 
\frac{1}{2} p - bTw -\frac{1}{2} A^t l -\frac{1}{4} A^tAw + w, \
Aw + l    
\right ) .
$$
Here $ T \colon \mathbb{R}^2 \rightarrow \mathbb{R}^2 $ is the anti-self adjoint
morphism $ T(x_1,x_2) = (x_2,-x_1) $. One checks that, in the above formulation, the image 
$$ \m{L}_{(A,g,B)} \ \colon = \ {\rm Im}  \left ( \varphi_{(A,g,B)}\right )  $$ 
with the induced inner product forms a lattice of rank $20$ embedded in the ambient space 
$ \mathbb{R}^{2,18} $. The lattice $ \m{L}_{(A,g,b)} $ is
isomorphic to $ H \oplus H \oplus \left ( - \Lambda \right )  $. A basis 
underlying this decomposition is given by:
\begin{equation} 
\label{l4} 
F_i \ \colon = \ \varphi_{(A,g,B)}( - e_i,0 , 0) \ = \ \left ( 
bTe_i + \frac{1}{4} A^tAe_i + e_i, \ 
bTe_i +\frac{1}{4} A^tAe_i - e_i, \
- Ae_i  \right )
\end{equation} 
\begin{equation} 
\label{l5} 
F^*_i \ \colon = \ \varphi_{(A,g,B)}(0,e^*_i , 0) \ = \ \left (
\frac{1}{2}e^*_i , \ \frac{1}{2}e^*_i, \ 0 
\right ) 
\end{equation} 
\begin{equation} 
\label{l6} 
 L \ \colon = \ \varphi_{(A,g,B)} (0,0 , l) \ = \ \left (
-\frac{1}{2}A^t l , \ -\frac{1}{2}A^t l, \ l
\right ). 
\end{equation}  
It satisfies: 
$$ F_i . F_j = 0 , \ \ F_i^* . F_j^* = 0 , \ \ F_i . F^*_j = \delta_{ij} $$ 
$$ F_i .  L = F^*_j . L = 0, \ \ L. L' = - \ l.l' . $$ 
One is interested in the behavior of
the oriented positive 2-plane $ \mathbb{R}^2 \subset \mathbb{R}^{2,18} $ with respect
to the lattice $ \m{L}_{(A,g,B) } $. Let us imagine that the lattice
$ \m{L} $ remains fixed and the oriented $ \mathbb{R}^2 $ is varying inside 
$ \m{L} \otimes \mathbb{R} $ parameterized by the heterotic variables. This 
provides an assignment:
\begin{equation}
\label{blah} 
\left \{ {\rm heterotic \ parameters} \ (A,g,B) \right \} \ \rightarrow \ 
O(2,18) / SO(2) \times O(18).     
\end{equation}        
Moreover, the target space in $ (\ref{blah}) $ has a natural holomorphic structure. One
can equivalently regard positive, oriented, two-planes in $ \m{L} \otimes \mathbb{R} $ as
complex lines $ \omega \subset \m{L} \otimes \Cee $ satisfying $ \omega.\omega =0 $ and 
$ \omega. \overline{\omega} > 0 $. There is then a bijective correspondence:
$$ 
O(2,18) / SO(2) \times O(18) \ \simeq \ 
\left \{ \ \omega \in \mathbb{P}\m{L}_{\Cee} \ \vert \ \omega.\omega =0, \ 
\omega. \overline{\omega} > 0   \ \right \} 
$$    
and the map $ (\ref{blah}) $ can be interpreted as sending triplets of heterotic
parameters to the $18$-dimensional complex period domain $ \Omega $ of section 
$ \ref{per1}$. 
\par One can describe explicitly this map. Let $ (A,g,B) $ be a heterotic triplet
determining:
$$ (\tau, z, v, b)  \ \in \ \mathbb{H} \times \Lambda_{\Cee} \times \mathbb{R}^*_+ 
\times \mathbb{R}. $$
Then, the complex line $ \omega $ is generated by:
\begin{equation}
\label{l10}
\omega \ = \ \sum \alpha_i F_i \ + \ \sum \beta_j F^*_j \ + \ \gamma 
\end{equation}  
with 
$$ \alpha_1 \ = \ -\tau, \ \  \alpha_2 \ = \ 1, \ \ \gamma = z $$
$$ \beta_1 \ = \ - 2(bv+iv) \ + \ 
\frac{(z,z) - (z,\bar{z}) }{2(\bar{\tau} - \tau)}  \       
$$
$$ \beta_2 \ = \ -2\tau(bv+iv) \ + \ 
\frac{\bar{\tau} (z,z) - \tau (z,\bar{z}) }{2(\tau -\bar{\tau})}.  \       
$$
\par Take the decomposition 
$ \m{L} = H \oplus H \oplus \left ( - \Lambda \right )  $ 
with a basis for $ H \oplus H $ given by $ \{ F_1, F_2, F^*_1, F_2^* \} $. 
The inner product on $ \m{L}_{\Cee} $ appears as:
$$ (a,b,c).(a',b',c') \  = \ (a,b') + (b,a') - (c,c'). $$
Let $ N \in {\rm End}(\m{L}) $ be the nilpotent anti-self adjoint endomorphism
$$ N(a,b,c) \ = \ (0,Ta,0) $$
and let $ {\rm exp}(tN)  =  I + t N $ be its exponential. The Narain correspondence
between heterotic parameters and period complex lines in  $ \Omega $ appears then
as:
\begin{equation}
\label{l100} 
\sigma_n \ \colon \ \mathbb{H} \times \Lambda_{\Cee} \times \mathbb{H}  
\ \rightarrow \ \Omega 
\end{equation}  
$$ \sigma_n(\tau,z,u) \ = \ \ {\rm exp} \left ( -2u \cdot N \right ) \left [ 
\ (-\tau, 1) , \ 
\frac{1}{2} \left ( 
\frac{(z,z) - (z,\bar{z} ) }{\bar{\tau} - \tau }, \ 
\frac{\bar{\tau} ( z,z) - \tau ( z , \bar{z} ) }{\bar{\tau} - \tau }  
\right ) , \ z \ \right ] \ . $$
The complex variable $u$ represents $ bv+iv  \in \Cee $. It is clear that 
$ (\ref{l100}) $ is not holomorphic. However one can move the non-holomorphic 
part of $ (\ref{l100}) $ to the exponential. Indeed:
$$ \frac{\bar{\tau} ( z,z) - \tau ( z , \bar{z} ) }{\bar{\tau} - \tau } \ = \ 
\tau \ \frac{(z,z) - (z,\bar{z} ) }{\bar{\tau} - \tau } + (z.z) $$ 
and therefore, one can rewrite:
$$ \sigma_n(\tau,z,u) \ = \ \ {\rm exp} \left ( \left ( -2u + 
\frac{(z,z) - (z,\bar{z} ) }{2(\bar{\tau} - \tau) } \right ) 
\cdot N \right ) \left [ 
\ (-\tau, 1) , \ 
\frac{1}{2} \left ( 0, \ ( z,z) \right ) , \ z \ \right ] \ . $$   
\vspace{.1in}
%
%
%  BIBLIOGRAPHY 
%
%
%\bibliography{parab.bib}

\begin{thebibliography}{10}

\bibitem{ash}
A.~Ash, D.~Mumford, M.~Rapoport, and Y.~Tai.
\newblock {\em Smooth {C}ompactification of {L}ocally {S}ymmetric {V}arieties}.
\newblock Math. Sci. Press, Brookline, Mass., 1975.
\newblock Lie Groups: History, Frontiers and Applications, Vol. IV.

\bibitem{morrison1}
Paul~S. Aspinwall and David~R. Morrison.
\newblock String {T}heory on ${K}3$ surfaces.
\newblock In {\em Mirror symmetry, II}, pages 703--716. Amer. Math. Soc.,
  Providence, RI, 1997.

\bibitem{bb}
W.~L. Baily, Jr. and A.~Borel.
\newblock Compactification of {A}rithmetic {Q}uotients of {B}ounded {S}ymmetric
  {D}omains.
\newblock {\em Ann. of Math. (2)}, 84:442--528, 1966.

\bibitem{bpv}
W.~Barth, C.~Peters, and A.~Van~de Ven.
\newblock {\em Compact {C}omplex {S}urfaces}.
\newblock Springer-Verlag, Berlin, 1984.

\bibitem{carl}
James~A. Carlson.
\newblock Extensions of {M}ixed {H}odge {S}tructures.
\newblock In {\em Journ\'ees de G\'eometrie Alg\'ebrique d'Angers, Juillet
  1979/Algebraic Geometry, Angers, 1979}, pages 107--127. Sijthoff \&\
  Noordhoff, Alphen aan den Rijn, 1980.

\bibitem{clingher}
Adrian Clingher.
\newblock Heterotic {S}tring {D}ata and {T}heta {F}unctions.
\newblock {\em PhD Thesis}, Columbia University.

\bibitem{demazure}
Michel Demazure and Henry~Charles Pinkham, editors.
\newblock {\em S\'eminaire sur les {S}ingularit\'es des {S}urfaces}, volume 777
  of {\em Lecture Notes in Mathematics}.
\newblock Springer, Berlin, 1980.

\bibitem{freed}
Daniel~S. Freed.
\newblock Dirac {C}harge {Q}uantization and {G}eneralized {D}ifferential
  {C}ohomology.
\newblock In {\em Surveys in differential geometry}, Surv. Differ. Geom., VII,
  pages 129--194. International Press, Somerville, MA, 2000.

\bibitem{frthesis}
Robert Friedman.
\newblock Hodge {T}heory, {D}egenerations and the {G}lobal {T}orelli {P}roblem.
\newblock {\em Thesis, Harvard University}, 1981.

\bibitem{friedman1}
Robert Friedman.
\newblock Global {S}moothings of {V}arieties with {N}ormal {C}rossings.
\newblock {\em Ann. of Math. (2)}, 118(1):75--114, 1983.

\bibitem{friedman}
Robert Friedman.
\newblock A {N}ew {P}roof of the {G}lobal {T}orelli {T}heorem for ${K}3$
  {S}urfaces.
\newblock {\em Ann. of Math. (2)}, 120(2):237--269, 1984.

\bibitem{friedman3}
Robert Friedman.
\newblock The {P}eriod {M}ap at the {B}oundary of {M}oduli.
\newblock In {\em Topics in transcendental algebraic geometry (Princeton, N.J.,
  1981/1982)}, pages 183--208. Princeton Univ. Press, Princeton, NJ, 1984.

\bibitem{rmw}
Robert Friedman, John Morgan, and Edward Witten.
\newblock Vector {B}undles and {${\rm F}$} theory.
\newblock {\em Comm. Math. Phys.}, 187(3):679--743, 1997.

\bibitem{delpezzo}
Robert Friedman and John~W. Morgan.
\newblock Principal {H}olomorphic {B}undles over {E}lliptic {C}urves {IV}: del
  {P}ezzo {S}urfaces.
\newblock {\em in preparation}.

\bibitem{mf1}
Robert Friedman, John~W. Morgan, and Edward Witten.
\newblock Principal ${G}$-bundles over {E}lliptic {C}urves.
\newblock {\em Math. Res. Lett.}, 5(1-2):97--118, 1998.

\bibitem{ginsparg}
Paul Ginsparg.
\newblock On {T}oroidal {C}ompactification of {H}eterotic {S}uperstrings.
\newblock {\em Phys. Rev. D (3)}, 35(2):648--654, 1987.

\bibitem{gs}
Phillip Griffiths and Wilfried Schmid.
\newblock Recent {D}evelopments in {H}odge {T}heory: {A} {D}iscussion of
  {T}echniques and {R}esults.
\newblock In {\em Discrete subgroups of Lie groups and applicatons to moduli
  (Internat. Colloq., Bombay, 1973)}, pages 31--127. Oxford Univ. Press,
  Bombay, 1975.

\bibitem{griff}
Phillip~A. Griffiths.
\newblock Periods of {I}ntegrals on {A}lgebraic {M}anifolds: {S}ummary of
  {M}ain {R}esults and {D}iscussion of {O}pen {P}roblems.
\newblock {\em Bull. Amer. Math. Soc.}, 76:228--296, 1970.

\bibitem{kac}
Victor~G. Kac and Dale~H. Peterson.
\newblock Infinite-{D}imensional {L}ie {A}lgebras, {T}heta {F}unctions and
  {M}odular {F}orms.
\newblock {\em Adv. in Math.}, 53(2):125--264, 1984.

\bibitem{kul}
Vik.~S. Kulikov.
\newblock Degenerations of ${K}3$ {S}urfaces and {E}nriques {S}urfaces.
\newblock {\em Izv. Akad. Nauk SSSR Ser. Mat.}, 41(5):1008--1042, 1199, 1977.

\bibitem{looijenga}
Eduard Looijenga and Chris Peters.
\newblock Torelli {T}heorems for {K}\"ahler ${K}3$ {S}urfaces.
\newblock {\em Compositio Math.}, 42(2):145--186, 1980/81.

\bibitem{manin}
Yu.~I. Manin.
\newblock {\em Cubic {F}orms:{A}lgebra, {G}eometry, {A}rithmetic}.
\newblock North-Holland Publishing Co., Amsterdam, 1986.

\bibitem{mumford2}
David Mumford.
\newblock {\em Tata {L}ectures on {T}heta. {I}}.
\newblock Birkh\"auser Boston Inc., Boston, MA, 1983.
\newblock With the assistance of C. Musili, M. Nori, E. Previato and M.
  Stillman.

\bibitem{narain1}
K.~S. Narain.
\newblock New {H}eterotic {S}tring {T}heories in {U}ncompactified {D}imensions
  $<10$.
\newblock {\em Phys. Lett. B}, 169(1):41--46, 1986.

\bibitem{narain2}
K.~S. Narain, M.~H. Sarmadi, and E.~Witten.
\newblock A {N}ote on {T}oroidal {C}ompactification of {H}eterotic {S}tring
  {T}heory.
\newblock {\em Nuclear Phys. B}, 279(3-4):369--379, 1987.

\bibitem{pinkham}
Ulf Persson and Henry Pinkham.
\newblock Degeneration of {S}urfaces with {T}rivial {C}anonical {B}undle.
\newblock {\em Ann. of Math. (2)}, 113(1):45--66, 1981.

\bibitem{shapiro}
I.~I. Pjatecki{\u\i}-{\v{S}}apiro and I.~R. {\v{S}}afarevi{\v{c}}.
\newblock Torelli's {T}heorem for {A}lgebraic {S}urfaces of {T}ype ${\rm
  {k}}3$.
\newblock {\em Izv. Akad. Nauk SSSR Ser. Mat.}, 35:530--572, 1971.

\bibitem{schmid}
Wilfried Schmid.
\newblock Variation of {H}odge {S}tructure: {T}he {S}ingularities of the
  {P}eriod {M}apping.
\newblock {\em Invent. Math.}, 22:211--319, 1973.

\bibitem{sen}
Ashoke Sen.
\newblock {$F$}-{T}heory and {O}rientifolds.
\newblock {\em Nuclear Phys. B}, 475(3):562--578, 1996.

\bibitem{sig}
Carl~Ludwig Siegel.
\newblock {\em Advanced {A}nalytic {N}umber {T}heory}.
\newblock Tata Institute of Fundamental Research, Bombay, second edition, 1980.

\bibitem{siu}
Yum~Tong Siu and G{\"u}nther Trautmann.
\newblock Deformations of {C}oherent {A}nalytic {S}heaves with {C}ompact
  {S}upports.
\newblock {\em Mem. Amer. Math. Soc.}, 29(238):iii+155, 1981.

\bibitem{vafa}
Cumrun Vafa.
\newblock Evidence for ${F}$-theory.
\newblock {\em Nuclear Phys. B}, 469(3):403--415, 1996.

\bibitem{witten1}
Edward Witten.
\newblock World-{S}heet {C}orrections via {D}-instantons.
\newblock {\em J. High Energy Phys.}, (2):Paper 30, 18, 2000.

\end{thebibliography}
%\bibliographystyle{plain}

\vspace{.1in}
 
{\bf School of Mathematics, Institute for Advanced Study, Princeton, NJ 08540} \

Email address: {\tt {\bf clingher@ias.edu}} \
\vspace{.1in}

{\bf Department of Mathematics, Columbia University, New York, NY 10027} \

Email address: {\tt {\bf jm@math.columbia.edu}} 
\end{document}